\documentclass{amsbook}

\usepackage{times}
\usepackage[T1]{fontenc}
\usepackage[latin1]{inputenc}
\usepackage{xypic}
\xyoption{all}

\usepackage{url} 
\usepackage{graphicx}        % standard LaTeX graphics tool
\makeindex 
\makeatletter

 \theoremstyle{plain}
 \newtheorem{theorem}{Theorem}[chapter]
 \newtheorem{definition}[theorem]{Definition}
 \newtheorem{conjecture}[theorem]{Conjecture}

 \numberwithin{figure}{section}

 \theoremstyle{plain}
 \newtheorem{corollary}[theorem]{Corollary}
 \newtheorem{lemma}[theorem]{Lemma}
 \newtheorem{proposition}[theorem]{Proposition}

 \theoremstyle{remark}
 
 \newtheorem{remark}[theorem]{Remark}
 
\theoremstyle{remark}
 \newtheorem{example}[theorem]{Example}
 
 \newtheorem{problem}[theorem]{Problem}
 \newtheorem{question}[theorem]{Question}

\newcommand{\bN}{\mathbb{N}} \newcommand{\bP}{\mathbb{P}}
\newcommand{\sA}{\mathcal{A}} 
\newcommand{\sF}{\mathcal{F}} \newcommand{\sG}{\mathcal{G}}
 
\newcommand{\sD}{\mathcal{D}} 
\newcommand{\sP}{\mathcal{P}}

\newcommand{\sC}{\mathcal{C}}\newcommand{\sO}{\mathcal{O}}
\renewcommand{\O}{\mathcal{O}}
\renewcommand{\P}{\mathbb{P}} 
\newcommand{\C}{\mathbb{C}}\newcommand{\R}{\mathbb{R}} 
\newcommand{\Q}{\mathbb{Q}}

\DeclareMathOperator{\Aut}{Aut} 
\DeclareMathOperator{\bir}{bir}
 \DeclareMathOperator{\Chow}{Chow}
\DeclareMathOperator{\codim}{codim} \DeclareMathOperator{\Hilb}{Hilb}
\DeclareMathOperator{\Hom}{Hom} 

\DeclareMathOperator{\length}{length} \DeclareMathOperator{\loc}{loc}
\DeclareMathOperator{\mult}{mult}

\DeclareMathOperator{\Pic}{Pic} \DeclareMathOperator{\red}{red}
\DeclareMathOperator{\RatCurves}{RatCurves}
\DeclareMathOperator{\reg}{reg} \DeclareMathOperator{\Sing}{Sing}
\DeclareMathOperator{\Sym}{Sym} \DeclareMathOperator{\rank}{rank}

\DeclareMathOperator{\id}{id} \DeclareMathOperator{\Var}{Var}
\DeclareMathOperator{\coker}{coker} \DeclareMathOperator{\ad}{ad}
\DeclareMathOperator{\rat}{RatCurves^n}
\DeclareMathOperator{\vertical}{vert}\DeclareMathOperator{\Ch}{Ch}

\def\factor#1.#2.{\left. \raise 2pt\hbox{$#1$} \right/\hskip
  -2pt\raise -2pt\hbox{$#2$}}

\begin{document}

\title[Rational curves and applications]{Existence of rational curves
  on algebraic varieties, minimal rational tangents, and applications}

\author{Stefan Kebekus}
\author{Luis Sol\'a Conde}

\thanks{Both authors were supported in full or in part by the priority
  program ``Globale Methoden in der komplexen Geometrie'' of the
  Deutsche Forschungsgemeinschaft, DFG}

\address{Stefan Kebekus, Mathematisches Institut, Universit\"at zu
  K\"oln, Weyertal 86--90, 50931 K\"oln, Germany}
\email{stefan.kebekus@math.uni-koeln.de}
\urladdr{http://www.mi.uni-koeln.de/$\sim$kebekus}

\address{Luis Sol\'a Conde, Mathematisches Institut, Universit\"at zu K\"oln,
  Weyertal 86-90, 50931 K\"oln, Germany}
\email{lsola@math.uni-koeln.de}

\maketitle
\tableofcontents

\section*{Introduction}

Over the last two decades, the study of rational curves on algebraic
varieties has met with considerable interest. Starting with Mori's\index{Mori, Shigefumi}
landmark works \cite{Mori79, Mori82} it has become clear that many of
the varieties met daily by the algebraic geometer contain rational
curves, and that a variety at hand can often be studied by looking at
the rational curves it contains. Today, methods coming from the study
of rational curves on algebraic varieties are applied to a broad
spectrum of problems in higher-dimensional algebraic geometry, ranging
from uniqueness of complex contact structures \index{Contact manifold!uniqueness of contact structure} to deformation rigidity
of Hermitian symmetric manifolds\index{Deformation rigidity!of Hermitian symmetric spaces}.

In this survey we would like to give an overview of some of the recent
progress in the field, with emphasis on methods developed in and
around the DFG Schwerpunkt ``Globale Methoden in der komplexen
Geometrie''. Accordingly, there is a large body of important work that
we could not cover here. Among the most prominent results of the last
years is the breakthrough work of Graber, \index{Graber, Tom} Harris \index{Harris, Joe} and Starr, \index{Starr, Jason} \index{GHS Theorem, the base of the MRC fibration is not uniruled}
\cite{GHS03}, where it is shown that the base of the rationally
connected fibration is itself not covered by rational
curves\footnote{see \cite{Araujo05} for a good introduction}.  Another
result not touched in this survey is the recent progress toward the
abundance conjecture, by Boucksom\index{Boucksom, S\'ebastien}, Demailly\index{Demailly, Jean-Pierre}, P\u aun\index{Paun, Mihai@P\u aun, Mihai} and Peternell\index{Peternell, Thomas},
\cite{BDPP03}\index{BDPP, Theorem of Boucksom, Demailly, Peternell, P\u aun} ---see the article of Jahnke-Peternell-Radloff\index{Jahnke, Priska}\index{Radloff, Ivo} in this
volume instead.

\subsection*{Outline of the paper}

We start this survey in Chapter~\ref{chap:existence} by reviewing
criteria that can be used to show that a given variety is covered by
rational curves.  After mentioning Mori's results, we discuss foliated
varieties in some detail and present a recent criterion that contains
Miyaoka's fundamental\index{Miyaoka, Yoichi}\index{Miyaoka, Yoichi!criterion of uniruledness} characterization of uniruledness,
\cite{Miy85}, as a special case.  Its proof is rather elementary, and
a number of known results follow as simple corollaries.
Keel\index{Keel, Se\'an}-McKernan's\index{McKernan, James} work \cite{KMcK} on rational curves on
quasi-projective varieties\index{KMcK, rational curves on quasi-projective varieties} will be briefly discussed.

We will then, in Chapter~\ref{chap:geomtry} discuss the geometry of
(higher-dimensional) varieties that are covered by rational curves,
and present some ideas how minimal degree rational curves can be used
to study these spaces. The hero of these sections is the ``variety of
minimal rational tangents''\index{Variety of minimal rational tangents}\index{Variety of minimal rational tangents!short explanation}, or VMRT. In a nutshell, if $X$ is a
projective variety covered by rational curves, and $x \in X$ a general
point, then the VMRT is the subvariety of $\P(T_X|_x^\vee)$ which
contains the tangent directions to minimal degree rational curves that
pass through $x$. If $X$ is covered by lines, this is a very classical
object that has been studied by Cartan\index{Cartan, \'Elie} and Fubini\index{Fubini, Guido}\index{Cartan-Fubini type theorem!for varieties covered by lines} in the past. In the
general case, the VMRT is an important variety, similar perhaps to a
conformal structure, whose projective geometry as a subset of
$\P(T_X|_x^\vee)$ encodes much of the information on the underlying
space $X$ and determines $X$ to a large extent.

\medskip

In Chapters~\ref{chap:moduli}--\ref{chap:hyperbolicity} we will apply
the methods and results of Chapters~\ref{chap:existence} and
\ref{chap:geomtry} in three different settings. To start, we discuss
the moduli space of stable rank-two vector bundles on a curve\index{Moduli space of vector bundles on a curve} in
Chapter~\ref{chap:moduli}. There exists a classical construction of
rational curves on these spaces, the so-called ``Hecke-curves''\index{Hecke curves}. These
have been used to answer a large number of questions about moduli of
vector bundles. We name a few of the applications and sketch a proof
for a result that helps to give an upper bound for the multiplicities of divisors at
general points of the moduli space; the result bears perhaps some
resemblance with the classical Riemann singularity theorem\index{Riemann singularity theorem}.

The existence results of Chapter~\ref{chap:existence} can also be used
to study varieties for which it is known \emph{a priori} that they are
not covered by rational curves. We conclude this paper by giving two
examples in Chapters~\ref{chap:hkp} and \ref{chap:hyperbolicity}.

Chapter~\ref{chap:hkp} deals with deformations of surjective morphisms\index{Deformation of a surjective morphism}
$f: X \to Y$, where the target is not covered by rational curves, or
at least not rationally connected. It turns out that there exists a
natural refinement of Stein factorization\index{Refinement of Stein factorization} which factors $f$ via an
intermediate variety $Z$, and that the existence results of
Chapter~\ref{chap:existence} can be used to show that the associated
component of the deformation space $\Hom(X,Y)$ is essentially the
automorphism group of $Z$.

In Chapter~\ref{chap:hyperbolicity} we apply the results of
Chapter~\ref{chap:existence} to the study of families of canonically
polarized varieties\index{Family of canonically polarized varieties}. Generalizing Shafarevich's\index{Shafarevich, Igor R.} hyperbolicity
conjecture\index{Shafarevich, Igor R.!Hyperbolicity conjecture}, it has been conjectured by Viehweg\index{Viehweg, Eckart}\index{Viehweg, Eckart!conjecture relating variation and log Kodaira dimension} that the base of a
smooth family of canonically polarized varieties is of log general
type if the family is of maximal variation\index{Family of canonically polarized varieties!variation of}. Using Keel-McKernan's\index{Keel, Se\'an}\index{McKernan, James}
existence results for rational curves on quasi-projective varieties,
we relate the \emph{variation} of a family to the \emph{logarithmic
  Kodaira dimension}\index{Logarithmic Kodaira dimension} of the base and sketch a proof for an affirmative
answer to Viehweg's\index{Viehweg, Eckart} conjecture for families over surfaces.
 
\medskip

Unless explicitly mentioned, we always work over the complex number
field.

\subsection*{Other references}

The reader who is interested in a broader perspective will obviously
want to consult the standard reference books \cite{K96} and
\cite{Debarre01}. Hwang's\index{Hwang, Jun-Muk} important survey \cite{Hwa00} explains by
way of examples how rational curves can be employed to study Fano
manifolds of Picard number one.  The article \cite{AK03} contains an
excellent introduction to the deformation theory of rational curves
and rational connectivity, also in the more general setting of
varieties defined over non-closed fields.

% Local Variables:
% TeX-master: "RC-arXive.tex"
% End:

\part{Existence of rational curves}

\chapter{Existence of rational curves}
\label{chap:existence}

The modern interest in rational curves on projective varieties started
with Shigefumi~Mori's \index{Mori, Shigefumi} fundamental work \cite{Mori79}. In his proof of
the Hartshorne-Frankel conjecture \index{Hartshorne-Frankel conjecture}\index{Hartshorne, Robin}\index{Frankel, Theodore}, he devised a new method to prove
the existence of rational curves on manifolds whose tangent bundle
possess certain positivity properties\index{Tangent bundle!positivity}. Though not explicitly
formulated like this, the following results appear in his papers.

\begin{theorem}[\protect{\cite{Mori79, Mori82}}]\label{thm:Mori1}\index{Rational curve!existence on Fano manifolds}
  Let $X$ be a complex-projective manifold. If $X$ is Fano, i.e., if
  $-K_X$ is ample, then $X$ is uniruled, i.e., covered by rational
  curves. More precisely, if $x \in X$ is any point, then there exists
  a rational curve $\ell \subset X$, such that $x \in \ell$, and
  $-K_X\cdot\ell \leq \dim X+1$. \qed
\end{theorem}
\begin{theorem}[\protect{\cite{Mori82}}]\label{thm:Mori2}\index{Rational curve!existence if $K_X$ is not nef}
  Let $X$ be a complex-projective manifold. If $K_X$ is not nef, then
  $X$ contains rational curves. More precisely, if $C \subset X$ is a
  curve with $K_X \cdot C < 0$, and $x \in C$ any point, then there
  exists a rational curve $\ell \subset X$ that contains $x$. \qed
\end{theorem}

We refer to \cite{CKM88} or \cite{Debarre01} for an accessible
introduction. These results, and the subsequently developed ``minimal
model program'' allowed, in dimension three, to give a positive answer
to a long-standing conjecture attributed to Mumford \index{Mumford, David} that characterizes
manifolds covered by rational curves as those without pluricanonical
forms\index{Mumford, David!conjecture relating uniruledness and Kodaira dimension}.

\begin{conjecture}\label{conj:mumford}
  A projective manifold $X$ is covered by rational curves if and
  only if $\kappa(X) = -\infty$.
\end{conjecture}

In higher dimensions, the conjecture is still open, although the
recent result of Boucksom\index{Boucksom, S\'ebastien}, Demailly\index{Demailly, Jean-Pierre}, Peternell\index{Peternell, Thomas} and P\u aun\index{Paun, Mihai@P\u aun, Mihai}\index{BDPP, Theorem of Boucksom, Demailly, Peternell, P\u aun} is considered
a serious step forward.

The first Chern class of $T_X$ used in
Theorems~\ref{thm:Mori1} and \ref{thm:Mori2} is, however, a rather
coarse measure of positivity. For instance, if $Y$ is a Fano manifold
and $Z$ a torus, then the tangent bundle of the product $X = Y \times
Z$ splits into a direct sum $T_X = \sF \oplus \sG$ where $\sF =
p_1^*(T_Y)$ has positivity properties and identifies tangent
directions to the rational curves contained in $X$. In general, rather
than looking at $K_X$, one would often like to deduce the existence of
rational curves from positivity properties of subsheaves $\sF \subset
T_X$ and relate the geometry of those curves to that of $\sF$. The
most fundamental result in this direction is Miyaoka's criterion of
uniruledness\index{Miyaoka, Yoichi!criterion of uniruledness}\index{Miyaoka, Yoichi}. In order to state Miyaoka's result\index{Miyaoka, Yoichi}, we recall the
theorem of Mehta-Ramanathan for normal varieties.

\begin{theorem}[Mehta-Ramanathan, Flenner, \protect{\cite{MR82, Flenner84}}]\label{thm:fmr}\index{Mehta-Ramanathan theorem}\index{Flenner, Hubert}\index{Mehta, Vikram B.}\index{Ramanathan, Annamalai}
  Let $X$ be a normal variety of dimension $n$, and $L_1, \ldots
  L_{n-1} \in \Pic(X)$ be ample line bundles. If $m_1, \ldots, m_{n-1}
  \in \mathbb N$ are large enough and $H_i$ are general elements of
  the linear systems $|m_i \cdot L_i|$, then the curve
  $$
  C = H_1 \cap \cdots \cap H_{n-1}
  $$
  is smooth, reduced and irreducible, and the restriction of the
  Harder-Narasimhan filtration of $T_X$ to $C$ is the Harder-Narasimhan
  filtration of $T_X|_C$. \qed
\end{theorem}

We refer to \cite{Langer04a, Langer04b} for a discussion and an
explicit bound for the $m_i$. A detailed account of slope,
semistability and of the Harder-Narasimhan filtration of vector
bundles on curves is found in \cite{Seshadri82}.

\begin{definition}\label{def:GCIC}
  We call a curve $C \subset X$ as in Theorem~\ref{thm:fmr} a
  \emph{general complete intersection curve}\index{General complete intersection curve} in the sense of
  Mehta-Ramanathan\index{Mehta, Vikram B.}\index{Ramanathan, Annamalai}.
\end{definition}

Miyaoka's\index{Miyaoka, Yoichi}\index{Miyaoka, Yoichi!criterion of uniruledness} result then goes as follows.

\begin{theorem}[\protect{\cite[thm.~8.5]{Miy85}}]\label{thm:Miyaoka1}
  Let $X$ be a normal projective variety, and $C \subset X$ a general
  complete intersection curve\index{General complete intersection curve} in the sense of Mehta-Ramanathan. Then
  $\Omega^1_X|_C$ is a semi-positive \index{Semipositive vector bundle}vector bundle unless $X$ is
  uniruled. \qed
\end{theorem}

We refer the reader to \cite{MP97} for a detailed overview of
Miyaoka's\index{Miyaoka, Yoichi} theory of foliations in positive characteristic\index{Foliation!in positive characteristic}. The
relation between negative directions of $\Omega^1_X|_C$ and tangents
to rational curves has been studied in
\cite[sect.~9]{SecondAsterisque}. We give a full account in
Section~\ref{sec:foliations}.

\subsection*{Outline of the section}

In Section~\ref{sec:RCfoliations}, we study the case where $X$ is a
complex manifold and $\sF \subset T_X$ is a (possibly singular)
foliation\index{Foliation}. The main result ---which appeared first in the preprint
\cite{BMcQ01} of Bogomolov\index{Bogomolov, Fedor A.} and McQuillan\index{McQuillan, Michael L.}--- gives a criterion to
guarantee that the leaves of $\sF$ are compact\index{Leaf of a foliation!algebraicity} and rationally
connected\index{Leaf of a foliation!rational connectivity}. Miyaoka's\index{Miyaoka, Yoichi} characterization of uniruledness\index{Miyaoka, Yoichi!criterion of uniruledness},
Theorem~\ref{thm:Miyaoka1}, and the statements of
\cite[sect.~9]{SecondAsterisque} follow as immediate corollaries.
Apart from a simple vanishing theorem for vector bundles in positive
characteristic\index{Vanishing result in characteristic $p$}, the proof employs only standard techniques of Mori
theory that are well discussed in the literature. In particular, it
will not be necessary to make any reference to the more involved
properties of foliations in characteristic $p$.  We also mention a
sufficient condition to ensure that all leaves of a given foliation
are algebraic\index{Leaf of a foliation!algebraicity}.

In Sections~\ref{sec:foliations} and \ref{sec:partRCQ} we discuss the
relation with the rationally connected quotient\index{Rationally connected quotient}. The results of
Section~\ref{sec:RCfoliations} are applied to show that $\mathbb
Q$-Fano varieties with unstable tangent bundles always admit a
sequence of partial rational quotients naturally associated to the
Harder-Narasimhan \index{Harder-Narasimhan filtration!of the tangent bundle} filtration of the tangent bundle.

We will later need to discuss an analog of Mumford\index{Mumford, David}'s
Conjecture~\ref{conj:mumford} \index{Mumford, David!conjecture relating uniruledness and Kodaira dimension}for quasi-projective varieties. Here the
logarithmic Kodaira dimension \index{Logarithmic Kodaira dimension}takes the role of the regular Kodaira
dimension, and rational curves are replaced by $\mathbb C$ or $\mathbb
C^*$. For surfaces, this setup has been studied by Keel\index{Keel, Se\'an} and McKernan\index{McKernan, James}.
We recall their results in Section~\ref{sec:KMcK}.

\section{Rationally connected foliations}
\label{sec:RCfoliations}

In the previous section we have mentioned that positivity properties
of $T_X$ imply the existence of rational curves in $X$. Here, we will
study how Mori's\index{Mori, Shigefumi} ideas can be applied to foliations on complex
varieties. We recall the notion of rational connectivity and fix
notation first.
\begin{figure}[tbp]
  \centering
  \begin{picture}(5,1.7)(0,0)
    \put(0, 0){\includegraphics[height=1.5cm]{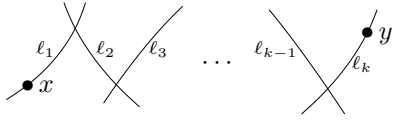}}
    \put(0.2,0.3){$\bullet \, x$}
    \put(2.6,0.6){$\cdots$}
    \put(4.72,1.0){$\bullet \, y$}
    \put(0.4,0.8){$\scriptstyle \ell_1$}
    \put(1.2,0.8){$\scriptstyle \ell_2$}
    \put(1.9,0.8){$\scriptstyle \ell_3$}
    \put(3.3,0.8){$\scriptstyle \ell_{k-1}$}
    \put(4.6,0.6){$\scriptstyle \ell_k$}
  \end{picture}
  \caption{The points $x$ and $y$ are joined by a chain of rational curves of length $k$}
  \label{fig:rcc}
\end{figure}
\begin{definition}\label{def:ratconn}
  A normal variety $X$ is \emph{rationally chain connected}\index{Rational chain connectivity} if any two
  general points $x, y \in X$ can be joined by a chain of rational
  curves, as shown in Figure~\ref{fig:rcc}. The variety $X$ is
  \emph{rationally connected}\index{Rational connectivity} if for any two general points $x, y \in
  X$ there exists a single rational curve that contains both.
\end{definition}
\begin{remark}
  For smooth varieties, the two notions \emph{rationally chain
    connected} and \emph{rationally connected} agree, see
  \cite[sect.~4.7]{Debarre01}.
\end{remark}

\begin{definition}
  In this survey, a \emph{foliation}\index{Foliation} $\sF$ on a normal variety $X$ is
  a saturated, integrable subsheaf of $T_X$. A
  \emph{leaf}\index{Leaf of a foliation} of $\sF$ is a maximal $\sF$-invariant connected subset
  of the set $X^\circ$ where both $X$ and $\sF$ are regular. A leaf is
  called \emph{algebraic}\index{Leaf of a foliation!algebraicity} if it is open in its Zariski closure.
\end{definition}

The main result of this section asserts that positivity properties of
$\sF$ imply algebraicity and rational connectivity of the leaves\index{Leaf of a foliation!rational connectivity}. In
particular, it gives a criterion for a manifold to be covered by
rational curves.

\begin{theorem}\label{thm:BMcQ}
  Let $X$ be a normal complex projective variety, $C \subset X$ a
  complete curve which is entirely contained in the smooth locus
  $X_{\reg}$, and $\sF \subset T_X$ a (possibly singular) foliation
  which is regular along $C$. Assume that the restriction $\sF|_C$ is
  an ample vector bundle on $C$. If $x \in C$ is any point, the leaf
  through $x$ is algebraic. If $x \in C$ is general, the closure of
  the leaf is rationally connected.\index{Leaf of a foliation!algebraicity} \index{Leaf of a foliation!rational connectivity}
\end{theorem}
 
The statement appeared first in the preprint \cite{BMcQ01} by
Bogomolov\index{Bogomolov, Fedor A.} and McQuillan\index{McQuillan, Michael L.}. Below we sketch a simple proof which
recently appeared in \cite{KST05}.

\begin{remark}
  In Theorem~\ref{thm:BMcQ}, if $x \in C$ is any point, it is
  \emph{not} generally true that the closure of the leaf through $x$
  is rationally connected ---this was wrongly claimed in \cite{BMcQ01}
  and in the first preprint versions of \cite{KST05}.\index{Leaf of a foliation!rational connectivity}
\end{remark}

The classical Reeb stability theorem \index{Reeb stability theorem}for foliations
\cite[thm.~IV.3]{CLN85}, the fact that rationally connected manifolds
are simply connected\index{Rational connectivity!implies simple connectivity} \cite[cor.~4.18]{Debarre01}, and the openness of
rational connectivity \index{Rational connectivity!openness of}\cite[cor.~2.4]{KMM92} immediately yield the
following\footnote{H\"oring has independently obtained similar
  results, \cite{Hor05}.}\index{Horing, Andreas@H\"oring, Andreas}.

\begin{theorem}[\protect{\cite[thm.~2]{KST05}}]
  In the setup of Theorem~\ref{thm:BMcQ}, if $\sF$ is regular, then
  all leaves are rationally connected submanifolds.\index{Leaf of a foliation!rational connectivity} \qed
\end{theorem}

In fact, a stronger statement holds that guarantees that most leaves
are algebraic and rationally connected if there exists a single leaf
through $C$ whose closure does not intersect the singular locus of
$\sF$, see \cite[thm.~28]{KST05}.

The following characterization of rational connectivity is a
straightforward corollary of Theorem~\ref{thm:BMcQ}.

\begin{corollary}\label{cor:charratconn}
  Let $X$ be a complex projective variety and let $f:C\rightarrow X$
  be a curve whose image is contained in the smooth locus of $X$ and
  such that $T_X|_C$ is ample. Then $X$ is rationally connected.\index{Rational connectivity!criterion}
  \qed
\end{corollary}

\subsection{Preparation for the proof of Theorem~\ref{thm:BMcQ}: The rational quotient map}

Before sketching a proof of Theorem~\ref{thm:BMcQ}, we recall a few
facts and notations associated with the rationally connected quotient\index{Rationally connected quotient}
of a normal variety, introduced by Campana\index{Campana, Fr\'ed\'eric} and Koll\'ar-Miyaoka-Mori\index{Koll\'ar, J\'anos}\index{Miyaoka, Yoichi}\index{Mori, Shigefumi},
\cite{Campana92, KMM92b}. See \cite{Debarre01}, \cite{K96} or
\cite{Araujo05} for an introduction.

\begin{theorem}[Campana, and Koll\'ar-Miyaoka-Mori]\label{thm:rq}
  Let $X$ be a smooth projective variety. Then there exists a rational
  map $q: X \dasharrow Q$, with the following properties.
  \begin{enumerate}
  \item The map $q$ is \emph{almost
      holomorphic}\index{Almost holomorphic map}\footnote{\cite{Debarre01} and \cite{K96} use the
      word ``fibration'' for an almost holomorphic map.}, i.e., there
    exists an open set $X^\circ \subset X$ such that $q|_{X^\circ}$ is
    a proper morphism.
  \item If $X_F \subset X$ is a general fiber of $q$, then $X_F$ is a
    compact, rationally chain connected manifold, i.e., any two points
    in $X_F$ can be joined by a chain of rational curves.
  \item If $X_F \subset X$ is a general fiber of $q$ and $x \in X_F$ a
    general point, then any rational curve that contains $x$ is
    automatically contained in $X_F$.
  \end{enumerate}
  Properties (1)--(3) define $q$ uniquely up to birational
  equivalence. \qed
\end{theorem}

\begin{definition}\label{def:ratQuot}
  Any map $q : X \dasharrow Q$ for which properties (1)--(3) of
  Theorem~\ref{thm:rq} hold is called a \emph{maximally rationally
    chain connected fibration}\index{Maximally rationally chain connected fibration} of $X$, or MRCC fibration, for short.
  
  Let $X$ be a normal projective variety, and let $q : X \dasharrow Q$
  be the rational map defined through the maximally rationally chain
  connected fibration of a desingularization of $X$. We call $q$ the
  \emph{rationally connected quotient}\index{Rationally connected quotient} of $X$.
\end{definition}

\begin{remark}
  If $X$ is normal, the rationally connected quotient is again defined
  uniquely up to birational equivalence.
\end{remark}
\begin{remark}
  Rational chain connectivity is \emph{not} a birational invariant.\index{Rational chain connectivity!is not a birational invariant}
  For instance, a cone over an elliptic curve is rationally chain
  connected, while the ruled surface obtained by blowing up the vertex
  is not. It is therefore important in Definition~\ref{def:ratQuot} to
  pass to a desingularization.
\end{remark}

In the proof of Theorem~\ref{thm:BMcQ}, we will need two important
properties of the rationally connected quotient.

\begin{theorem}[Universal property of the maximally  rationally chain connected fibration, \protect{\cite[thm.~IV.5.5]{K96}}]\label{thm:mrcuniv}
  Let $X_1$, $X_2$ be projective manifolds, and $f_X : X_1 \dasharrow
  X_2$ dominant. If $q_i: X_i \dasharrow Q_i$ are the maximally
  rationally chain connected fibrations, then there exists a
  commutative diagram as follows.\index{Maximally rationally chain connected fibration!universal property}
  $$
  \xymatrix{
    X_1 \ar@{-->}[r]^{f_X} \ar@{-->}[d]_{q_1} & X_2 \ar@{-->}[d]^{q_2} \\
    Q_1 \ar@{-->}[r]_{f_Q} & Q_2
  }
  $$
  \qed
\end{theorem}

\begin{theorem}[Graber-Harris-Starr, \protect{\cite{GHS03}}]\label{thm:ghs}\index{GHS Theorem, the base of the MRC fibration is not uniruled}\index{Graber, Tom}\index{Harris, Joe}\index{Starr, Jason}
  If $X$ is a normal projective variety and $q : X \dasharrow Q$ the
  rationally connected quotient, then $Q$ is not covered by rational
  curves. \qed
\end{theorem}

\subsection{Sketch of proof of Theorem~\ref{thm:BMcQ}}

\subsection*{Step 1: Reduction to the case where $C$ is transversal to $\sF$}

\begin{figure}[tbp]
  \centering
     \begin{picture}(11.5,7)(0,0)
      \put(0, 3){\includegraphics[height=4cm]{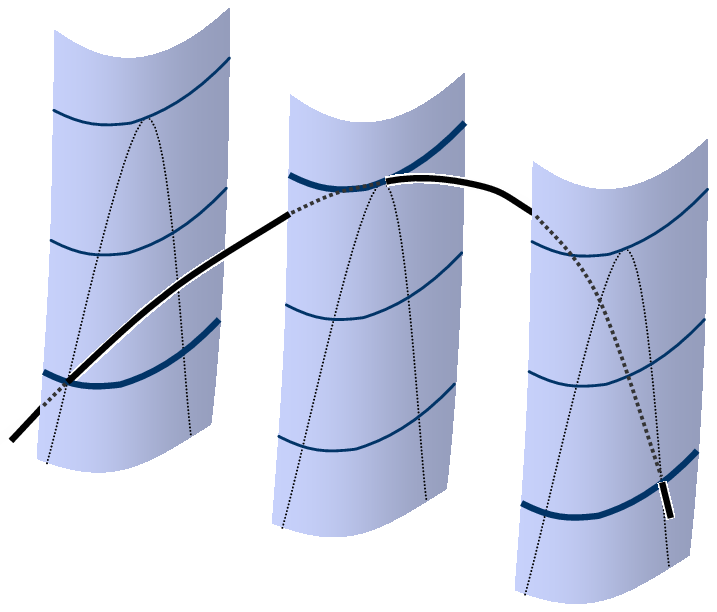}}
      \put(2,1.5){\vector(0,1){1.5}}
      \put(0.7,2.2){${\scriptstyle \sigma :=(\nu,id)}$}
      \put(2.5,3){\vector(0,-1){1.5}}
      \put(2.7,2.2){${\scriptstyle p_2=:\pi}$}
      \put(5.5,5){\vector(1,0){2}}
      \put(6.2,5.2){${\scriptstyle p_1}$}
      \put(8, 3){\includegraphics[height=3cm]{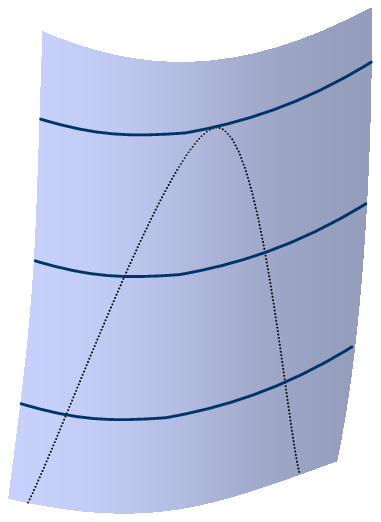}}     
      \put(1.5,0.8){$\tilde{C}$}
      \put(2.5,6.7){$Y=X\times\tilde{C}$}
      \put(9,5.9){$X$}
      \put(9.45,4.85){${\scriptstyle C}$}
      \put(1.6,5.55){${\scriptstyle C'}$}
      \put(10.5,4.7){\tiny leaves}
      \put(10.4,4.95){\vector(-1,2){0.3}}
      \put(10.4,4.75){\vector(-1,0){0.3}}
      \put(10.4,4.55){\vector(-1,-2){0.3}}
       \put(0.2,0){\includegraphics[width=4cm]{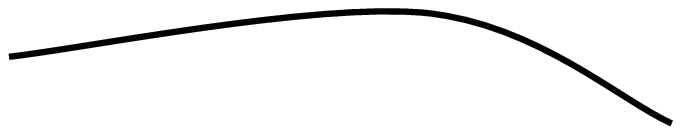}}  
     \end{picture}
   \caption{Reduction to the case of a normal foliation}
   \label{fig:normFol}
 \end{figure}

Following an idea of Bogomolov\index{Bogomolov, Fedor A.} and McQuillan\index{McQuillan, Michael L.}, we consider a
non-constant morphism $\nu:\tilde{C} \to C$ from a smooth curve
$\tilde{C}$ of positive genus $g(\tilde{C}) > 0$. Let $Y$ denote the
product $X\times \tilde{C}$ with projections $p_1$ and $p_2$ and
consider the following diagram, depicted in Figure~\ref{fig:normFol}.
$$
\xymatrix{
  Y \ar[r]^{p_1} \ar[d]^{p_2}& X \\
  \tilde C \ar@/^0.3cm/[u]^{\sigma = (\nu,id)}
}
$$
It is obviously enough to show Theorem~\ref{thm:BMcQ} for the
variety $Y$, the curve $C':= \sigma(\tilde{C})$ and the foliation
$\sF_Y: = p_1^*(\sF) \subset T_{Y|\tilde C} \subset T_Y$, which is
ample along $C'$. The advantage lies in the smoothness of $C'$, and in
the transversality of $\sF_Y$ to $C'$. Both properties are required
for the next step.

\subsection*{Step 2: Algebraicity of the leaves}\index{Leaf of a foliation!algebraicity}

Since $C'$ is everywhere transversal to $\sF_Y$, we can apply the
classical Frobenius theorem: there exists an analytic submanifold
$W\subset Y$ which contains $C'$ and has the property that its
fibers over $\tilde C$ are analytic open sets of the leaves of
$\sF_Y$. Let $\overline{W}$ be the Zariski closure of $W$.

Using the transversality of $C'$ and $\sF_Y$ and the fact that
$\sF_Y|_C$ is ample, a theorem of Hartshorne\index{Hartshorne, Robin}, \cite[thm.~6.7]{Ha68},
asserts that $\dim \overline{W} = \dim W$. Accordingly, every leaf is
algebraic\index{Leaf of a foliation!algebraicity}. Step 2 again follows \cite{BMcQ01}, but see also
\cite[thm.~3.5]{Bost01}.

\subsection*{Step 3: Setup of notation}

Replacing $X$ by a desingularization of the normalization of
$\overline{W}$, we are reduced to prove Theorem~\ref{thm:BMcQ} under
the following extra assumptions: $X$ is smooth, $C$ is a smooth curve
of genus $g(C)>0$, and there exists a morphism $\pi : X \to C$ such
that
\begin{itemize}
\item $\pi$ has connected fibers and is smooth along $C$,
\item the foliation $\sF$ is the foliation associated to $\pi$,
  i.e.~$\sF = T_{X|C}$ wherever $\pi$ is smooth, and 
\item $\pi$ admits a section, $\sigma : C \to X$.
\end{itemize}

We will have to show that the general $\pi$-fiber is rationally \index{Leaf of a foliation!rational connectivity}
connected. To this end, consider the rationally connected quotient
$q:X\dashrightarrow Z$ of $X$\index{Rationally connected quotient}. The universal property of the maximally rationally
chain connected fibration\index{Maximally rationally chain connected fibration!universal property}, Theorem~\ref{thm:mrcuniv}, then yields a
diagram as follows.
\begin{equation}\label{eq:extdiag}
  \xymatrix{
    X \ar[d]^{\pi} \ar@{-->}[r]^{q} & Z \ar@{-->}@/^0.3cm/[ld]^{\beta} \\
    C \ar@/^0.3cm/[u]^{\sigma}}
\end{equation}
We finally fix a very ample line bundle $H_Z$ on $Z$, and we
denote by $H_X$ its pull-back to $X$. We can then consider $q$ as the
rational map associated to a certain linear subsystem of $H^0(X,
H_X)$.

Observe that to prove Theorem~\ref{thm:BMcQ}, it suffices to show that
$\dim Z=1$. Namely, if $\dim Z = 1$, then $\pi$ is itself a rationally
connected quotient, and its general fiber will therefore be rationally
chain connected, hence rationally connected.

We assume the contrary and suppose that $\dim Z>1$. Below we will then
show the following.

\begin{proposition}\label{prop:pfBMcQ} 
  Assume that $\dim(Z)\geq 2$.  Then $Z$ is uniruled with curves of
  $H_Z$-degrees at most $d:=2\deg\sigma^*(H_X)\cdot\dim X$.
\end{proposition}

This clearly contradicts Theorem~\ref{thm:ghs}, and concludes the
proof of Theorem~\ref{thm:BMcQ}.\index{GHS Theorem, the base of the MRC fibration is not uniruled}

\subsection*{Step 4: Strategy of proof for Proposition~\ref{prop:pfBMcQ} and Theorem~\ref{thm:BMcQ}}

Assume for the moment that the morphism $\sigma$ admits a large number
of deformations. More precisely, assume that there exists a component
$\mathcal H \subset \Hom(C,X)$ that contains $\sigma$, and an open set
$\Omega\subset\mathcal H$ such that the following holds:
\begin{enumerate}
\item If $T \subset X$ is any set of codimension $\leq 2$, then the
  set of morphisms that avoid $T$, $\mathcal A :=\{\tau \in \Omega
  \,|\, \tau^{-1}(T) = \emptyset\}$, is a non-empty open set in
  $\Omega$.
\item If $\tau \in \Omega$ and $x \subset C$ a point, then the
  evaluation morphism associated to the set $\Omega_{\tau(x)}=\{\tau'
  \in \Omega \,|\, \tau'(x) = \tau(x)\}$ still dominates $X$.
\end{enumerate}

If $(1)$ and $(2)$ hold true, choosing $T$ to be the indeterminacy
locus of $q$ the natural morphism $\mathcal A\to
\Hom(C,Z)$ provides a family $\mathcal C \subset
\Hom(C,Z)$ verifying:
\begin{itemize}
\item The evaluation morphism associated to $\sC$ dominates $Z$, and
\item if $x \in C$ and $\tau \in \mathcal C$ are general points, then
  the evaluation morphism associated to the set $\mathcal C_x =\{\tau'
  \in \mathcal C \,|\, \tau'(x) = \tau'(x)\}$ still dominates $Z$.
\end{itemize}

If $\dim Z \geq 2$, we can then use Mori's\index{Mori, Shigefumi} Bend-and-Break argument.\index{Bend-and-Break}

\begin{definition}
  Let $f : C \to Z$ be any morphism and $B \subset C$ a subscheme of
  finite length. The space of morphisms that agree with $f$ on $B$ is
  denoted as $\Hom(C,Z,f|_B)$.
\end{definition}

\begin{theorem}[Mori's Bend-and-Break, \protect{\cite{MM86}, \cite[prop.~3.3]{Kol91}}]\label{thm:bab2}\index{Mori, Shigefumi}\index{Bend-and-Break}
  Let $Z$ be a projective variety and let $H_Z$ be a nef $\R$-divisor
  on $Z$. Let $C$ be a smooth, projective and irreducible curve,
  $B\subset Z$ a finite subscheme and $f:C\rightarrow Z$ a
  non-constant morphism. Assume that $Z$ is smooth along $f(C)$. If
  $\dim_{[f]}\Hom(C,Z,f|_B)\geq 1$ then there exists a rational curve
  $R$ on $Z$ meeting $f(B)$ and such that
  $$
  H_Z\cdot R\leq\frac{2H_Z\cdot C}{\# B}. 
  $$
  If $H_Z$ is ample, the above inequality can be made strict. \qed
\end{theorem}
\begin{remark}
  Theorem~\ref{thm:bab2} works for varieties defined over
  algebraically closed fields of arbitrary characteristic.
\end{remark}
In our situation, Theorem~\ref{thm:bab2} with $B=(*)$ and $f$ a
general element of $\mathcal C$ would show that $Z$ is uniruled and
hence prove Proposition~\ref{prop:pfBMcQ} by contradiction with
Theorem~\ref{thm:ghs} ---for this, a simpler version of Mori's
Bend-and-Break\index{Bend-and-Break} would suffice, but we will need the full force of
Theorem~\ref{thm:bab2} later.

Unfortunately, the ampleness of the normal bundle of $C$ in $X$ does
not guarantee the existence of large family of deformations of
$\sigma$ that satisfy conditions (1) and (2). We can circumvent this
problem using reduction modulo $p$.

\subsection*{Step 5: Reduction modulo $p$}
 
In view of Mori's\index{Mori, Shigefumi} standard argument using reduction modulo $p$\index{Mori, Shigefumi!Reduction argument modulo $p$}, see
\cite{Mori79}, \cite{Debarre01} or \cite{CKM88}, it is enough to prove
Proposition~\ref{prop:pfBMcQ} over an algebraically closed field $k$
of characteristic $p$, for $p$ large enough. We use a subindex $k$ to
denote the reductions modulo $p$ of all the objects defined above, and
let $F$ denote the $k$-linear Frobenius morphism.

\subsection*{Step 6: A vanishing result in characteristic $p$}

We briefly recall the language of $\Q$-twisted vector bundles, as
explained in \cite[II,~6.2]{Laz04}. This notion is a generalization of
the concept of a $\Q$-divisor to higher rank. It allows us to make a
finer use of the positivity of a vector bundle.

We identify rational numbers $\delta$ with numerical classes
$\delta\cdot[P]\in N^1_{\Q}(C)$, where $P$ is a point in $C$. For
every $\delta\in\Q$, the {\it $\Q$-twist} $E\langle\delta\rangle$ is
defined as the ordered pair of $E$ and $\delta$. A $\Q$-twisted vector
bundle is said to be {\it ample} if the class $c_1 \bigl(\O_{\mathbb
  P_C(E)}(1) \bigr)+\pi^*(\delta)$ is ample on the projectivized
bundle $\P_C(E)$, where $\pi$ denotes the natural projection. One
defines the degree
$\deg(E\langle\delta\rangle):=\deg(E)+\rank(E)\delta$. A quotient of
$E\langle\delta\rangle$ is a $\Q$-twisted vector bundle of the form
$E'\langle\delta\rangle$ where $E'$ is a quotient of $E$.  Pull-backs
of $\Q$-twisted vector bundles are defined in the obvious way. 

We can now formulate the following vanishing result in characteristic
$p$ \index{Vanishing result in characteristic $p$}that will be used later on.

\begin{proposition}[\protect{\cite[prop.~9]{KST05}}]\label{prop:vanishp} 
  Let $C_k$ be a curve defined over an algebraically closed field of
  characteristic $p>0$.  Let $E_k$ be a vector bundle of rank $r$ over
  $C_k$, and $\delta$ a positive rational number. Assume that
  $E_k\langle -\delta\rangle$ is ample and that the ``vanishing
  threshold'' \index{Vanishing threshold}
  $$
  b_p(\delta):=p\delta-2g(C)+1
  $$
  is non-negative. Let $F: C_k[1] \to C_k$ be the $k$-linear
  Frobenius morphism. Then for every subscheme $B\subset C_k[1]$ of
  length smaller than or equal to $b_p(\delta)$ we have
  $$
  H^1 \bigl(C_k[1], F^*(E_k)\otimes I_B \bigr) = \{ 0\}.
  $$
  Further, $F^*(E_k)\otimes I_B$ is globally generated.
\end{proposition}

\begin{proof}
  To prove both vanishing and global generation, it is enough to show
  that
  $$
  H^1 \bigl(C_k[1], F^*(E_k)\otimes I_B \bigr) = \{ 0\} \,
  \mbox{ for } \#(B)\leq b_p(\delta)+1.
  $$
  Since $F$ is finite, the pull-back $F^*(E_k\langle
  -\delta\rangle)=F^*(E_k)\langle -p\delta\rangle$ is also ample, and
  so is every quotient of rank one. In particular,
  $$
  \Hom_{C_k} \bigl( F^*(E_k),\, \O_{C_k[1]}(B) \otimes
  \omega_{C_k[1]}\bigr) = \{0\} \, \mbox{ if
  }\deg(\O_{C[1]}(B)\otimes\omega_{C[1]})\leq p\delta.
  $$
  The proof is concluded by applying Serre duality.
\end{proof}

\subsection*{Step 7: Proof in characteristic $p$}

Back to the proof of Theorem~\ref{thm:BMcQ}. We would like to apply
the vanishing result of Proposition~\ref{prop:vanishp}\index{Vanishing result in characteristic $p$} to describe the
space of relative deformations of 
$$
\tau := \sigma_k \circ F : C_k[1] \to X_k.
$$
over $C_k$. To this end, recall the following standard description
of the Hom-scheme.

\begin{theorem}\label{thm:descrHom}
  Let $H:=\Hom_{C_k}(C_k[1],X_k)$ be the space of relative
  deformations of $\tau$, and let $\nu \in H$ be any element.
  \begin{enumerate}
  \item If $H^1(C_k[1], \nu^*(T_{X_k|C_k})) = 0$, then $H$ is smooth
    at $\nu$, and has dimension $H^0(C_k[1], \nu^*(T_{X_k|C_k}))$.
    
  \item Let $T \subset X_k$ be any set of codimension $\leq 2$. If
    $\tau^*(T_{X_k|C_k})$ is globally generated and $H^1(C_k[1],
    \tau^*(T_{X_k|C_k})) = 0$, then the set of morphisms whose images
    avoid $T$ is a non-empty open subset of $H$.
    
  \item If $B \subset C_k$ is any subscheme of finite length, if
    $\nu^*(T_{X_k|C_k}) \otimes I_B$ is globally generated and
    $H^1(C_k[1], \nu^*(T_{X_k|C_k})\otimes I_B) = 0$, then the images
    of morphisms $\nu'$ that agree with $\nu$ along $B$ dominate $X$.
  \end{enumerate}
  \qed
\end{theorem}

By Hartshorne's \index{Hartshorne, Robin} characterization of ampleness\index{Hartshorne, Robin!characterization of ampleness},
\cite[Thm.~6.4.15]{Laz04}, the ampleness of $\sigma^*(T_{X|C})$ is
equivalent to the ampleness of the $\Q$-twisted vector bundle
$\sigma^*(T_{X|C})\langle -1/\dim X\rangle$. Note also that the
ampleness of $\sigma^*(T_{X|C})\langle -1/\dim X\rangle$ is preserved
by the general reduction modulo $p$, for $p$ sufficiently large.

Apply Proposition~\ref{prop:vanishp}\index{Vanishing result in characteristic $p$} to the vector bundle
$\sigma_k^*(T_{X_k|C_k})\langle -1/\dim X\rangle$. Observe that by
semicontinuity, the vanishing and global generation results obtained
in Proposition~\ref{prop:vanishp} extend to general deformations of
$\tau=\sigma_k\circ F$. Theorem~\ref{thm:descrHom} then immediately yields
the following.
\begin{corollary}\label{cor:rdc}
  There is an open neighborhood $\Omega \subset \Hom_{ C_k}(C_k[1],
  X_k)$ of $\tau$ such that
  \begin{enumerate}
  \item If $[\nu] \in \Omega$ is any morphism and $B \subset
    C_k[1]$ any subscheme of length $\#(B) \leq b_p(1/\dim X)$, then
    the relative deformations of $\nu$ over $C_k$ fixing $B$
    dominate $X_k$.
  \item If $T \subset X$ is the set of fundamental points of the
    birational map $q$ then the subset
    $$
    \Omega^0 = \{ [\nu] \in \Omega \,\, | \,\,
    (\nu)^{-1}(T_k) = \emptyset\}
    $$
    of morphisms whose images avoid $T_k$ is again open in
    $\Hom_{C_k}(C_k[1], X_k)$.
  \end{enumerate}
  \qed
\end{corollary}
 
In particular, Corollary~\ref{cor:rdc} states the following: given a
general element $[\nu]\in\Omega^0$ and a subscheme $B \subset
C_k[1]$ of length $\#(B) \leq b_p(1/\dim X)$, the deformations of
$q\circ\nu$ that fix $B$ dominate $Z_k$. Using Mori's\index{Mori, Shigefumi}
Bend-and-Break\index{Bend-and-Break}, Theorem~\ref{thm:bab2}, we obtain that $Z_k$ is
uniruled in curves of $H_{Z,k}$-degree at most
 \begin{align*}
   2\deg(q\circ\nu)^*(H_{Z,k})/\lfloor b_p(1/\dim X)
   \rfloor & =2\deg\nu^*(H_{X,k})/\lfloor b_p(1/\dim X) \rfloor \\
   & =2p\cdot \deg\sigma^*(H_X)/\lfloor b_p(1/\dim X) \rfloor   .
 \end{align*}
 The proof of Proposition~\ref{prop:pfBMcQ} and Theorem~\ref{thm:BMcQ}
 is finished if we note that this number is smaller than or equal to
 $d$, for $p>>0$\index{Vanishing threshold}.\qed

\section{An effective version of Miyaoka's criterion}
\label{sec:foliations}

As an immediate corollary of Theorem~\ref{thm:BMcQ}, we deduce an
effective version of Miyaoka's\index{Miyaoka, Yoichi} characterization of uniruledness\index{Miyaoka, Yoichi!criterion of uniruledness}. It
asserts that positive parts in the restriction of $T_X$ to a general
complete intersection curve are tangent to rational curves.  More
precisely, we use the following definition.

\begin{definition}
  Let $X$ be a normal projective variety, and $C \subset X$ a
  subvariety which is not contained in the singular locus of $X$, and
  not contained in the indeterminacy locus of the rationally connected
  quotient\index{Rationally connected quotient} $q: X \dasharrow Q$. If $\sF \subset T_X|_C$ is any
  subsheaf, we say that $\sF$ is \emph{vertical with respect to the
    rationally connected quotient}\index{Rationally connected quotient!subsheaf vertical with respect to}, if $\sF$ is contained in $T_{X|Q}$
  at the general point of $C$.
\end{definition}

The effective version of Miyaoka's\index{Miyaoka, Yoichi} criterion\index{Miyaoka, Yoichi!criterion of uniruledness} is then formulated as
follows.

\begin{corollary}\label{cor:Miyaoka}\index{Miyaoka, Yoichi!criterion of uniruledness}
  Let $X$ be a normal complex projective variety and $C \subset X$ a
  general complete intersection curve. Assume that the restriction
  $T_X|_C$ contains an ample locally free subsheaf $\sF_C$. Then
  $\sF_C$ is vertical with respect to the rationally connected
  quotient of $X$.
\end{corollary} 

This statement appeared first implicitly in
\cite[chap.~9]{SecondAsterisque}, but we believe there are issues with
the proof, see \cite[rem.~23]{KST05}. To our best knowledge, the
argument presented here gives the first complete proof of this
important result.

\subsection{Vector bundles over complex curves}

The proof of Corollary~\ref{cor:Miyaoka} relies on a number of facts
about the Harder-Narasimhan filtration of vector bundles on curves,
which are  known to the experts. For lack of an adequate
reference we include full proofs here. To start, we show that any
vector bundle on a smooth curve contains a maximally ample subbundle\index{Maximally ample subbundle}.

\begin{proposition}\label{prop:HN1-new}
  Let $C$ be a smooth complex-projective curve and $E$ a vector bundle
  on $C$, with Harder-Narasimhan filtration
  $$
  0 = E_0 \subset E_1 \subset \ldots \subset E_r = E.
  $$
  Let $\mu_i:=\mu (E_i/E_{i-1})$ be the slopes of the
  Harder-Narasimhan quotients.  Suppose that $\mu_1>0$ and let $k :=
  \max\{\, i\, |\, \mu_i>0 \}$.  Then $E_i$ is ample for all $1\leq i
  \leq k$ and every ample subsheaf of $E$ is contained in $E_k$.
\end{proposition}

\begin{definition}\label{def:masb}
  In the setup of Proposition~\ref{prop:HN1-new}, the bundle $E_k$ is
  called the \emph{maximal ample subbundle of $E$}.\index{Maximally ample subbundle}
\end{definition}

\begin{proof}[Proof of Proposition~\ref{prop:HN1-new}]
  Hartshorne's\index{Hartshorne, Robin} characterization of ampleness,
  \cite[thm.~2.4]{Hartshorne71}\index{Hartshorne, Robin!characterization of ampleness}, says that $E_i$ is ample iff all its
  quotients have positive degree. Dualizing, we have to prove that
  every subbundle of $E_i^\vee$ has negative degree, or, equivalently,
  that its maximal destabilizing subsheaf has negative slope, see
  \cite[1.3.4]{HL97}. This, however, holds because the uniqueness of
  the Harder-Narasimhan filtration of $E_i^\vee$, \cite[1.3.5]{HL97},
  implies that
  $$
  \mu_{\max}(E_i^\vee)=-\mu_i<0.
  $$ 
  
  To show the second statement, let $F \subset E$ be any ample
  subsheaf of $E$ and set
  $$
  j := \min\{\, i\, |\, F\subset E_i, \, 1\leq i\leq r \}.
  $$
  We need to check that $j\leq k$.  By the definition of $j$ and
  the ampleness of $F$, the image of $F$ in $E_j/E_{j-1}$ has positive
  slope. The semi-stability of $E_j/E_{j-1}$ therefore implies
  $\mu_j>0$ and $j\leq k$.
\end{proof}

Proposition~\ref{prop:HN1-new} says that the first few terms in the
Harder-Narasimhan filtration are ample. The following, related
statement will be used in the proof of Corollary~\ref{cor:Miyaoka}\index{Miyaoka, Yoichi!criterion of uniruledness} to
construct foliations on $X$.

\begin{proposition}\label{prop:ampleness}
  In the setup of Proposition~\ref{prop:HN1-new}, the vector bundles
  $E_j \otimes \bigl(\factor E.E_i.\bigr)^\vee$ are ample for all $0<j
  \leq i<r$.  In particular, if $E_i$ is any ample term in the
  Harder-Narasimhan Filtration of $E$, then $\Hom\bigl(E_i, \factor
  E.E_i.\bigr)$ and $\Hom\bigl(E_i\otimes E_i, \factor E.E_i.\bigr)$
  are both zero.\index{Harder-Narasimhan filtration!ample terms in}
\end{proposition}

\begin{remark}
  If $X$ is a polarized manifold whose tangent bundle contains a
  subsheaf of positive slope, Proposition~\ref{prop:ampleness} shows
  that the first terms in the Harder-Narasimhan filtration\index{Harder-Narasimhan filtration!of the tangent bundle} of $T_X$
  are special foliations\index{Foliation} in the sense of Miyaoka\index{Miyaoka, Yoichi}\index{Foliation!special in the sense of Miyaoka},
  \cite[sect.~8]{Miy85}. By \cite[thm.~8.5]{Miy85}, this already
  implies that $X$ is dominated by rational curves that are tangent to
  these foliations.
\end{remark}

\begin{proof}[Proof of Proposition~\ref{prop:ampleness}]
  As a first step, we show that the vector bundle
  $$
  F_{i,j} := \bigl(\factor E_j.E_{j-1}.\bigr)\otimes \bigl(\factor
  E.E_i.\bigr)^\vee
  $$
  is ample.  Assume not. Then, by Hartshorne's\index{Hartshorne, Robin} ampleness criterion\index{Hartshorne, Robin!characterization of ampleness}
  \cite[prop.~2.1(ii)]{Hartshorne71}, there exists a quotient $A$ of
  $F_{i,j}$ of degree $\deg_C A \leq 0$.  Equivalently, there exists a
  non-trivial subbundle
  $$
  \alpha : B \to F_{i,j}^\vee = \bigl(\factor
  E_j.E_{j-1}.\bigr)^\vee\otimes\bigl(\factor E.E_i.\bigr)
  $$
  with $\deg_C B \geq 0$. Replacing $B$ by its maximally
  destabilizing subbundle, if necessary, we can assume without loss of
  generality that $B$ is semistable. In particular, $B$ has
  non-negative slope $\mu(B\bigr) \geq 0$.  On the other hand, we have
  that $\bigl(\factor E_j.E_{j-1}.\bigr)$ is semistable.  The slope of
  the image of the induced morphism
  $$
  B\otimes \bigl(\factor E_j.E_{j-1}.\bigr) \to \bigl(\factor E.E_i.\bigr)
  $$
  will thus be larger than $\mu_{\max}\bigl(\factor E.E_i.\bigr)=
  \mu \bigl(\factor E_{i+1}.E_i.\bigr)$.  This shows that $\alpha$
  must be zero, a contradiction which proves the amplitude of
  $F_{i,j}$.
  
  With this preparation we will now prove
  Proposition~\ref{prop:ampleness} inductively. If $j=1$, then the
  above claim and the statement of Proposition~\ref{prop:ampleness}
  agree. Now let $1< j \leq i<r$ and assume that the statement was
  already shown for $j-1$.  Then consider the sequence
  $$
  0 \to \underbrace{E_{j-1}\otimes \bigl(\factor
    E.E_i.)^\vee}_{\text{ample}} \to E_j \otimes \bigl(\factor E.E_i.)^\vee
  \to \underbrace{(\factor E_j.E_{j-1}. )\otimes(\factor
    E.E_i.)^\vee}_{\text{ample}} \to 0
  $$
  But then also the middle term is ample, which shows
  Proposition~\ref{prop:ampleness}.
\end{proof}

\subsection{Proof of Corollary~\ref{cor:Miyaoka}}\index{Miyaoka, Yoichi!criterion of uniruledness}

We will show that the sheaf $\sF_C$, which is defined only on the
curve $C$ is contained in a foliation\index{Foliation} $\sF$ which is regular along $C$
and whose restriction to $C$ is likewise ample.
Corollary~\ref{cor:Miyaoka}\index{Miyaoka, Yoichi!criterion of uniruledness} then follows immediately from
Theorem~\ref{thm:BMcQ}.

An application of Proposition~\ref{prop:HN1-new} to $E := T_X|_C$
yields the existence of a locally free term $E_i \subset T_X|_C$ in
the Harder-Narasimhan filtration\index{Harder-Narasimhan filtration!ample terms in} of $T_X|_C$ which contains $\sF_C$
and is ample. The choice of $C$ then guarantees that $E_i$ extends to
a saturated subsheaf $\sF \subset T_X$.  The proof is thus finished if
we show that $\sF$ is a foliation, i.e.~closed under the Lie bracket.
Equivalently, we need to show that the associated
O'Neill tensor\index{O'Neill tensor}\footnote{The Lie bracket is of course not
  $\O_X$-linear. However, an elementary computation show that $N$ is
  well-defined and linear.}
$$
N : \sF \otimes \sF \to \factor T_X.\sF.
$$
vanishes. By Proposition~\ref{prop:ampleness}, the restriction of
the bundle
$$
Hom\left(\sF \otimes \sF, \factor T_X.\sF.\right) \cong (\sF \otimes
\sF)^{\vee} \otimes \factor T_X.\sF.
$$
to $C$ is anti-ample. In particular, 
$$
N|_C \in H^0\left( C, Hom\left(\sF \otimes \sF, \factor T_X.\sF.\right)|_C \right) = \{ 0\}.
$$
Ampleness is an open property, \cite[cor.~9.6.4]{EGA4-3}, so that the
restriction of $N$ to deformations $(C_t)_{t\in T}$ of $C$ stays zero
for most $t \in T$.  Since the $C_t$ dominate $X$, the claim follows.
This ends the proof of Corollary~\ref{cor:Miyaoka}\index{Miyaoka, Yoichi!criterion of uniruledness}. \qed

\section{The stability of the tangent bundle, and partial rationally connected quotients}
\label{sec:partRCQ}\index{Partial rationally connected quotients}

Recall that a complex variety $X$ is called $\Q$-Fano if a
sufficiently high multiple of the anticanonical divisor $-K_X$ is
Cartier and ample. The methods introduced above immediately yield that
$\Q$-Fano varieties whose tangent bundles are unstable allow sequences
of rational maps with rationally connected fibers.

\begin{corollary}\label{cor:fgnb}
  Let $X$ be a normal complex $\Q$-Fano variety and $L_1, \ldots,
  L_{\dim X-1} \in \Pic(X)$ be ample line bundles. Let 
  $$
  \{0\} = E_{-1} = E_0 \subset E_1 \subset \cdots \subset E_m = T_X
  $$
  be the Harder-Narasimhan filtration\index{Harder-Narasimhan filtration!of the tangent bundle} of the tangent sheaf with
  respect to $L_1, \ldots, L_{\dim X-1}$ and set
  $$
  k := \max \{ 0 \leq i \leq m \,|\, \mu(E_i/E_{i-1}) > 0\}.
  $$
  Then $k > 0$, and there exists a commutative diagram of dominant
  rational maps
  \begin{equation}
    \label{eq:ratlquot}
    \xymatrix{
      X \ar@{-->}[d]_{q_1} \ar@{=}[r] & X \ar@{-->}[d]_{q_2} \ar@{=}[r] & \cdots \ar@{=}[r] & X \ar@{-->}[d]_{q_k} \\
      Q_1 \ar@{-->}[r] & Q_2 \ar@{-->}[r] & \cdots \ar@{-->}[r] & Q_k, }
  \end{equation}
  with the following property: if $x \in X$ is a general point, and
  $F_i$ the closure of the $q_i$-fiber through $x$, then $F_i$ is
  rationally connected, and its tangent space at $x$ is exactly $E_i$,
  i.e., $T_{F_i}|_x = E_i|_x$.
\end{corollary}
\begin{proof}
  Let $C \subset X$ be a general complete intersection curve\index{General complete intersection curve} with
  respect to $L_1, \ldots, L_{\dim X-1}$. Since $c_1(T_X) \cdot C >
  0$, Proposition~\ref{prop:HN1-new} implies $k>0$ and that the
  restrictions $E_1|_C, \ldots, E_k|_C$ are ample vector bundles. We
  have further seen in Theorem~\ref{thm:BMcQ} that the $(E_i)_{1 \leq
    i \leq k}$ give a sequence of foliations with algebraic\index{Leaf of a foliation!algebraicity} and
  rationally connected leaves\index{Leaf of a foliation!rational connectivity}.
  
  To end the construction of Diagram~\eqref{eq:ratlquot}, let $q_i : X
  \dasharrow \Chow(X)$ be the map that sends a point $x$ to the
  $E_i$-leaf through $x$, and let $Q_i := {\rm Image}(q_i)$.
\end{proof}

\begin{remark}
  Corollary~\ref{cor:fgnb} also holds in the more general setup where
  $X$ is a normal variety whose anti-canonical class is represented by
  a Weil divisor with positive rational coefficients.
\end{remark}

\subsection{Open Problems}

It is of course conjectured that the tangent bundle of a Fano manifold
$X$ with $b_2(X)=1$ is stable. We are therefore interested in a
converse to Corollary~\ref{cor:fgnb} and ask the following.

\begin{question}\index{Partial rationally connected quotients!(in)dependence of polarization}
  Given a $\Q$-Fano variety and a sequence of rational maps with
  rationally connected fibers as in Diagram~\eqref{eq:ratlquot}, when
  does the diagram come from the unstability of $T_X$ with respect to
  a certain polarization? Is Diagram~\eqref{eq:ratlquot} characterized
  by universal properties?
\end{question}

\begin{question}
  To what extent does Diagram~\eqref{eq:ratlquot} depend on the
  polarization chosen?\index{Partial rationally connected quotients}
\end{question}

\begin{question}
  If $X$ is a uniruled manifold or variety, is there a polarization
  such that the rational quotient map comes from the Harder-Narasimhan
  filtration of $T_X$?\index{Partial rationally connected quotients}\index{Rationally connected quotient}\index{Harder-Narasimhan filtration!of the tangent bundle}
\end{question}

\section{Rational curves on quasi-projective manifolds}
\label{sec:KMcK}

Quasi-projective varieties appear naturally in a number of settings,
e.g., as moduli spaces or modular varieties. In this setup, it is
often not reasonable to ask if a given variety $S^\circ$ contains
complete rational curves\footnote{But see \cite[cor.~5.9]{KMcK} for
  criteria that can sometimes be used if $S^\circ$ can be compactified
  by a finite number of points.}.  Instead, one is interested in
hyperbolicity properties\index{Hyperbolicity problems} of $S^\circ$, i.e., one asks if there are
non-constant morphisms $\C \to U$, or if $S^\circ$ is dominated by
images of such morphisms ---we refer to \cite{Siu04} for a general
discussion, and to Chapter~\ref{chap:hyperbolicity} for a very brief
overview of the hyperbolicity question for moduli of canonically
polarized manifolds and for applications\index{Family of canonically polarized varieties}.

In this respect, a famous conjecture attributed to
Miyanishi\index{Miyanishi, Masayoshi!conjecture on affine lines on quasi-projective varieties}\index{Miyanishi, Masayoshi}\footnote{See \cite{GZ94} for a partial case.} suggests that
the following logarithmic analog of Conjecture~\ref{conj:mumford}
holds true.

\begin{conjecture}[Miyanishi]\label{conj:Miyanishi}\index{Miyanishi, Masayoshi}\index{Miyanishi, Masayoshi!conjecture on affine lines on quasi-projective varieties}
  Let $S^\circ$ be a smooth quasi-projective variety. Then
  $\kappa(S^\circ) = -\infty$ if and only if $S^\circ$ is dominated by
  images of $\C$.
\end{conjecture}

We briefly recall the definition of the logarithmic Kodaira dimension.

\begin{definition}\label{def:log-Kdim}\index{Logarithmic Kodaira dimension}
  Let $S^\circ$ be a smooth quasi-projective variety and $S$ a smooth
  projective compactification of $S^\circ$ such that $D := S \setminus
  S^\circ$ is a divisor with simple normal crossings. The
  \emph{logarithmic Kodaira dimension} of $S^\circ$, denoted by
  $\kappa(S^\circ)$, is defined as the Kodaira-Iitaka dimension
  $\kappa(K_S+D)$ of the line bundle $\O_S(K_S + D) \in \Pic(S)$.
  
  The variety $S^\circ$ is called \emph{of log general type} if
  $\kappa(S^\circ)=\dim S^\circ$, i.e., if the divisor $K_S+D$ is big.
\end{definition}

It is a standard fact in logarithmic geometry that a compactification
$S$ with the described properties exists, and that the logarithmic
Kodaira dimension\index{Logarithmic Kodaira dimension} does not depend on the choice of the
compactification, \cite[chap.~11]{Iitaka82}.

Conjecture~\ref{conj:Miyanishi} was studied by Miyanishi\index{Miyanishi, Masayoshi} and Tsunoda\index{Tsunoda, Shuichiro}
in \cite{Miyanishi-Tsunoda82, Miyanishi-Tsunoda84} and a number of
further papers, and by Zhang\index{Zhang, De-Qi} \cite{Zhang88}. Complete results for
surfaces were obtained by Keel-McKernan\index{Keel, Se\'an}\index{McKernan, James}, as follows.

\begin{theorem}[\protect{\cite[thm.~1.1]{KMcK}}]\label{thm:KMcK1}\index{Miyanishi, Masayoshi!conjecture on affine lines on quasi-projective varieties}
  Let $S^\circ$ be a smooth quasi-projective variety of dimension at
  most two. Then $\kappa(S^\circ) = -\infty$ if and only if $S^\circ$
  is dominated by images of $\C$. \qed
\end{theorem}

Keel-McKernan\index{Keel, Se\'an}\index{McKernan, James} also give conditions that guarantee that $U$ is
dominated by images of $\mathbb C^*$.\index{KMcK, rational curves on quasi-projective varieties}

\begin{theorem}[\protect{\cite[prop.~1.4]{KMcK}}]\label{thm:KMcK2}
  Let $S$ be a normal projective surface, let $D \subset S$ be a
  reduced curve, and set $U := S \setminus (D \cup {\rm Sing}(S))$.
  Consider the following conditions:
  \begin{enumerate}
  \item $K_S + D$ is numerically trivial, but not log-canonical.
  \item $K_S + D$ is numerically trivial, and $D \not = \emptyset$.
  \item $K_S$ is numerically trivial, $D = \emptyset$, and $S$ has a
    singularity which is not a quotient singularity.
  \end{enumerate}
  If any of the above hold, then $U$ is dominated by images of $\C^*$.
  \qed
\end{theorem}

\subsection{Open Problems}

Theorems~\ref{thm:KMcK1} and \ref{thm:KMcK2} were shown in \cite{KMcK}
using deformation theory on non-separated algebraic spaces and a
rather involved case-by-case analysis of possible curve configurations
on $S$. The analysis alone covers more than a hundred pages, and a
generalization to higher dimensions seems out of the question. We
would therefore like to pose the following problem.

\begin{problem}
  Find a more conceptual proof of Theorem~\ref{thm:KMcK1}, perhaps
  using methods introduced in Section~\ref{sec:RCfoliations}.
\end{problem}

% Local Variables:
% TeX-master: "RC-arXive.tex"
% End:

\part{Geometry of uniruled varieties}

\chapter{Geometry of rational curves on projective manifolds}
\label{chap:geomtry}

In Chapter~\ref{chap:existence} we have reviewed criteria to guarantee
that a given variety is covered by rational curves. In this chapter,
we assume that we are given a variety $X$ that is covered by rational
curves, and study the geometry of curves on $X$ in more detail.

\subsection*{Outline of the chapter}

In Section~\ref{sec:setupgeom} we review a number of known concepts
concerning rational curves on $X$. In particular, we give the
definition of a family of rational curves\index{Family of rational curves} and fix the notation that
will be used later on.

Section~\ref{sec:singminrat1} deals with the locus of singular curves
of a dominating family of minimal rational curves\index{Family of rational curves!dominating of minimal degrees|(}\index{Family of rational curves!dominating of minimal degrees!has few singular members}. We apply the
results of \cite{Kebekus02a} in Section~\ref{sec:singminratVMRT} to
ensure the existence and finiteness of the tangent morphism\index{Tangent morphism}\index{Tangent morphism!existence}\index{Tangent morphism!finiteness}, which
maps every minimal rational curve passing by a point $x$ to its tangent
direction at $x$. The image, i.e., the set of tangent directions to
which there exists a minimal degree rational curve, is called
\emph{variety of minimal rational tangents}\index{Variety of minimal rational tangents}, or \emph{VMRT}. Recently,
Hwang\index{Hwang, Jun-Muk} and Mok\index{Mok, Ngaiming} have shown that the general minimal degree rational
curve on $X$ is uniquely determined by its tangent direction at a
point\index{Tangent morphism!birationality}. This result gives more information about the VMRT.  We sketch
their argument in Section~\ref{sec:birtgtmor}.

We have claimed in the introduction that the VMRT determines the
geometry of $X$ to a large degree. In order to substantiate this
claim, we mention, in Section~\ref{sec:importance}, a number of
results in that direction. One particular result, an estimate for the
minimal number of rational curves required to connect two points on a
Fano manifold\index{Length!of a Fano manifold}, is explained in Section~\ref{sec:secants} in more
detail.

\section{The space of rational curves, setup of notation}
\label{sec:setupgeom}\index{Family of rational curves|(}\index{Space of rational curves on a variety}

We begin by defining some well known concepts about parameter spaces
of rational curves in projective varieties. We refer the reader to
\cite[II.2]{K96} for a detailed account.

\begin{definition}
  Let $X$ be a normal complex projective variety, and $\RatCurves(X)
  \subset \Chow(X)$ be the quasi-projective subvariety whose points
  correspond to irreducible and generically reduced rational curves in
  $X$. Let $\rat(X)$ be its normalization and $U$ the normalization of
  the universal family over $\RatCurves(X)$. We obtain a diagram as
  follows.
  \begin{equation}\label{eq:d2}
    \xymatrix{
      U \ar[d]_{\pi} \ar[rr]^{\iota}_{\text{evaluation morphism}} && X \\
      \rat(X)
    }
  \end{equation}
\end{definition}

\begin{remark}
  The normalization morphism $\rat(X) \to \Chow(X)$ is finite and
  generically injective, but not necessarily injective. It is
  therefore possible that points in $\rat(X)$ are not in 1:1
  correspondence with actual curves in $X$. We use the notation
  $[\ell]$ to denote points in $\rat(X)$, and $\ell$ for the
  associated curves.
\end{remark}

A standard cohomological argument for families of irreducible and
generically reduced rational curves shows that the morphism $\pi$ is
in fact a $\P^1$-bundle. The space of rational curves is often
described in terms of the Hom-scheme. The following theorem
establishes the link.

\begin{theorem}[\protect{\cite[II~thm.~2.15]{K96}}]\index{Space of rational curves on a variety!relation with the Hom-scheme}\index{Space of rational curves on a variety}
 There exists a diagram as follows
 $$
 \xymatrix{ \Hom^n_{\bir}(\P^1,X)\times \P^1
   \ar[rrr]^(.65){\text{quotient by natl.}}_(.65){\text{action of
     }\Aut(\P^1)} \ar@/^.7cm/[rrrrr]^{\text{univ. morphism}}
   \ar[d]_{\text{projection}} &&&
   U \ar[d]_{\pi}^{\P^1-\text{bundle}} \ar[rr]^{\iota}_{\text{evaluation morphism}} && X \\
   \Hom^n_{\bir}(\P^1,X) \ar[rrr]^{\text{quotient by
       natl.}}_{\text{action of }\Aut(\P^1)} &&& \rat(X) }
 $$
 where $\Hom_{\bir}(\P^1,X) \subset \Hom(\P^1,X)$ is the scheme
 parametrizing birational morphisms from $\P^1$ to $X$ and
 $\Hom^n_{\bir}(\P^1,X)$ its normalization. \qed
\end{theorem}

\begin{definition}\label{def:ratscheme}
  A \emph{maximal family of rational curves}\index{Family of rational curves!maximal} is an irreducible
  component $H \subset \rat(X)$. A maximal family $H$ is called
  \emph{dominating}\index{Family of rational curves!dominating}, if $\iota|_{\pi^{-1}(H)}$ dominates $X$. A
  dominating family $H$ is a \emph{dominating family of rational
    curves of minimal degrees}\index{Family of rational curves!dominating of minimal degrees} if the degrees of the associated
  rational curves on $X$ are minimal among all dominating families.
\end{definition}

\begin{example}\label{ex:clebsch}
  Let $X \subset \P^3$ be a general cubic surface. It is classically
  known that $X$ contains 27 lines, and that there exists a
  dominating family of conics, i.e., rational curves of degree 2, on
  $X$. In this case, the family of conics would be a dominating family
  of rational curves of minimal degree.
\end{example}

Example~\ref{ex:clebsch} shows that a dominating family of rational
curves of minimal degrees needs not be proper ---there are sequences of
conics on $X$ whose limit cycle is a union of two lines. It is,
however, true that given a sequence of rational curves on any $X$, the
limit cycle is composed of rational curves. This observation
immediately gives the following properness statement for families of
curves of minimal degrees.

\begin{theorem}[\protect{\cite[IV~cor.~2.9]{K96}}]\index{Family of rational curves!dominating of minimal degrees!are generically unsplit}\index{Space of rational curves on a variety!through a given point is proper if degree is minimal}
  Let $H \subset \rat(X)$ be a dominating family of rational curves of
  minimal degrees. If $x \in X$ is a general point, then the subspace
  of curves through $x$,
  $$
  H_x := \pi\left(\iota^{-1}(x)\right) \cap H = \{ \ell \in H \,|\, x \in \ell\}
  $$
  is proper. \qed
\end{theorem}

\begin{definition}
  Let $H \subset \rat(X)$ be a dominating family of rational curves.
  If $H$ is proper, it is often called \emph{unsplit}\index{Family of rational curves!unsplit}. If for a
  a general point $x \in X$ the associated subspace $H_x$ is
  proper, the family is called \emph{generically unsplit}.\index{Family of rational curves!generically unsplit}
\end{definition}

One of the key tools in the description of rational curves on a
manifold $X$ is an analysis of the restriction of $T_X$ to the
rational curves in question ---recall that any vector bundle on $\P^1$
can be written as a sum of line bundles. The following proposition
summarizes the most important facts.

\begin{proposition}[\protect{\cite[Prop.~4.14]{Debarre01}, \cite[IV Cor.~2.9]{K96}}]\label{prop:gen=stand}
  Let $X$ be a smooth complex projective variety and $H
  \subset\rat(X)$ be a dominating family of rational curves of minimal
  degrees. If $x\in X$ is a general point and $[\ell] \in H_x$ any
  element, with normalization $\eta : \P^1 \to \ell$, then the
  following holds.
  \begin{itemize}
  \item $H$ is smooth at $[\ell]$.
  \item The normalization $\tilde H_x$ of $H_x$ is smooth.
  \item $\ell$ is \emph{free}\index{Rational curve!free}, i.e. the vector bundle $\eta^*(T_X)$ is
    nef.
  \item If $[\ell]$ is a general point of $H_x$, then $\ell$ is
    \emph{standard}\index{Rational curve!standard}, i.e., there exists a number $p$ such that
    $$
    \eta^*(T_X) \cong \sO_{\P^1}(2) \oplus \sO_{\P^1}(1)^{\oplus p} \oplus \sO_{\P^1}^{\oplus
      \dim(X)-1-p}.
    $$\index{Rational curve!splitting type}
  \end{itemize}
  \qed
\end{proposition}

\subsection{Open Problems}

\begin{problem}
  To what extent are dominating families of rational curves of minimal
  degrees characterized by their generic unsplitness?
\end{problem}

\begin{problem}
  Let $X$ be smooth and $H$ be a dominating family of rational curves
  of minimal degrees. Assume that for general $x \in X$, the space
  $H_x$ is positive dimensional. Is it true that $H_x$ is irreducible?
  See \cite[sect.~5.1]{KK04} for a partial case.
\end{problem}
\index{Family of rational curves|)}

\section{Singular rational curves}\label{sec:singminrat1}\index{Rational curve!singular|(}\index{Family of rational curves!dominating of minimal degrees!has few singular members|(}

Let $X \subset \P^n$ be a projective manifold that is covered by
lines. If $x \in X$ is any point, we can consider the space $\mathcal
C_x \subset \P(T_X|^\vee_x)$ of tangent directions that are tangent to
lines through $x$. The so-defined \emph{variety of minimal rational
  tangents}\index{Variety of minimal rational tangents} is a very classical and important object that has been
studied in the past by Cartan\index{Cartan, \'Elie} and Fubini\index{Fubini, Guido}, and a number of other
projective differential geometers.

We will, in this section, give a similar construction for families of
rational curves of minimal degrees. The main obstacle is that these
curves may be singular, which makes it difficult to properly define
tangents to them. It is, however, well understood that minimal degree
curves have only mild singularities at the general point of $X$.

\begin{definition}
  We say that a curve $C$ is \emph{immersed}\index{Rational curve!immersed}, if its normalization
  morphism $\tilde{C}\rightarrow C$ has rank $1$ at every point.
\end{definition}

\begin{theorem}[\protect{\cite[thm.~3.3]{Kebekus02a}}]\label{thm:dimsingcurves}\index{Rational curve!singular}\index{Family of rational curves!dominating of minimal degrees!has few singular members}\index{Family of rational curves!dominating of minimal degrees|)}
  Let $X$ be normal and $H \subset \rat(X)$ a dominating family of
  rational curves of minimal degrees. Further, let $x \in X$ be a
  general point, and consider the closed subvarieties
  \begin{align*}
    H_x^{\Sing} & := \{ [\ell] \in H_x \,|\, \ell \text{ is singular}\} \mbox{ and} \\
    H_x^{\Sing,x} & := \{ [\ell] \in H_x \,|\, \ell \text{ is singular
      at }x\}.
  \end{align*}
  Then the following holds.
  \begin{enumerate}
  \item The space $H_x^{\Sing}$ has dimension at most one, and the
    subspace $H_x^{\Sing,x}$ is at most finite. Moreover, if
    $H_x^{\Sing,x}$ is not empty, the associated curves are immersed.
  \item If there exists a line bundle $L \in \Pic(X)$ that intersects
    the curves with multiplicity $2$ then $H_x^{\Sing}$ is at most
    finite and $H_x^{\Sing,x}$ is empty.
  \end{enumerate}
\end{theorem}

The main idea in the proof of Theorem~\ref{thm:dimsingcurves} is the
observation that an arbitrary family of singular rational curves, like
any family of higher genus curves, is hardly ever projective ---see
\cite{Kebekus02c} for worked examples. An analysis of the projectivity
condition yields the statement.

\subsection{Sketch of proof of Theorem~\ref{thm:dimsingcurves}}

As usual, we subdivide the proof into several steps. We will only give
an idea how to show that the subspace $H_x^{\Sing,x}$ is at most
finite.

\subsection*{Step 1: Dimension count}

We assume that $H^{\Sing,x}_x$ is not empty because otherwise there is
nothing to prove.  A technical dimension count ---which we are not
going to detail in this sketch--- shows that
$$
\dim H_{x}^{\Sing }\geq \dim H_{x}^{\Sing , x}+1.
$$
Thus, the assumption implies that $\dim H_{x}^{\Sing }\geq 1$. We
fix a proper 1-dimensional subfamily $H' \subset H^{\Sing }_{x}$.

\subsection*{Step 2: A partial resolution of singularities}\index{Partial resolution of singularities}

Recall that $H' \subset \rat(X)$ has a natural morphism into the
Chow variety of $X$. Let $\pi : U' \to H'$ be the pull-back of the
universal family. We aim to replace $U'$ by a family where all fibers
are singular plane cubics. For that, consider the normalization
diagram.
\begin{equation}\label{eq:norm_diag}
\xymatrix{
  {\tilde U}\ar[rr]^{\eta}_{\txt{\tiny normalization}} \ar@/_/ [rd]_
  {\tilde \pi} & & {U'} \ar@/^/[ld]^{\pi} \\ & {H'}& } 
\end{equation}
After performing a series of finite base changes, if necessary, we can
assume that the following holds:
\begin{enumerate}
\item $H'$ is smooth.
\item $\tilde{U}$ is a $\P^{1}$-bundle over $H'$ ---see
  \cite[thm.~II.2.8]{K96}.
\item There exists a curve $s\subset U'_{\Sing }$ contained in the
  singular locus of $U'$ such that $\pi|_{s}$ is an isomorphism. For this,
  let $s$ be the normalization of a suitable component of $U'_{\Sing}$.
\item There exists a subscheme $\tilde{s}\subset \eta ^{-1}(s)$ whose
  restriction to all $\tilde{\pi }$-fibers is of length 2. For this,
  let $\tilde{s}$ be the normalization of a curve in $\Hilb _{2}(\eta
  ^{-1}(s)/H')$ and note that the relative $\Hilb $-functor commutes
  with base change.
\end{enumerate}
\begin{figure}
  \begin{center}
    $$
    \xymatrix{ 
      {\txt{$\tilde U$ \\ \includegraphics[width=3.5cm]{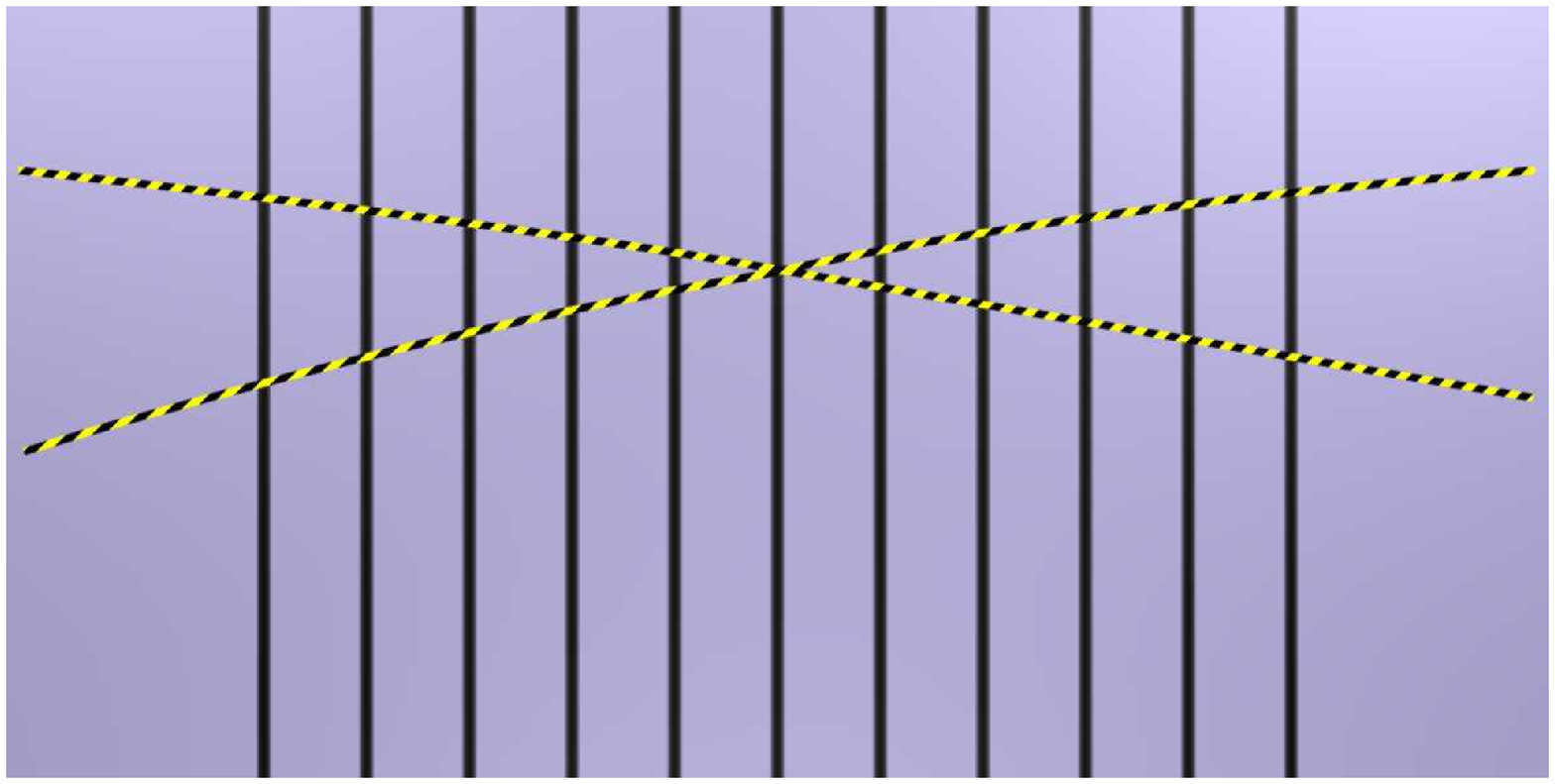} \\ \scriptsize $\P^1$-bundle}  }
      \ar[r]^(.45){\alpha}
      \ar@/_/ [rd]_ {\tilde \pi} & 
      { \txt{$\hat U$ \\ \includegraphics[width=4.5cm]{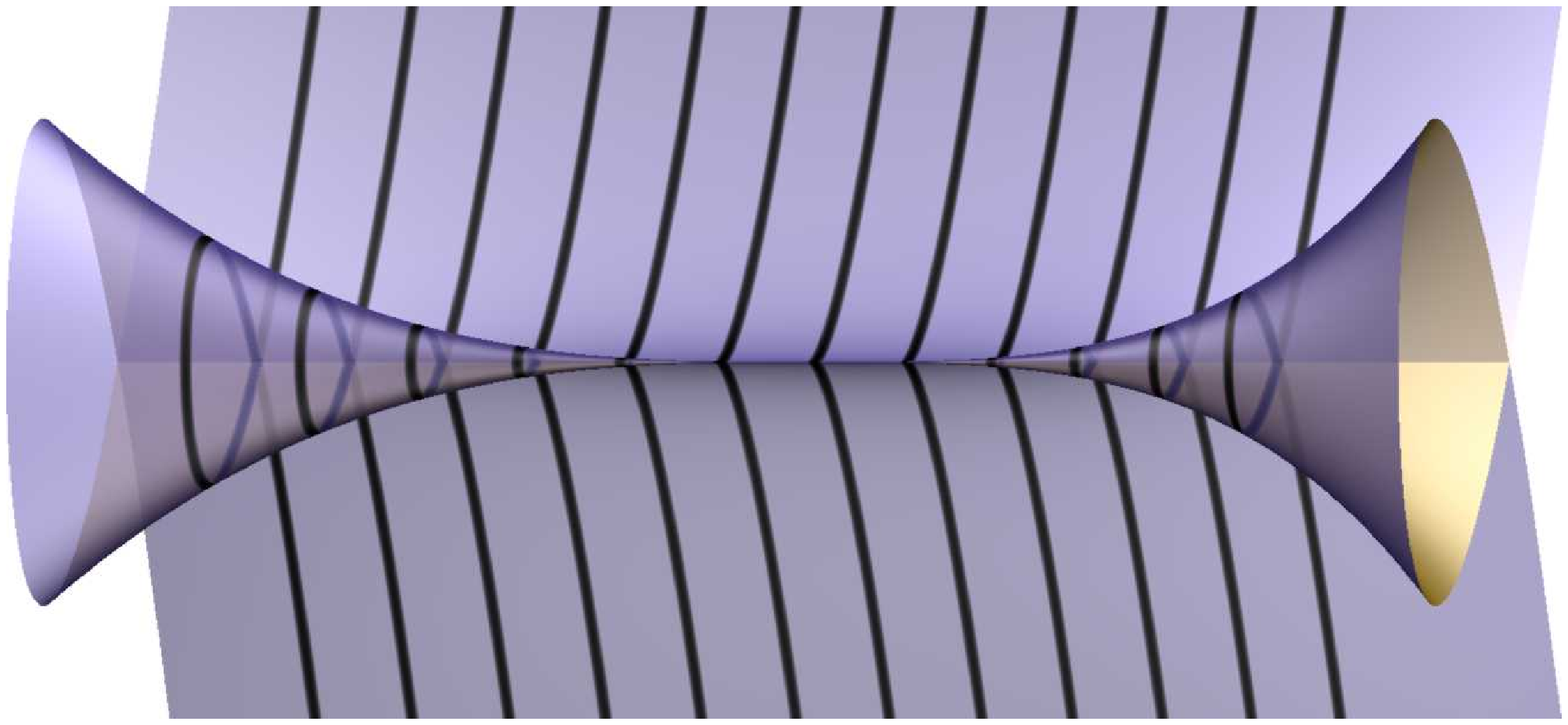} \\ \scriptsize family of plane cubics}  }
      \ar[r]^(.65){\beta}
      \ar[d]^(.7){\hat \pi}   & 
      {\txt{\ \\ \ \\ $U'$ \\ \includegraphics[width=1cm]{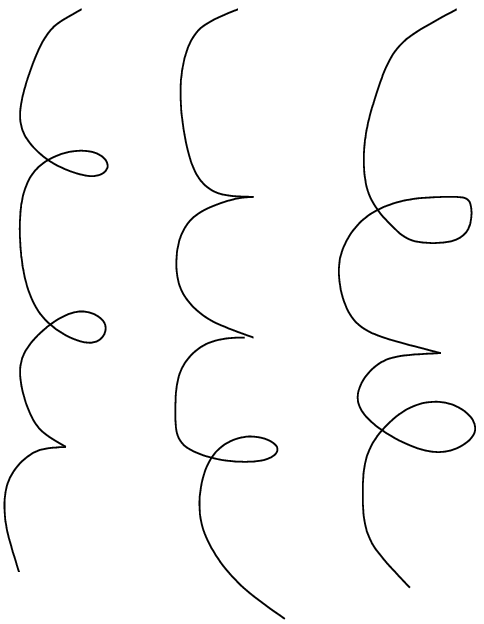} \\ \scriptsize family of \\ \scriptsize very singular \\ \scriptsize curves}} 
      \ar@/^/[ld]^(.6){\pi} \\ 
      & {\includegraphics[width=3cm]{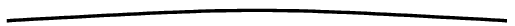} H' }& } 
    $$
    \caption{Replacing singular curves by plane cubics}
    \label{fig:plane_cubic}
  \end{center}
\end{figure}
We would like to extend the diagram~\eqref{eq:norm_diag} to
$$
\xymatrix{ {\tilde U} \ar@{-->}[rrd]_{\exists \alpha}
  \ar[rrrr]^{\eta}_{\txt{\tiny normalization}} \ar@/_/ [rrdd]_ {\tilde
    \pi} & & & & {U'} \ar@/^/[lldd]^{\pi}~ \\ & & {\hat U}
  \ar@{-->}[rru]_{\exists \beta} \ar@{-->}[d]_{\exists \hat\pi} & & \\ 
  & & {H'}& &}
$$
where all fibers of $\hat \pi$ are rational curves with a single
cusp or node, i.e., isomorphic to a plane cubic.
Figure~\ref{fig:plane_cubic} depicts this setup. Here we explain only
how to do this locally.

Knowing that $\tilde{U}$ is a $\P^1$-bundle over $H'$, we find an
(analytic) open set $V\subset H'$ with coordinate $v$, identify an
open subset of $\tilde{\pi }^{-1}(V)$ with $V\times \C $, choose a
bundle coordinate $u$ and write
$$
\tilde{s}=\{u^{2}=f(v)\}
$$
where $f$ is a function on $V$. We would then define $\alpha$ to be
$$
\begin{array}{cccl}
\alpha :  & V\times \C  & \longrightarrow  & V\times \C ^{2}\\
 & (v,u) & \mapsto  & (v,u^{2}-f(v),u(u^{2}-f(v)))
\end{array}
$$
A direct calculation shows that these locally defined morphisms
glue together to give a global morphism $\alpha :\tilde{U}\to \hat U$,
that a morphism $\beta: \hat U \to U'$ exists and that the induced map
$\hat \pi := \pi \circ \beta $ has the desired properties.

\subsection*{Step 3: Ruling out several cases}

In order to conclude in the next step, we have to rule out several
possibilities for the geometry of $\hat U$. We do this in every case
by a reduction to the absurd.

\subsection*{Step 3a: The case where all curves are immersed }

\begin{figure}
  \begin{center}
    $$
    \xymatrix{ 
      {\txt{$\tilde U$ \\ \includegraphics[width=4cm]{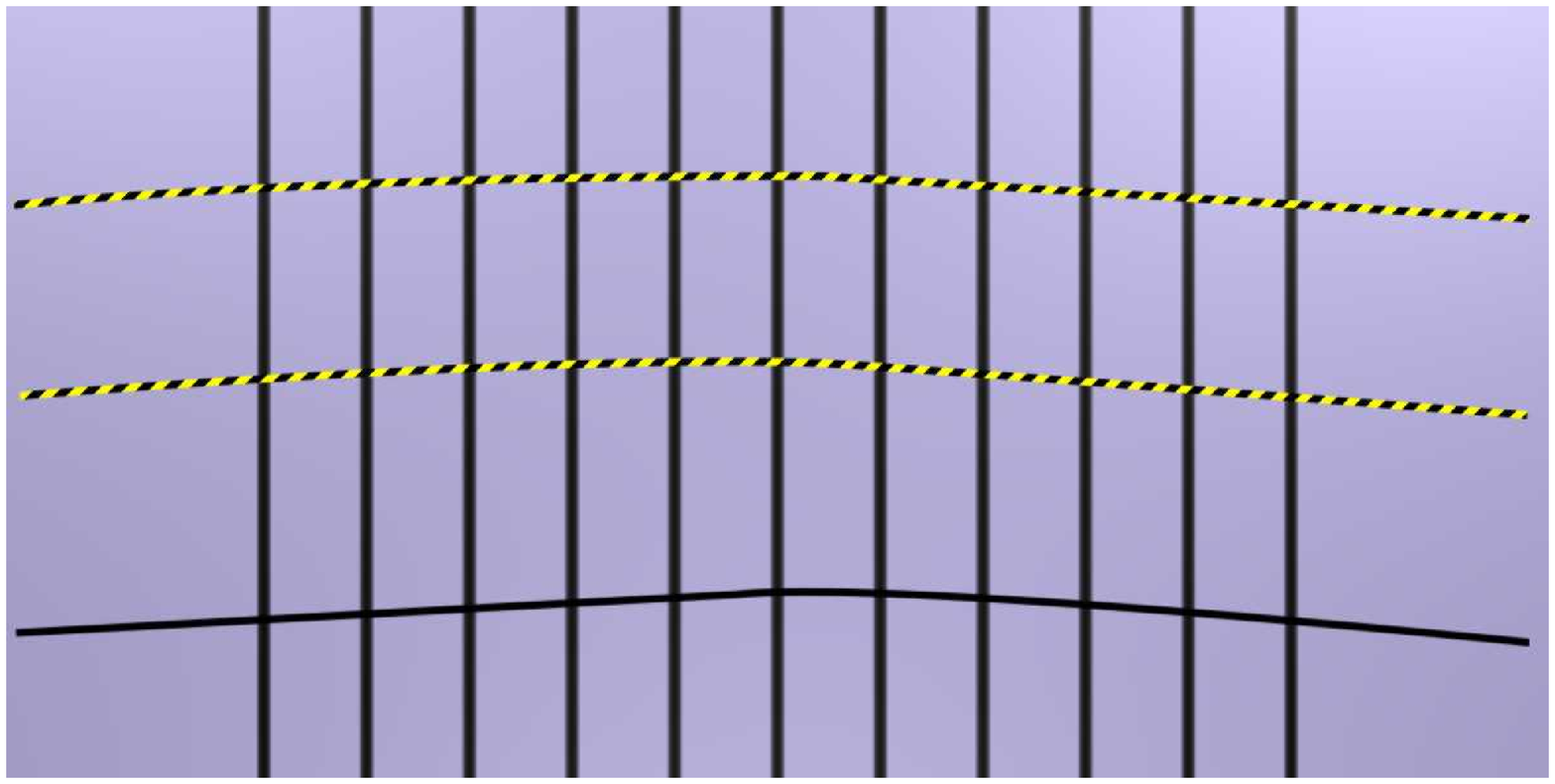} \\ \scriptsize $\P^1$-bundle}  }
      \ar[rr]^(.47){\alpha}_(.47){\txt{\tiny normalization}} 
      \ar@/_/[rd]_(.65){\tilde \pi} & & 
      {\txt{$\hat U$ \\ \includegraphics[width=5cm]{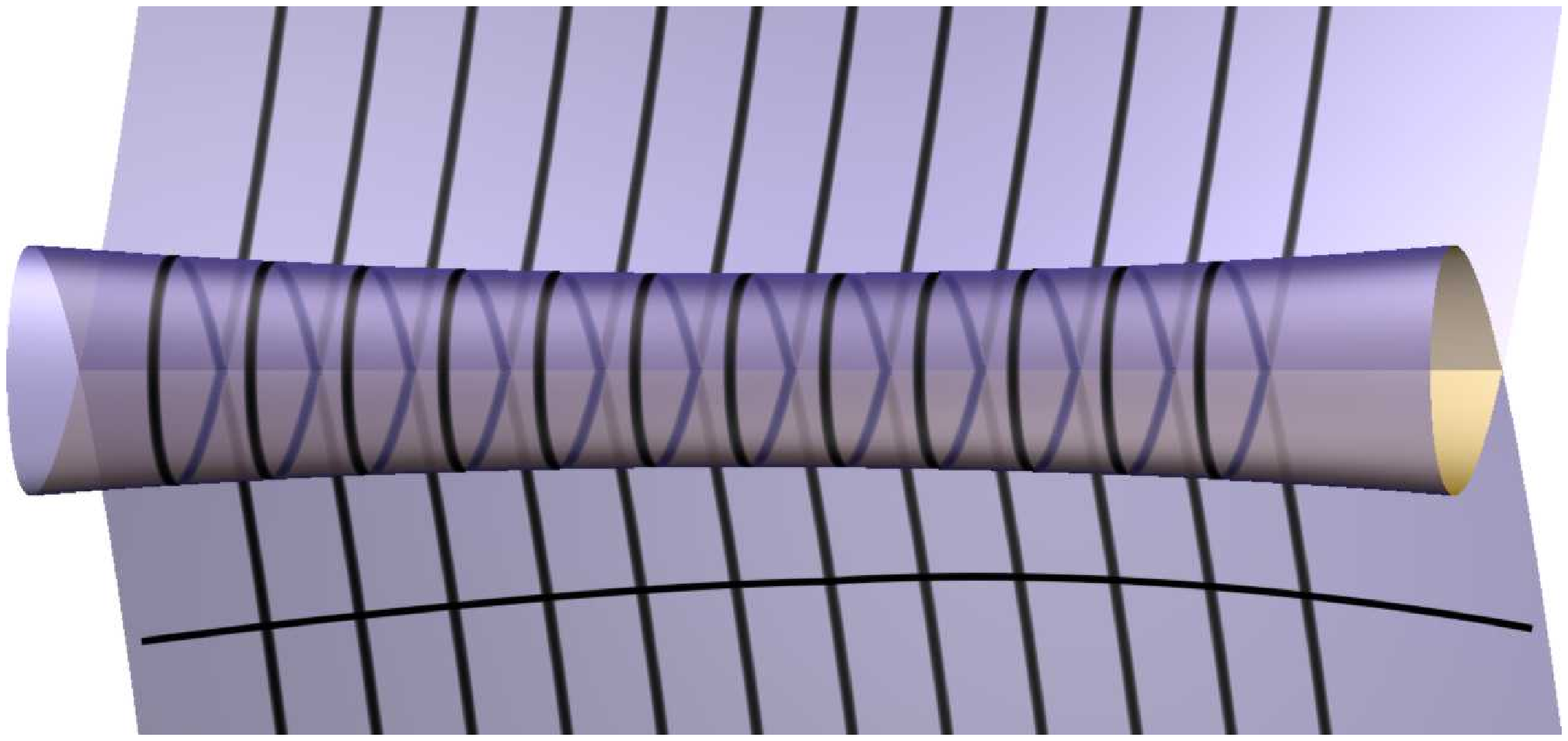} \\ \scriptsize family of immersed curves}} 
      \ar@/^/[ld]^(.65){\hat \pi} \\ 
      & {
        \begin{picture}(.5,.3)(0,0)
          \put(-1.2,0){\includegraphics[width=3cm]{base}}
          \put(1.9,-0.1){$H'$}
        \end{picture}
      }& } 
    $$
    \caption{Family of immersed curves}
    \label{fig:immersed}
  \end{center}
\end{figure}

If all curves associated with $H'$ are immersed, then the construction
outlined above will automatically give a family $\hat U$ where all
fibers are isomorphic to nodal plane cubics ---see
Figure~\ref{fig:immersed}.  Let $\sigma_{\infty }\subset \tilde{U}$ be
the section which is contracted to the point $x\in X$ (drawn as a
solid line) and consider the preimage of the singular locus
$\alpha^{-1}(\hat U_{\Sing})$. After another finite base change, if
necessary, we may assume that this set decomposes into two disjoint
components $\alpha^{-1}(\hat U_{\Sing})=\sigma_0\cup \sigma_1$, drawn
as dashed lines. That way we obtain three sections $\sigma_0$,
$\sigma_1$ and $\sigma_{\infty}$ in $\tilde{U}$, where
$\sigma_{\infty}$ can be contracted to a point and $\sigma_0$,
$\sigma_1$ are disjoint. But then it follows from an elementary
calculation with intersection numbers that either $\sigma_0 =
\sigma_{\infty}$ or that $\sigma_1 = \sigma_{\infty}$. This however,
is impossible, because then the Stein factorization of $\iota_x\circ
\eta$ would contract both $\sigma_0$ and $\sigma_1$, but ruled
surfaces allow at most a single contractible section.

\begin{remark}
  This setup has already been considered by several authors. See
  e.g.~\cite{CS95}.
\end{remark}

\subsection*{Step 3b: The case where no curve is immersed}

\begin{figure}
  \begin{center}
    $$
    \xymatrix{ 
      {\txt{$\tilde U$ \\ \includegraphics[width=4cm]{bundle1} \\ \scriptsize $\P^1$-bundle}  }
      \ar[rr]^(.45){\alpha}_(.45){\txt{\tiny normalization}} 
      \ar@/_/[rd]_(.65){\tilde \pi} & & 
      {\txt{$\hat U$ \\ \includegraphics[width=5cm]{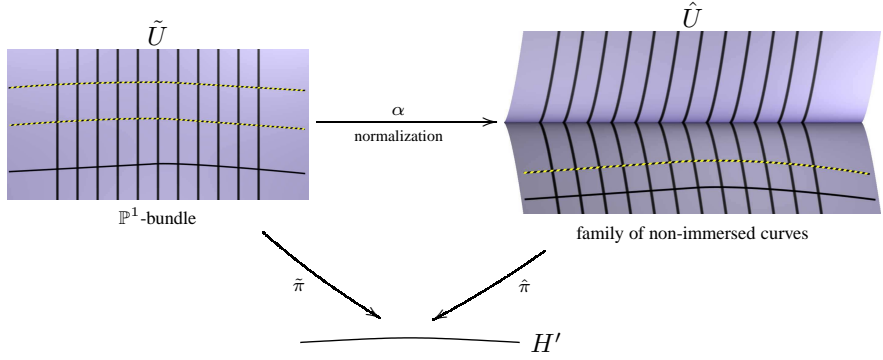} \\ \scriptsize family of non-immersed curves}} 
      \ar@/^/[ld]^(.65){\hat \pi} \\ 
      & {
        \begin{picture}(.5,.3)(0,0)
          \put(-1.2,0){\includegraphics[width=3cm]{base}}
          \put(1.9,-0.1){$H'$}
        \end{picture}
        }& } 
    $$
    \caption{Family of non-immersed curves}
    \label{fig:non-immersed}
  \end{center}
\end{figure}

If none of the curves associated with $H'$ is immersed, then the
curves $s$ and $\tilde{s}$ in the construction can be chosen so that
$\hat U \to H'$ is a family of cuspidal plane cubics, see
Figure~\ref{fig:non-immersed}. Again, let $\sigma_{\infty} \subset
\tilde{U}$ be the section which can be contracted and let $\sigma
_{0}$ be the preimage of the singularities.  That way we obtain two
sections.

In order to obtain a third one, remark that if $C\subset \P^2$ is a
cuspidal plane cubic and $H\in \Pic (C)$ a line bundle of positive
degree $k>0$, then there exists a unique point\footnote{If $H =
  \O_{\P^2}(1)|_C$, this is the classical inflection point. In
  general, this will be the hyperosculating point associated with the
  embedding given by $H$.} $y \in C$, contained in the smooth locus of
$C$, such that $\O_C(ky)\cong H$.  Thus, using the pull-back of an
ample line bundle $L\in \beta ^{*}\Pic (X)$, we obtain a third section
$\sigma_1$ which is disjoint from $\sigma_0$. Now conclude as above.
This time, however, it is not obvious that neither $\sigma_0$ nor
$\sigma_1$ coincides with $\sigma_{\infty}$.  Actually, this is true
because $x\in X$ was chosen to be a general point.  The proof of this
is rather technical therefore omitted.

\subsection*{Step 4: End of proof}

To end the proof, we argue by contradiction and assume that $\dim
H_x^{\Sing,x}\geq 1$. We have seen in Step~1 that this
implies $\dim H_x^{\Sing}\geq 2$. The argumentation of Step~3
implies that points of $\dim H_x^{\Sing}\geq 2$ correspond to both
immersed and non-immersed singular rational curves in $X$. 

Elementary deformation theory, however, shows that the closed
subfamily of non-immersed curves is always of codimension at least
one. Thus, the subfamily $H_x^{\Sing, \rm ni} \subset H_x^{\Sing }$ of
non-immersed curves is proper and positive-dimensional. This again has
been ruled out in Step~3b. \qed

\subsection{Open Problems}

\begin{problem}
  Examples show that Theorem~\ref{thm:dimsingcurves} is optimal for
  normal varieties $X$. It is not clear to us if
  Theorem~\ref{thm:dimsingcurves} can be improved if $X$ is smooth.
\end{problem}
\index{Rational curve!singular|)}\index{Family of rational curves!dominating of minimal degrees!has few singular members|)}

\section{The tangent morphism and the variety of minimal rational tangents}\label{sec:singminratVMRT}\index{Tangent morphism|(}\index{Variety of minimal rational tangents|(}

We apply Theorem~\ref{thm:dimsingcurves} to show the existence and
finiteness of the \emph{tangent morphism}, an important tool in the
study of uniruled manifolds that encodes the infinitesimal
behavior of $H$ near a general point $x \in X$.

\begin{definition}
  Given a dominating family $H$ of rational curves of minimal degrees
  on $X$, and a general point $x\in X$ we consider a rational
  \emph{tangent map}\index{Tangent map}
  $$
  \begin{array}{rccc}
    t_x: & H_x & \dashrightarrow & \P(T_X|^\vee_x) \\
    & \ell & \to & \P(T_{\ell}|^\vee_x)
  \end{array}
  $$
  associating to a curve through $x$ its tangent direction at $x$.
\end{definition}

\begin{theorem}[\protect{\cite[thm.~3.4]{Kebekus02a}}]\label{thm:tgtmor}
  With the notation as above, let $x \in X$ be a general point, let
  $\tilde H_x$ be the normalization of $H_x$, and let
  $$
  \tau_x : \tilde H_x \dashrightarrow \P(T_X|^\vee_x)
  $$
  be the composition of $t_x$ with the normalization morphism. Then
  $\tau_x$ is a finite morphism, called \emph{tangent morphism}.\index{Tangent morphism}
\end{theorem}
\begin{proof}
  Consider the pull-back of Diagram~\eqref{eq:d2} from
  page~\pageref{eq:d2},
  $$
  \xymatrix{
    U_x \ar[d]_{\pi_x} \ar[rr]^{\iota_x}_{\text{evaluation}} && X \\
    \tilde H_x }
  $$
  Theorem~\ref{thm:dimsingcurves}\index{Rational curve!singular} asserts that the preimage
  $\iota_x^{-1}(x)$ contains a reduced section $\sigma_\infty \cong
  \tilde H_x$, and at most finitely many points.  Since all curves are
  immersed at $x$, the tangent morphism of $\iota_x$ gives a
  nowhere vanishing morphism of vector bundles,
  \begin{equation}\label{eq:d3}
  T\iota_x : T_{U_x|\tilde H_x}|_{\sigma_\infty} \to \iota_x^*(T_X|_x).
  \end{equation}
  The tangent morphism is then given by the projectivization of
  \eqref{eq:d3}. 

  Assuming that $\tau_x$ is not finite, Equation~\eqref{eq:d3} asserts
  that we can find a curve $C \subset \tilde H_x$ such that
  $N_{\sigma_\infty,U_x}$ is trivial along $C$. But $\sigma_\infty$ can be
  contracted, and the normal bundle must thus be negative. 
\end{proof}

With these preparations, we can now introduce one of the central
objects of this survey, the \emph{variety of minimal rational
  tangents}, or \emph{VMRT}.\index{Variety of minimal rational tangents}

\begin{definition}\label{def:varmintgts}\index{Variety of minimal rational tangents}
  If $x \in X$ is a general point, we call the image $ \sC_x :=
  \tau_x(\tilde H_x) \subset \P(T_X|^\vee_x) $ the \emph{variety of
    minimal tangents} of $H$ at $x$. The subvariety
  $$
  \sC:=\mbox{closure of}\bigcup_{x \mbox{ general in
    }X}\sC_x\subset\P(T_X^\vee)
  $$
  is called the \emph{total variety of minimal rational tangents}\index{Variety of minimal rational tangents!total}
  of $H$.
\end{definition}

\begin{remark}
The projectivized tangent map of the evaluation morphism yields a diagram
$$
\xymatrix{ & & \sC \ar[d]^{\rho, \text{ projection}} &
  \hspace{-0.9cm}\subset\P(T_X^\vee) \\
  U \ar[rr]^{\iota}_{\text{evaluation}} \ar@/^.3cm/@{-->}[urr]^{\tau} & & X  }
  $$
  that we will later also use to describe the tangent map $\tau_x$
  at a general point. We call $\tau$ the \emph{global tangent map}\index{Tangent map!global}.
\end{remark}

The variety of minimal rational tangents has been extensively studied
by several authors, including Hwang\index{Hwang, Jun-Muk} and Mok\index{Mok, Ngaiming}. It can be computed in a
number of examples of uniruled varieties, such as Fano hypersurfaces,
rational homogeneous spaces\index{Rational homogeneous space} and moduli spaces of vector bundles\index{Moduli space of vector bundles on a curve}. We
refer the reader to \cite{Hwa00} for examples, but see also
Chapter~\ref{chap:moduli} below.

For varieties covered by lines, the situation is particularly easy.

\begin{proposition}[\protect{\cite[prop.~1.5]{Hwa00}}]\label{prop:Cxforlines}\index{Variety of minimal rational tangents!tangent space for varieties covered by lines}
  Suppose that there exists an embedding $X \subset \P^N$ such that
  the curves associated with points of $H$ are lines. If $x \in X$ is
  a general point, then $\tau_x$ is an embedding and $\sC_x$ is
  smooth. 
  
  If $[\ell] \in H_x$ is any line through $x$, then the projective
  tangent space to $\sC_x$ at $\P(T_{\ell}|_x)$ is exactly
  $\P\left((T_X|_\ell^+)^\vee\right)$, where $T_X|_\ell^+ \subset
  T_X|_\ell$ is the maximal ample subbundle.
\end{proposition}

In general, there is a direct and well understood relation between the
splitting type of $T_X|_\ell^+$ and the Zariski tangent space of
$\sC_x$. See \cite[sect.~1]{Hwa00} for details.

\subsection{Open Problems}

\begin{problem}
  In all smooth examples that we are aware of, the variety $\mathcal
  C_x$ is smooth, and has good projective-geometrical properties. What
  can be said in general? See \cite{Hwa00} for more detailed lists of
  problems.
\end{problem}

\section{Birationality of the tangent morphism}\label{sec:birtgtmor}\index{Tangent morphism!birationality|(}

We will later see that the projective geometry of the VMRT $\mathcal
C_x \subset \P(T_X|_x^\vee)$ encodes a lot of the geometrical
properties of $X$. One of the main difficulties in the applications is
that the variety $\mathcal C_x$ of minimal rational
tangents might be singular or reducible.  To overcome this difficulty one
often studies the tangent morphism $\tau_x : \tilde H_x \to \mathcal
C_x$. Since $\tilde H_x$ is smooth, one asks if $\tau_x$ is injective,
and if it has maximal rank ---this question can sometimes be answered
in the examples.

In general, it has been shown by Hwang\index{Hwang, Jun-Muk} and Mok\index{Mok, Ngaiming} that the normalization
of $\sC_x$ is smooth.\index{Variety of minimal rational tangents!has smooth normalization} This is a direct consequence of
Proposition~\ref{prop:gen=stand} and the main result of \cite{HM04}.

\begin{theorem}[\protect{\cite[thm.~1]{HM04}}]\label{thm:birtangmap}\index{Tangent map!generic injectivity}
  With the same notation as above, the global tangent map $\tau$ is
  generically injective. In particular, the normalization of $\sC_x$
  is smooth for general $x\in X$.
\end{theorem}

See \cite{KK04} for related criteria to guarantee that $\tau_x$ is in fact injective.

\begin{remark}
  A line in $\P^n$ is specified by giving a point $x \in \P^n$ and a
  tangent direction at $x$. Theorem~\ref{thm:birtangmap} says that a
  similar statement holds for minimal degree curves.\index{Rational curve!of minimal degree!is determined by a point and a tangent direction} See
  \cite[sect.~3.3]{Kebekus02a} for the related question if a minimal
  degree curve can be specified by two points.\index{Rational curve!of minimal degree!is determined by two points}
\end{remark}

Theorem~\ref{thm:birtangmap} was known to be true in the case where
$\sC_x = \P(T_X^\vee|_x)$, where it follows as a by-product of a
characterization of $\P^n$.

\begin{theorem}[\protect{\cite{CMSB02}, \cite{Kebekus02b}}]\label{thm:CMSB02}\index{Characterization of $\P^n$}
  Let $X$ be an irreducible normal projective variety of dimension
  $n$. Let $H \subset \rat(X)$ be a dominating family of rational
  curves of minimal degrees and assume that $\mathcal C =
  \P(T_X^\vee)$.  Then there exists a finite morphism $\P^n\rightarrow
  X$, \'etale over $X \setminus \Sing(X)$ that maps lines in $\P^n$ to
  curves parametrized by $H$. In particular, $X \cong \P^n$ if $X$ is
  smooth. \qed
\end{theorem}

Hwang\index{Hwang, Jun-Muk} and Mok\index{Mok, Ngaiming} have applied the theory of differential systems to $\sC$
in order to reduce the general case to that of
Theorem~\ref{thm:CMSB02}. We give a very rough sketch of the proof and
refer to \cite{HM04} for details.

\subsection{Sketch of  proof of Theorem~\ref{thm:birtangmap}}

\subsection*{Step 1: Differential systems in uniruled varieties}\index{Differential system}

\begin{definition}\label{def:distri}
  A \emph{distribution} on $X$ is a saturated subsheaf of $T_X$.
  Given a distribution $\sD$, its \emph{Cauchy characteristic
    distribution}\index{Cauchy characteristic of a distribution} is the integrable subdistribution
  $\Ch(\sD)\subset\sD$ whose fiber at the general point $x\in X$ is
  $$
  \Ch(\sD)_x:=\left\{v\in\sD_x;\,N(v,\sD_x)=0\right\},
  $$
  where $N$ denotes the O'Neill tensor \index{O'Neill tensor}, i.e., the $\O_X$-linear map
  $N : \sD\otimes\sD\rightarrow T_X/\sD$ induced by the Lie bracket.
\end{definition}

On $\sC \subset \P(T_X^\vee)$ we can consider a natural distribution
$\sP$ defined in the general point $\alpha\in\sC$ by:
$$
\sP_{\alpha} := (T\rho)^{-1} (T_{\hat{\sC}_{x,\alpha}}),
$$
where $T_{\hat{\sC}_{x,\alpha}} \subset T_X|_x$ is the tangent
space to the affine cone of $\sC_x$ along the ray $\C \cdot \alpha
\subset T_X|_x$ determined by $\alpha$.  The following proposition is
the technical core of Theorem \ref{thm:birtangmap}. It is based in a
detailed study of the distribution $\sP$ and its relation with the
family $H$, which is beyond the purpose of this survey. We refer to
\cite{HM04} for a detailed account.

\begin{proposition}\label{lem:Cauchysub}
  With notation as above, the following holds.
  \begin{enumerate}
  \item\label{it:contcurv} The distribution $\Ch(\sP) \subset
    T_{\mathcal C}$ contains the tangent directions of the images in
    $\sC$ of the curves parametrized by $H$. More precisely, if
    $[\ell] \in H$ is a general point, then the morphism
    $\tau|_{\pi^{-1}([\ell])} : \P^1 \to \mathcal C$ is immersive and
    its image is tangent to $\Ch(\sP)$.
  \item\label{it:leaf} Let $y \in \mathcal C$ be a general point and
    $S$ the associated leaf of $\Ch(\sP)$. Then $S$ is algebraic, and
    there exists a dense open set $S_0$ such that $W := \pi(S_0)
    \subset X$ is quasi-projective, and such that 
    $$
    \pi|_{S_0}:S_0\rightarrow W
    $$
    is a bundle of projective spaces, isomorphic to
    $\P(T_{W}^\vee)$.
  \item\label{it:family} Let ${\overline W}\subset X$ be the Zariski
    closure of $W$.  The subvariety 
    $$
    H_W := \left\{ [\ell] \in H \, \bigl| \, \ell \subset \overline
      W \right\}
    $$
    is a dominating family of rational curves of minimal
    degrees\footnote{Strictly speaking, we have defined the notion of
      a family of rational curves only for normal varieties because
      the notion of the Chow-variety and its universal property is a
      little delicate for non-normal spaces. Although $\overline W$
      need not be normal, we ignore this (slight) complication in this
      sketch for simplicity.} in $\overline{W}$.  \qed
  \end{enumerate}
\end{proposition}
The subschemes of the form $\overline{W}$ are called \emph{Cauchy
  subvarieties of}\index{Cauchy subvariety} $X$ with respect to $H$.

\subsection*{Step 2: End of sketch of proof}

Let $x\in X$ be a general point and $c \in \sC_x \subset
\P(T_X|_x^\vee)$ be a general minimal rational tangent. We assume to
the contrary and take two general curves $[\ell_1]$ and $[\ell_2] \in
H_x$ that have $c$ as tangent direction. Let $S$ be the leaf of
$\Ch(\sP)$ that contains $c$ and let ${\overline W} \subset X$ be the
corresponding Cauchy\index{Cauchy subvariety} subvariety. By
Proposition~\ref{lem:Cauchysub}.(1), ${\overline W}$ contains both
$\ell_1$ and $\ell_2$.

By Proposition~\ref{lem:Cauchysub}.(2) and (3), the tangent map
associated with $H_W$,
$$
\tau_{W,x} : H_{W,x} \dasharrow \P(T_{\overline W}|_x^\vee) \subset
\P(T_X|_x^\vee),
$$
is surjective. Theorem~\ref{thm:CMSB02} then applies to the
normalization of $\overline W$ and asserts that $\tau_W$ must be
birational. By general choice, $\ell_1$ and $\ell_2$ must then be
equal, a contradiction.\qed
\index{Tangent morphism!birationality|)}
\index{Tangent morphism|)}

\section{The importance of the VMRT}
\label{sec:importance}\index{Variety of minimal rational tangents!importance}

Given a uniruled variety $X$ and a dominating family $H \subset
\rat(X)$ of rational curves of minimal degrees, we have claimed that
geometry of $X$ is determined to a large degree by the projective
geometry of the associated VMRT $\sC \subset \P(T_X^\vee)$. In this
section we would like to name a few results that support the claim.
Section~\ref{sec:secants} discusses one particular example in more
detail.

In view of the Minimal Model Program, we restrict ourselves to Fano
manifolds of Picard number $1$. In this case, Hwang and Mok\index{Hwang, Jun-Muk}\index{Mok, Ngaiming} have shown
under some technical assumptions, that $X$ is \emph{completely
  determined} by the family of minimal rational tangents over an
analytic open set. This ``Cartan-Fubini'' type result\index{Cartan-Fubini type theorem}\index{Cartan, \'Elie}\index{Fubini, Guido} is stated as follows.

\begin{theorem}[\protect{\cite{HM01}}]\label{thm:cartanfubini}
  Let $X$ and $X'$ be Fano manifolds of Picard number $1$ defined over
  the field of the complex numbers, and let $H$ and $H'$ be
  dominating families of minimal rational curves in $X$ and $X'$,
  respectively. Let $\sC_x$ and $\sC'_{x'}$ denote the associated
  VMRT\index{VMRT} at $x\in X$ and $x'\in X'$.
  Assume that $\sC_x$ is positive-dimensional and that the Gauss map
  of the embedding $\sC_x\subset\P(\Omega_X|_x)$ is finite for general
  $x\in X$.  Then any biholomorphic map $\phi:V\rightarrow V'$ between
  analytic open neighborhoods $V\subset X$ and $V'\subset X'$ inducing
  an isomorphism between $\sC_x$ and $\sC'_{\phi(x)}$ for all $x\in V$
  can be extended to a biholomorphic map $\Phi:X\rightarrow X'$. \qed
\end{theorem}

A detailed proof of Theorem~\ref{thm:cartanfubini} is given in the
survey article \cite{Hwa00}. There it is also discussed under what
conditions even stronger results can be expected.

The VMRT\index{VMRT} has also been used to attack the following problems.

\begin{description} \index{Variety of minimal rational tangents!example applications}
\item[\emph{Stability of the tangent bundle:}] For Fano manifolds of
  Picard number $1$, this property can be easily restated in terms of
  projective properties of $\sC_x\subset\P(T_X|_x^\vee)$ for general
  $x$, which can be checked in some cases, see e.g.~\cite{Hwa98, HW01}
  for low-dimensional Fano manifolds and moduli spaces of vector
  bundles\index{Moduli space of vector bundles on a curve}, or \cite{Kebekus05} for contact manifolds\index{Contact manifold}.
\item[\emph{Deformation rigidity:}]\index{Deformation rigidity!of certain varieties} The VMRT\index{VMRT} have been used by Hwang \index{Hwang, Jun-Muk}
  and Mok\index{Mok, Ngaiming} to prove deformation rigidity of various types of varieties
  and morphisms. See, e.g., \cite{Hwa97, Hwa00} or \cite{HM98}.
\item[\emph{Uniqueness of contact structures:}]\index{Contact manifold!uniqueness of contact structure} It has been shown in
  \cite{Kebekus01} that contact structures on Fano manifolds of Picard
  number $1$ are unique since their projectivization at a point
  coincides with the linear span of the VMRT\index{VMRT}.
\item[\emph{The Remmert--Van de Ven / Lazarsfeld problem:}]\index{Remmert--Van de Ven / Lazarsfeld problem} Theorem
  \ref{thm:cartanfubini} can be used to classify smooth images of
  surjective morphisms from rational homogeneous spaces of Picard
  number $1$, see \cite{HM99}.
\end{description}

\section{Higher secants and the length of a uniruled manifold}\index{Length!of a uniruled manifold|(}
\label{sec:secants}

As before, let $X$ be a complex projective manifold and $H \subset
\rat(X)$ be a dominating family of rational curves of minimal degrees.
If $b_2(X) = 1$, then $X$ is rationally connected ---see
Definition~\ref{def:ratconn} on page~\pageref{def:ratconn}. In
particular, if $x$, $y \in X$ are any two general points, there exists
a connected chain of rational curves $\ell_1, \ell_2, \ldots, \ell_k
\in H$ such that $x \in \ell_1$, $y \in \ell_k$ ---see
figure~\ref{fig:rcc} on page \pageref{fig:rcc}.

The number $k$, i.e., the minimal length of chains of $H$-curves
needed to connect two general points is called the \emph{length of $X$
  with respect to $H$}. The length is an important invariant; among
other applications, it was used by Nadel\index{Nadel, Alan M.} in the proof of the
boundedness of the degrees of Fano manifolds of Picard number one,
\cite{Nadel91}.

The main aim of this section is to introduce an effective method that
allows to compute the length for a number of interesting varieties. We
relate the length of $X$ to the projective geometry of the variety of
minimal rational tangents. The length of $X$ can then be computed in
situations where the secant defect of the VMRT is known.\index{Variety of minimal rational tangents!secant defect}

We will employ these results in Chapter~\ref{chap:moduli} below to give
a bound on the multiplicities of divisors at a general point of the
moduli of stable bundles of rank two on a curve.

\subsection{Statement of result}

To formulate the result precisely, let us fix our assumptions first.
For the remainder of the present section, let $X$ be a complex
projective manifold of arbitrary Picard number and $H \subset \rat(X)$
a dominating family of rational curves of minimal degrees, as in
Section~\ref{sec:setupgeom}. If $x \in X$ is a general point, assume
additionally that the space $H_x$ is irreducible\footnote{The
  irreducibility assumption is posed for simplicity of exposition. It
  holds for all examples we will encounter. See \cite[thm.~5.1]{KK04}
  for a general irreducibility criterion.}.

We also need to consider a few following auxiliary spaces.

\begin{definition}\label{def:sectoc}
  Consider the irreducible subvarieties\index{$\loc^k(x)$, locus of chains of rational curves through a point}
  \begin{align*}
    \loc^1(x) & :=\mbox{closure of}\bigcup_{\mbox{\scriptsize general
        $C\in H_x$}}C, & \mbox{ and } \\
    \loc^{k+1}(x) &:=\mbox{closure of}
    \bigcup_{\begin{array}{c} \mbox{\scriptsize general $C\in H_y$}\\
         \mbox{\scriptsize general $y\in\loc^k(x)$}\end{array}}C.
  \end{align*}
  Set $d_k:=\dim(\loc^k(x))$; we call it the $k$-th \emph{spanning
    dimension}\index{Family of rational curves!$k$th-spanning dimension of} of the family $H$. Finally, let $\sC^k_x \subset
  \P(T_X|_x^\vee)$ be the tangent cone to $\loc^k(x)$ at $x$.
\end{definition}

\begin{remark}\label{rem:fail_of_strict}
  Even though $\sC_x$ is assumed to be irreducible, $\sC^1_x$ can have
  several components. This might be the case if some of this curves
  associated with points in $H_x$ are nodal, see
  Figure~\ref{fig:node}.
  \begin{figure}[tbp]
    \centering
    \begin{picture}(6,3.5)(-1,0)
      \put(0, 0){\includegraphics[height=3.5cm]{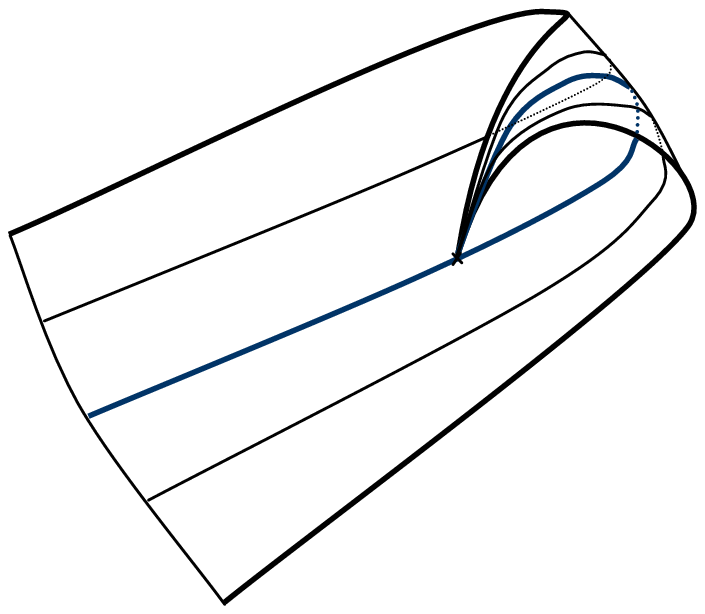}}
      \put(2.6,1.8){$x$}
      \put(-0.45,1.2){\tiny nodal}
      \put(0.2,1.25){\vector(1,0){0.5}}
      \put(0.05,0.2){\tiny smooth}
      \put(0.85,0.25){\vector(1,0){0.6}}
      \put(0.85,0.25){\vector(1,1){0.5}}
      \put(-0.8,1.75){\tiny smooth}
      \put(0,1.8){\vector(1,0){0.45}}
      \put(0,1.8){\vector(1,1){0.5}}
    \end{picture}
    \caption{Nodal curves causing reducibility of $\sC_x^1$}
    \label{fig:node}
  \end{figure}
\end{remark}

It is clear that the spanning dimensions of $H$ do not depend on the
choice of the general point $x$. The tangent cones $\sC^k_x \subset
\P(T_X|_x^\vee)$ have pure dimension $\dim \sC^k_x = d_k-1$ ---we
refer to \cite[lect.~20]{Harris95} for an elementary introduction to
tangent cones.

\begin{remark}
  The spanning dimensions are natural invariants of the pair $(X, H)$.
  In many situations, however, there is a canonical choice of $H$, so
  that we can actually view the $d_k$ as invariants of $X$.
\end{remark}
  
\begin{definition}
  If there exists a number $l > 0$ such that $d_l = \dim(X)$ and
  $d_{l-1} < \dim(X)$, we call $X$ \emph{rationally connected by
    $H$-curves}, write $\length_{H}(X) = l$ and say that $X$ has
  \emph{length} $l$ with respect to $H$.\index{Length!of a uniruled manifold}
\end{definition}

\begin{remark}[\protect{\cite{Nadel91}}]
  If there exists a number $k$ such that $d_k = d_{k+1}$, it is clear
  from the definition that $\loc^{k}(x) = \loc^{k+1}(x)= \loc^{k+2}(x)
  = \cdots$. In particular, if $X$ is rationally connected by
  $H$-curves, then $\length_H(X) \leq \dim X$.
\end{remark}

\begin{definition}[\protect{\cite{CC02, Zak93}}]\index{Variety of minimal rational tangents!secant variety and spanning dimension}\index{Variety of minimal rational tangents!secant variety and spanning dimension!for varieties covered by lines}
  Let $S^k\sC_x$ be the $k$-th secant variety of $\sC_x$, i.e., the
  closure of the union of $k$-dimensional linear subspaces of
  $\P(T_X|_x^\vee)$ determined by general $k+1$ points on $Y$.
\end{definition}

With this notation, the main result is formulated as follows.

\begin{theorem}[\protect{\cite[thms.~3.11--14]{HK05}}]\label{thm:chainsratcurves}
  Let $X$ be as above. Then the following holds.
  \begin{enumerate}
  \item We have $\sC_x\subset \sC_x^1 \subset S^1\sC_x$. If none of
    the curves of $H$ has nodal singularities at $x$ then
    $\sC_x=\sC^1_x$.
  \item For each $k\geq 1$, we have $S^k\sC_x \subset\sC^{k+1}_x$. In
    particular, $d_{k+1} \geq \dim(S^k\sC_x)+1$.
  \item If $X$ admits an embedding $X \hookrightarrow \P^N$ such that
    the curves parametrized by $H$ are mapped to lines, then
    $S^1\sC_x=\sC^2_x$, and so $d_2=\dim(S^1\sC_x)+1$.
  \end{enumerate}
\end{theorem}

We will give a rough sketch of a proof below. In view of
Theorem~\ref{thm:chainsratcurves}, one might be tempted to conjecture
that $C^{k+1}_x=S^k\sC_x$. Remark~\ref{rem:fail_of_strict} and the
following remark show that this is not the case.

\begin{remark}\label{rem:fail_of_strict2}
  The inequality in Theorem~\ref{thm:chainsratcurves}.(2) can be
  strict.  For an example, let $X \not \cong \P^n$ be a Fano manifold
  with Picard number one which carries a complex contact\index{Contact manifold}
  structure\footnote{See e.g.~\cite{Beauville99, KPSW00} and the
    references therein for an introduction to complex contact
    manifolds.}, and $H$ a dominating family of rational curves of
  minimal degrees. In this setup, $\sC_x$ is linearly degenerate. In
  particular, we have that
  $$
  \dim(S^k\sC_x) + 1 < \dim X = \length_{H}(X)
  $$
  for all $k$. 
\end{remark}

Generalizing Nadel's product theorem\index{Nadel, Alan M.}\index{Nadel, Alan M.!product theorem}, \cite{Nadel91}, the knowledge of
the length can be used to bound multiplicities of sections of vector
bundles at general points of $X$, as follows.

\begin{proposition}[\protect{\cite[prop.~2.6]{HK05}}]\label{prop:multbound}\index\index{Nadel, Alan M.!product theorem!generalization}
  Assume that $X$ is rationally connected by curves of the family $H$ and let $V$ be a
  vector bundle on $X$. Consider a general curve $C \in H$, let
  $\nu : \P^1 \to C$ be its normalization, and write
  $$
  \nu^*(V) \cong\O(a_1)\oplus\cdots\oplus\O(a_r),\quad \text{with }a_1 \geq
  \cdots \geq a_r.
  $$
  
  If $x \in X$ is a general point and $\sigma\in H^0(X,V)$ any
  non-zero section, then the order of vanishing of $\sigma$ at $x$
  satisfies
  $\mult_x(\sigma)\leq \length_{H}(X)\cdot a_1$. \qed
\end{proposition}

\subsection{Sketch of Proof of Theorem~\ref{thm:chainsratcurves}.(1)}

Let $x\in X$ be a general point. The first inclusion,
$\sC_x\subset\sC_x^1$, is immediate from the definition.

For the other inclusion, let again $\tilde H_x$ be the normalization
of $H_x$ and consider the pull-back of Diagram~\eqref{eq:d2} from
page~\pageref{eq:d2},
$$
\xymatrix{
  U_x \ar[d]_{\pi_x} \ar[rr]^{\iota_x}_{\text{evaluation}} && X \\
  \tilde H_x }
$$
If none of the curves of $H_x$ has nodal singularities at $x$,
Theorem~\ref{thm:dimsingcurves} asserts that the preimage
$\iota_x^{-1}(x)$ is exactly a reduced section $\sigma_\infty \cong
\tilde H_x$. By the universal property of blowing up, the evaluation
morphism $\iota_x$ factors via the blow-up of $X$ at $x$ and
the equality $\sC_x^1=\sC_x$ follows.
  
If some of the curves in $H_x$ are nodal, $\iota_x^{-1}(x)$ contains
the section $\sigma_\infty$ and finitely many points. A somewhat
technical analysis of the tangent morphism $T \iota_x$, and a
comparison between positive directions in the restriction of $T_X$ to
a nodal curve and the Zariski tangent space to the VMRT\index{VMRT} then yields
the inclusion $\sC_x^1\subset S^1\sC_x$.

\subsection{Sketch of Proof of Theorem~\ref{thm:chainsratcurves}.(2)}

Here we only give an idea of why the first secant variety of the VMRT\index{VMRT},
$S^1\sC_x$ is contained in the tangent cone $\sC^2_x$ to the locus of
length-2 chains of rational curves.

If $v$ and $w \in \sC_x$ are any two general elements, we need to show
that the line in $\P(T_X|_x^\vee)$ through $v$ and $w$ is contained in
the tangent cone to $\loc^2(x)$. To this end, let $[\ell_v]$ and
$[\ell_w] \in H_x$ be two rational curves that have $v$ and $w$ as
tangent directions, respectively. By general choice of $v$ and $w$,
the curves $\ell_v$ and $\ell_w$ are smooth at $x$.

\begin{figure}[tbp]
  \centering
  \begin{picture}(8,4.2)(0,0)
    \put(0, 0){\includegraphics[height=4cm]{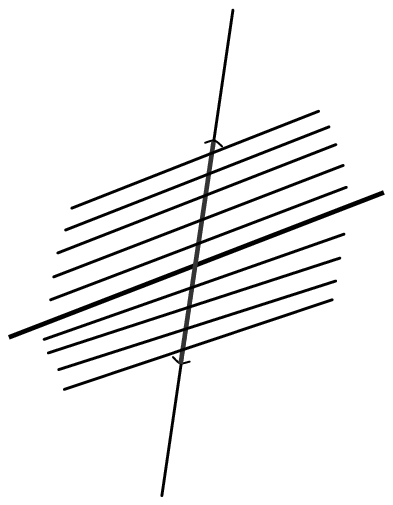}}
    \put(1.5,3.8){${\scriptstyle l_v}$}
    \put(0,1.5){${\scriptstyle l_w}$}
    \put(1.1,2.9){${\scriptstyle \Delta}$}
    \put(1.25,2.8){\vector(1,-2){0.3}}
    \put(2.5,3.1){${\scriptstyle S\subset\loc^2(x)}$}
    \put(1.33,1.89){${\scriptstyle x}$}
    \put(0.7,0.2 ){$X$}
    \put(6, 0){\includegraphics[height=4cm]{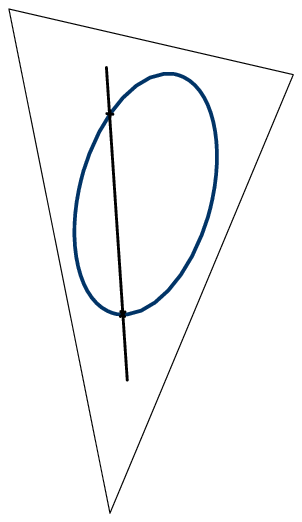}}     
    \put(5.5,0.2){$\P(T_X|_x^\vee)$}
    \put(6.2,3.5){${\scriptstyle \P(T_S|_x^\vee)}$}
    \put(7.72,3){${\scriptstyle \sC^1_x}$}
    \put(6.9,3){${\scriptstyle v}$}
    \put(7,1.7){${\scriptstyle w}$}
  \end{picture}
  \caption{Secants to $\sC^1_x$}                            
  \label{fig:tancone}
\end{figure}

Consider a small unit disk $\Delta \subset \ell_v$, centered about
$x$. By general choice of $x$, we can find a smooth holomorphic arc
$\gamma: \Delta \to H$ such that
\begin{itemize}
\item $\gamma(0) = [\ell_w]$, 
\item for any point $y \in \Delta \subset \ell_v \subset X$, the curve
  $\gamma(y)$ contains $y$, i.e., $\gamma(y) \in H_y$.
\end{itemize}
The situation is described in Figure~\ref{fig:tancone} above. The
unions of the curves $\gamma(\Delta)$ then forms a surface $S \subset
\loc^2(x)$. By general choice, it can be seen that $S$ is smooth at
$x$ and contains both $v$ and $w$ are tangent directions ---see
Figure~\ref{fig:tancone}. As a consequence, we obtain that the line in
$\P(T_X|_x^\vee)$ through $v$ and $w$ is also contained in
$\P(T_S|_x^\vee)$. Since $\P(T_S|_x^\vee) \subset \sC^2_x$, this shows
the claim.

\subsection{Sketch of Proof of Theorem~\ref{thm:chainsratcurves}.(3)}
\label{sec:chains}

Let $v \in \sC^2_x$ be any element in the tangent cone to $\loc^2(x)$.
We need to show that $v$ is contained in the first secant variety to
$\sC_x$, i.e., that $v \in S^1\sC_x$.

Recall from \cite[lect.~20]{Harris95} that the tangent cone $\sC^2_x$
to $\loc^2(x)$ is set-theoretically exactly the union of the tangent
lines to holomorphic arcs $\gamma : \Delta \to \loc^2(x)$ that are
centered about $x$. Recall also that if the arc $\gamma$ is not smooth
at 0, then the tangent line at 0 is the limit of the tangent lines to
$\gamma$ at points where $\gamma$ is smooth.

\begin{figure}[tbp]
  \centering
  \begin{picture}(12,5)(0,0)
    \put(0, 0.5){\includegraphics[height=3.5cm]{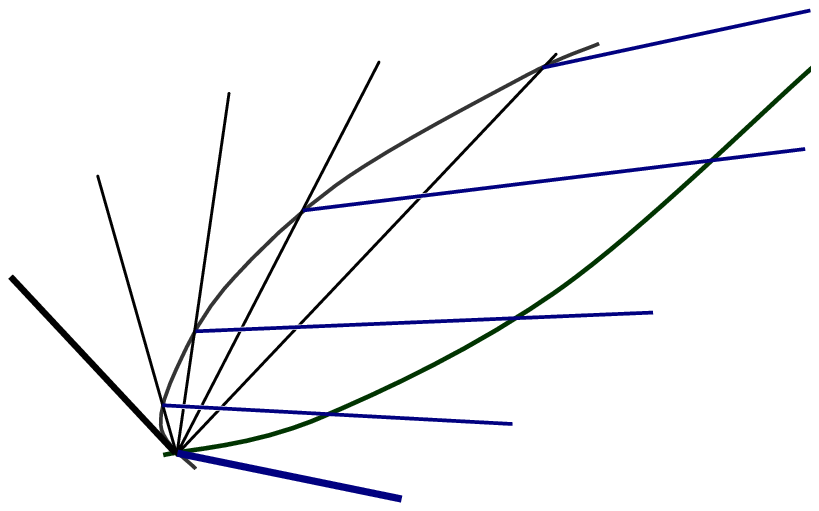}}
    \put(0.45,0.7){${\scriptstyle \gamma(0)=x}$}
    \put(1.65,2.53){${\scriptstyle g(t)}$}
    \put(2.5,3.63){${\scriptstyle G(t)}$}
    \put(5.53,2.9){${\scriptstyle F(t)}$}
    \put(0,2.15){${\scriptstyle G(0)}$}
    \put(2.8,0.4){${\scriptstyle F(0)}$}
    \put(4.5,3){${\scriptstyle \gamma(t)}$}
    \put(2.5,0){\tiny Case (i)}
    \put(6.5,0.5){\includegraphics[height=3.5cm]{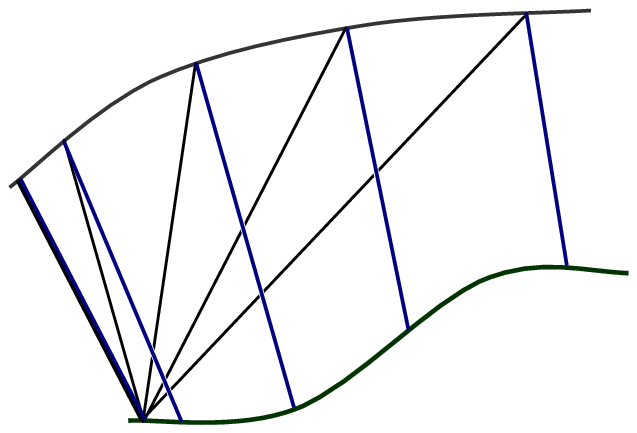}}
    \put(6.8,0.4){${\scriptstyle \gamma(0)=x}$}
    \put(9.1,3.85){${\scriptstyle g(t)}$}
    \put(8.78,2.6){${\scriptstyle G(t)}$}
    \put(6.6,1.5){${\scriptstyle G(0)}$}
    \put(6.55,1.2){${\scriptstyle =F(0)}$}
    \put(9.7,2){${\scriptstyle F(t)}$}
    \put(9.7,1.1){${\scriptstyle \gamma(t)}$}
    \put(5.9,2.3){${\scriptstyle g(0)=y}$}
    \put(8.5,0){\tiny Case (ii)} 
  \end{picture}
  \caption{Arcs in $\loc^2(x)$}
  \label{fig:arcsloc2}
\end{figure}

We can thus find an arc $\gamma: \Delta \to \loc^2(x)$ with
$\gamma(0)=x$ and tangent $v$. Recall that $\loc^2(x)$ is the locus of
chains of rational curves of length $2$ that contain $x$. Replacing
$\Delta$ by a finite covering, if necessary, we can then find arcs
$$
G:  \,\Delta \to H_x \quad \text{ and } \quad F:  \,\Delta \to  H
$$
such that for all $t \in \Delta$ the curves $G(t) \cup F(t)$ form
a chain of length two that contains both $x$ and $\gamma(t)$ ---this
construction is depicted in Figure~\ref{fig:arcsloc2}. Further, we
consider arcs, defined for general $t$ as follows.
\begin{align*}
    g: & \,\Delta \to \loc^1(x) & P: & \,\Delta \to  {\rm Grass}(2,\P^N) \\
    & \,\, t \,\, \to  G(t)\cap F(t) & & \,\, t \,\, \to \text{plane spanned by $G(t)$ and $F(t)$} 
\end{align*}
Observe that this makes $g$ and $P$ well-defined because the target
varieties are proper.

We distinguish two cases.

\begin{description}
\item[Case (i), $G(0)\neq F(0)$:] In this case, $x=g(0)=G(0)\cap
  F(0)$, and $[v] \in \P(T_{P(0)}|_x^\vee)$. Since for general $t \in
  \Delta$, the line $\P(T_{P(t)}|_x^\vee)$ is contained in
  the secant variety $S^1\sC_{g(t)}$, we have $\P(T_{P(0)}|_x^\vee)
  \subset S^1\sC_x$. This shows the claim.
  
\item[Case (ii), $G(0)=F(0)$:] Since curves of $H$ are lines by
  assumption, Proposition~\ref{prop:Cxforlines} applies to the point
  $y=g(0)$ and the line $\ell = G(0)$.

  Since $\P(T_{P(t)}|_{g(t)}^\vee)$ is secant to $\sC_{g(t)}$ for
  general $t$, we obtain that $\P(T_{P(0)}|_{y}^\vee)$ is secant to
  $\sC_y$. Moreover, since $\ell = G(0) = F(0)$, the line
  $\P(T_{P(0)}|_{y}^\vee)$ lies in projective tangent space to $\sC_y$
  at the point corresponding to $[\ell]$.
  
  To end, observe that the plane $P(0)$ is tangent to $X$ all along
  $\ell$. Its tangent directions determine a subbundle of
  $T_X|_{\ell}^+$. This implies $\P(T_{P(0)}|_{x}^\vee) \subset
  \P(T_X|_{\ell}^+|_x^\vee)$ and concludes the proof. \qed
\end{description}

\subsection{Open Problems}

\begin{problem}
  In view of Remarks~\ref{rem:fail_of_strict} and
  \ref{rem:fail_of_strict2}, is it possible to improve on
  Theorem~\ref{thm:chainsratcurves} and find a formula that computes
  the spanning dimensions in all cases? Perhaps one needs to take into
  account how the VMRT\index{VMRT} (or its linear span) deforms as the base point
  changes.
\end{problem}
\index{Length!of a uniruled manifold|)}
\index{Variety of minimal rational tangents|)}

% Local Variables:
% TeX-master: "RC-arXive.tex"
% End:

\chapter{Examples of uniruled varieties: moduli spaces of vector bundles}\index{Moduli space of vector bundles on a curve|(}
\label{chap:moduli}

In this chapter we apply the results of Chapter~\ref{chap:geomtry} to
one particular example of a uniruled variety, namely the moduli space
of stable vector bundles on a curve\index{Moduli space of vector bundles on a curve}. We first recall a classic
construction of rational curves on the moduli space, and then show by
example how the understanding of the VMRT\index{Variety of minimal rational tangents} can be used to study
geometric properties of the moduli space.

\section{Setup of notation. Definition of Hecke curve}\index{Hecke curves}
\label{sec:Heckedef}

Let $C$ be a smooth projective curve of genus $g\geq 3$ and $L$ be a
line bundle of $C$ of degree $d$. We will denote by $M^r(L)$, or
simply by $M^r$\index{Moduli space of vector bundles on a curve} if there is no possible confusion, the moduli scheme
of semistable vector bundles of rank $r$ and determinant $L$. We can
assume without loss of generality that $0 \leq d \leq r-1$, otherwise
twist with a line bundle.  For simplicity, let us assume that
$(r,d)=1$, in which case $M^r$ is a smooth Fano manifold\index{Moduli space of vector bundles on a curve!is sometimes Fano} of Picard
number one. We remark that most of the arguments and constructions
presented here work the same for the general case.

\begin{definition}
  Given two integers $k$ and $l$, we say that a vector bundle $E$ of
  rank $r$ is \emph{$(k,l)$-stable}\index{$(k,l)$-stable sheaf} if every locally free subsheaf $F
  \subset E$ verifies:
  $$
  \frac{\deg(F)+k}{\rank(F)}<\frac{\deg E+k-l}{\rank(E)}.
  $$
  Let $M^r(k,l)\subset M^r$ the subset parameterizing
  $(k,l)$-stable vector bundles of rank $r$.
\end{definition}
\begin{remark}
  The subsets $M^r(k,l)$ are open in $M^r$. For instance
  $M^r(0,0)\subset M^r$ is the open set parametrizing stable bundles.
\end{remark}

Narasimhan\index{Narasimhan, Mudumbai S.} and Ramanan\index{Ramanan, S.} showed that the subset $M^r(1,1)$ is non-empty
for $g\geq 3$ ---see \cite[prop.~5.4]{NR78} for a precise statement.
For every $(1,1)$-stable bundle $E$ they constructed a $1$-dimensional
family deformations of $E$ that corresponds to a rational curve in
$M^r$. The so-constructed curves are called \emph{Hecke curves}\index{Hecke curves}. We
briefly review their definition below.

\subsection*{Construction of the Hecke curves}
Let $E$ be a vector bundle of rank $r$ and let $x$ be any point of
$C$. Every element $p=[\alpha] \in \P(E_x)$ provides an exact sequence
of $\O_C$-modules,
\begin{equation}\label{eq:elemtrans}
  0\rightarrow E^p\stackrel{\phi}{\longrightarrow} E\longrightarrow\O_x\rightarrow 0,
\end{equation}
which is induced by the map $\alpha : E_x \to \C$. The kernel $E^p$
depends only on the class $p=[\alpha]$ and it is called the {\it
  elementary transform}\index{Elementary transform} of $E$ at $p$. Let $H_p(E) := (E^p)^\vee$ be
its dual. 

\begin{figure}[tbp]
  \centering
  \begin{picture}(12,5.5)(0,0)
    \put(0, 0.4){\includegraphics[height=2.3cm]{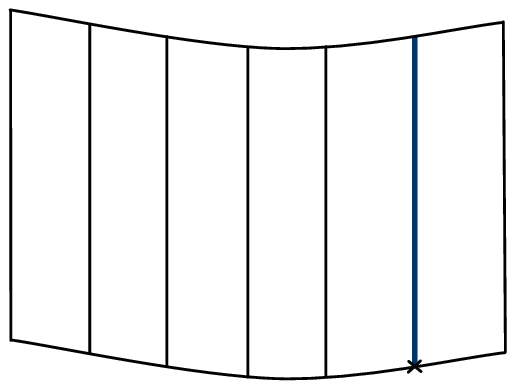}}
    \put(3.6, 2.8){\includegraphics[height=2.7cm]{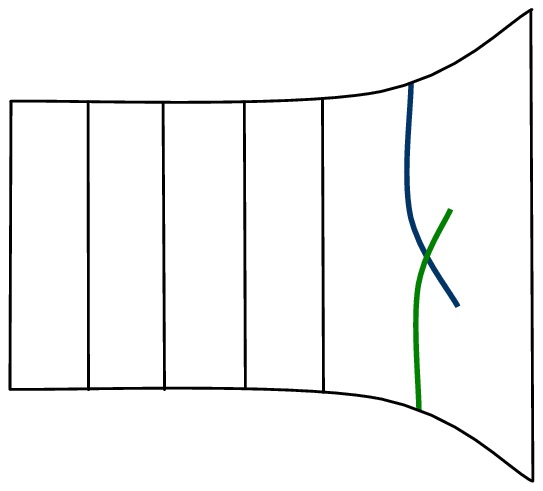}}
    \put(8, 0){\includegraphics[height=2.7cm]{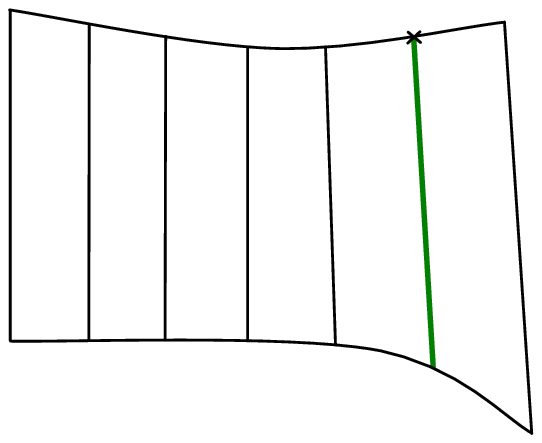}}  
    \put(2.5,0.3){\tiny $p$}
    \put(10.5,2.6){\tiny $p'$}
    \put(3.1,0.2){$\P(E)$}
    \put(7.2,0.2){$\P(E_p)$}
    \put(3.2,2.2){\tiny $\P(E_x)$}
    \put(3.15,2.2){\line(-2,-1){0.63}}
    \put(3.5,3.7){\vector(-3,-2){1.7}}
    \put(1.4,3.25){\tiny blowing-up $p$}
    \put(6.8,3.9){\vector(3,-2){1.8}}
    \put(7.9,3.3){\tiny blowing-down $\widetilde{\P(E_x)}$}
    \put(5.95,4.6){\tiny $\widetilde{\P(E_x)}$}
    \put(6,3.5){\tiny Exc.}
  \end{picture}
  \caption{Elementary transform of a ruled surface. Moving the point $p$ along  $\P(E_x)$ gives a rational curve in $M^r$.}
  \label{fig:elemtrans}
\end{figure}

Considering the associated projective bundles, elementary transforms
can be described in terms of blow-ups and blow-downs. We have depicted
the case $r=2$ in Figure~\ref{fig:elemtrans}.

Dualizing the exact sequence~\eqref{eq:elemtrans} we obtain the
following.
\begin{equation}\label{eq:elemtransdual}
  0\rightarrow E^\vee\stackrel{\phi^\vee}{\longrightarrow} H_p(E)\longrightarrow\O_x\rightarrow 0
\end{equation}
Sequence~\eqref{eq:elemtransdual} expresses $E^\vee$ as the elementary
transformation of $H_p(E)$ at the point
$p'=\coker\phi^\vee_x\in\P(H_p(E)_x)$, that is $H_{p'}(H_p(E))=E$.  We
can then consider the $H_q(H_p(E))$ as deformations of $E$
parametrized by $q\in\P(H_p(E)_x)$.

It is easy to see that if $E$ is $(k,l)$-stable, then $H_p(E)$ is
$(l-1,k)$-stable ---see \cite[lem.~5.5]{NR78}. In particular, when $E$
is $(1,1)$-stable, then every element of the form $H_q(H_p(E))$ is
stable. A line $L\subset\P(H_p(E)_x)$ through $p'$ will therefore
determine a rational curve through $[E]$ in $M^r$. We denote this
curve by $C(E,p,L)$.

\begin{definition}\label{def:Heckedefn}\index{Hecke curves}
  A curve of the form $C(E,p,L)$ constructed as above is called a
  \emph{Hecke curve} through $[E]$.
\end{definition}

\section{Minimality of the Hecke curves}\label{sec:Heckearemin}\index{Hecke curves!have minimal degree}

Notice that lines in $\P(H_p(E)_x)$ through $p'$ are in one-to-one
correspondence with the set of hyperplanes in $\P(E_x)$ that contain
$p$, i.e., with $\P(T_{\P(E_x)|_p})$. The next proposition states
basic properties of Hecke curves that have been originally proved by
Narasimhan\index{Narasimhan, Mudumbai S.} and Ramanan\index{Ramanan, S.}---we refer the interested reader to
\cite[sect.~5]{NR78} for details.
\begin{proposition}\label{prop:Hecke} 
  With the same notation as above, the following holds.
  \begin{itemize}
  \item \emph{\cite[5.13]{NR78}}: The map that sends $(p,L)$ to the
    Hecke curve $C(E,p,L)$ is injective for every $[E] \in M^r(1,1)$.
    In particular, Hecke curves through $[E]$ are naturally
    parametrized by the points of the projectivization of the relative
    tangent bundle $\P(T_{\P(E),C})$.
  \item \emph{\cite[5.9, 5.15, 5.16]{NR78}}: Hecke curves\index{Hecke curves!are free and smooth}\index{Hecke curves!are rational} are free
    smooth rational curves of anticanonical degree $2r$ ---see
    Proposition~\ref{prop:gen=stand} on page \pageref{prop:gen=stand}
    for the notion of a free rational curve.  \qed
  \end{itemize}
\end{proposition}

\begin{definition}
  Set $H_E := \P(T_{\P(E),C})$, and let $\O(1)$ be the associated
  tautological line bundle. Further, let $\rho : H_E \to C$ be the
  natural projection.  
  
  If $[E]\in M^r(1,1)$ is any point, we can view $H_E$ as a subset of
  $\rat(M^r)$. Let $H \subset \rat(M^r)$ be the closure of the union
  of the $H_E$, for all $[E]\in M^r(1,1)$.
\end{definition}

Minimality of the Hecke\index{Hecke curves!have minimal degree} curves was shown by Hwang\index{Hwang, Jun-Muk} in
\cite[prop.~9]{HW00} for the case $r=2$. Recently Sun\index{Sun, Xiaotao} has proven that
Hecke curves are minimal curves in $M^r$ for any $r$, and that the
converse is also true:

\begin{theorem}[\protect{\cite[thm.~1]{Sun05}}]\label{them:Heckeminuni}\index{Moduli space of vector bundles on a curve!rational curves of minimal degree}
  Any rational curve $C\subset M^r$ passing through a general point
  of $M^r$ has anticanonical degree at least $2r$. If $g\geq 3$, then the
  anticanonical degree of $C$ is $2r$ if and only if it is a Hecke
  curve. \qed
\end{theorem}

In particular, Sun\index{Sun, Xiaotao} shows that the variety $H$ is a dominating family
of rational curves of minimal degrees in $M^r$\index{Hecke curves!form a dominating family}\index{Hecke curves!have minimal degree}.  If $[E]\in M^r(1,1)$
is a general point, the associated subspace of curves through $[E]$ is
exactly $H_E$.

\section{VMRT associated to $M^r$}\label{sec:VMRTmoduli}\index{Moduli space of vector bundles on a curve!variety of minimal rational tangents}

The variety of minimal rational tangents to $M^r$ at $[E]\in M^r(1,1)$
is the image of the tangent map\index{Moduli space of vector bundles on a curve!tangent map for Hecke curves}:
$$
\begin{array}{rccc}
  \tau_E: & H_E & \longrightarrow&\P(T^\vee_{M^r,[E]})\\
  &(p,L)&\mapsto&[T_{C(E,p,L),[E]}] 
\end{array}
$$

We would like to  understand $\tau_E$ in two ways, namely
\begin{itemize}
\item in terms of the Kodaira-Spencer map\index{Kodaira-Spencer map} associated to Hecke curves,
  as studied by Narasimhan\index{Narasimhan, Mudumbai S.} and Ramanan\index{Ramanan, S.} in \cite{NR78}, and
\item in terms of linear systems on $H_E$.
\end{itemize}
For the first item, recall the standard description of the tangent
space to the moduli scheme,
$$
T_{M^r,[E]} \cong H^1 \bigl(C,\ad(E)\bigr) \cong
H^0\bigl(C,K_C\otimes\rho_*(\O(1))\bigr),
$$
where $\ad(E)$ is the sheaf of traceless endomorphisms of $E$. The
morphism $\tau_E$ is then described as follows.

\begin{proposition}[\protect{\cite[Lemma 5.10]{NR78}}]\label{prop:kodaspenhecke}
  The Kodaira-Spencer map\index{Kodaira-Spencer map} $T_{L,p'} \to H^1\bigl(C,\ad(E)\bigr)$
  coincides, up to sign, with the composition of the following two
  morphisms.
  \begin{itemize}
  \item The natural map $\alpha:T_{L,p'}\rightarrow E_x$, and
  \item the connecting morphism $\beta: E_x\rightarrow H^1(E\otimes
    E^\vee)$ of the exact sequence~\eqref{eq:elemtransdual} tensored
    with $E$. \qed
  \end{itemize}
\end{proposition}

For the other description of $\tau_E$, we mention another result of
Hwang\index{Hwang, Jun-Muk}.

\begin{theorem}[\protect{\cite[thms.~3-4]{HW01}, \cite[prop.~11]{HW00}}]\label{thm:Hwangveryample}
  With the same notation as above, the morphism $\tau_E$ is given by
  the complete linear system associated to the line bundle
  $\rho_*(K_C)\otimes\O(1)$. If moreover $g>2r+1$ then $\tau_E$ is an
  embedding. \qed
\end{theorem}

By the minimality of Hecke curves\index{Hecke curves!have minimal degree}, Theorem~\ref{thm:tgtmor} implies
that $\tau_E$ is a finite map, that is, that the line bundle
$\rho_*(K_C)\otimes\O(1)$ is ample.

\section{Applications}\label{sec:applimoduli}

The description of the VMRT\index{VMRT} given in Section~\ref{sec:VMRTmoduli} has
interesting consequences.  For instance, Hwang\index{Hwang, Jun-Muk} has used the projective
properties of the tangent morphism $\tau_E$ to deduce the following
 properties of the tangent bundle to the moduli space.

\begin{theorem}[\protect{\cite[thm.~1]{HW00}, \cite[cor.~3]{HW01}}]\label{thm:Hwangstable}\index{Moduli space of vector bundles on a curve!stability of the tangent bundle}
  Let $M^2$ be the moduli space of stable bundles of rank $2$ with a
  fixed determinant of odd degree over an algebraic curve of genus
  $g\geq 2$. Then the tangent bundle of $M^2$ is stable. Let $M^r$ be
  the moduli space of semistable bundles of rank $r$ with fixed
  determinant and $(M^r)^0$ its smooth locus. Then $T_{(M^r)^0}$ is
  simple for $g\geq 4$. \qed
\end{theorem}

The above description of the variety of minimal rational tangents\index{Variety of minimal rational tangents} also
allows to deduce some of its projective properties, for instance its
secant defect. This has been done in the case $r=2$.

\begin{proposition}[\protect{\cite[prop.~6.10]{HK05}}]\label{thm:HwangKebekus05}\index{Moduli space of vector bundles on a curve!variety of minimal rational tangents!secant defect}
  Let $M^2$ be the moduli space of stable bundles of rank $2$ with a
  fixed determinant of odd degree over an algebraic curve of genus
  $g\geq 4$. The variety of minimal rational tangents at a general
  point of $M^2$ has no secant defect. \qed
\end{proposition}

This result, combined with Proposition~\ref{prop:multbound} from page
\pageref{prop:multbound}, provides a bound on the multiplicity of
divisors in the moduli space, perhaps similar in spirit to the
classical Riemann singularity theorem\index{Riemann singularity theorem}.

\begin{corollary}[\protect{\cite[cor.~6.12]{HK05}}]\label{cor:HwangKebekus05}\index{Nadel, Alan M.!product theorem!generalization}
  With the same notation as above, let $x\in M^2$ be a general point,
  and $L$ be the ample generator of $\Pic(M^2)$, and $D\in|mL|$,
  $m\geq 1$ be any divisor. Then
  $$
  mult_x(D)\geq 2m(g-1).
  $$
  \qed
\end{corollary}

The fact that the only minimal rational curves at the general point
are Hecke curves\index{Hecke curves!have minimal degree} has also very important corollaries, as pointed out
by Sun\index{Sun, Xiaotao}. To begin with, it allows us to state a Torelli-type theorem
for moduli spaces:

\begin{theorem}[\protect{\cite[cor.~1.3]{Sun05}}]\label{thm:SunCor1}\index{Moduli space of vector bundles on a curve!Torelli-type theorem}
  Let $C$ and $C'$ be two smooth projective curves of genus $g\geq 4$. Let
  $M^r$ and $(M')^r$ be two irreducible components of the moduli
  schemes of vector bundles of rank $r$ over $C$ and $C'$,
  respectively. If $M^r\cong (M')^r$ then $C\cong C'$. \qed
\end{theorem}

Second, it can be used to describe the automorphism group of $M^r$.

\begin{theorem}[\protect{\cite[cor.~1.4]{Sun05}}]\label{SunCor2}\index{Moduli space of vector bundles on a curve!automorphism group}
  Let $C$ be a smooth projective curve of genus $g\geq 4$, $L$ a line
  bundle on $C$. If $r>2$, then the group of automorphisms of $M^r(L)$
  is generated by:
  \begin{itemize}
  \item automorphisms induced by automorphisms of $C$, and
  \item automorphisms of the form $E\mapsto E\otimes L'$, where $L'$
    is an $r$-torsion element of $\Pic^0(C)$.
  \end{itemize}
  For $r=2$, we need additional generators of the form $E\mapsto
  E^\vee\otimes L'$ where $L'$ is a line bundle verifying
  $(L')^{\otimes 2} \cong L^{\otimes 2}$. \qed
\end{theorem}
\index{Moduli space of vector bundles on a curve|)}
% Local Variables:
% TeX-master: "RC-arXive.tex"
% End:

\part{Consequences of non-uniruledness}

\chapter{Deformations of surjective morphisms}\index{Deformation space of a surjective morphism|(}
\label{chap:hkp}

Let $f: X \to Y$ be a surjective morphism between normal
complex projective varieties. A classical problem of complex geometry
asks for a criterion to guarantee the (non-)existence of deformations
of the morphism $f$, with $X$ and $Y$ fixed. More generally, one is
interested in a description of the connected component $\Hom_f (X,Y)
\subset \Hom (X,Y)$ of the space of morphisms.\index{Deformation space of a surjective morphism}

For instance, if $X$ is of general type, it is well-known that the
automorphism group is finite. It is more generally true that
surjective morphisms between projective manifolds $X$, $Y$ of general
type are always infinitesimally rigid so that the associated connected
components of $\Hom(X,Y)$ are reduced points. Similar questions were
discussed in the complex-analytic setup by Borel\index{Borel, Armand} and Narasimhan\index{Narasimhan, Raghavan},
\cite{BN67}. We will show here how Miyaoka's\index{Miyaoka, Yoichi} characterization of
uniruledness\index{Miyaoka, Yoichi!criterion of uniruledness}, or the existence result for rational curves contained in
Theorem~\ref{thm:BMcQ}/Corollary~\ref{cor:Miyaoka} can be used to give
a rather satisfactory answer in the projective case. Before giving an
idea of the methods employed, we will first, in
Sections~\ref{sec:HKP1-result} and \ref{sec:HKP2-result}, state and
explain the result.

\section{Description of $\Hom_f(X,Y)$ if $Y$ is not uniruled}\index{Deformation space of a surjective morphism!when target is not uniruled|(}
\label{sec:HKP1-result}

If $f : X \to Y$ is as above, it is obvious that $f$ can always be
deformed if the target variety $Y$ has positive-dimensional
automorphism group; this is because the composition morphism
\begin{equation}
  \label{eq:authom}
  \begin{array}{rccc}
    f^\circ : & \Aut^0(Y) & \to & \Hom_f(X,Y) \\
    & g & \mapsto & g \circ f
  \end{array}
\end{equation}\index{Deformation space of a surjective morphism}
is clearly injective. One could na\"ively hope that all deformations
of $f$ come from automorphisms of the target. While this hope does not
hold true in general, we will show, however, that it is \emph{almost}
true: if $Y$ is not uniruled, $f$ always factors via an intermediate
variety $Z$ whose automorphism group is positive-dimensional and
induces all deformations of $f$. More precisely, the following holds.

\begin{theorem}[\protect{\cite[thm.~1.2]{hkp03}}]\label{thm:hkp-main}\index{Deformation space of a surjective morphism!when target is not uniruled}
  Let $f: X \to Y$ be a surjective morphism between normal
  complex-projective varieties, and assume that $Y$ is not uniruled.
  Then there exists a factorization of $f$,
  $$
  \xymatrix{ X \ar[r]_{\alpha} \ar@/^0.5cm/[rr]^{f} & Z
    \ar[r]_{\beta} & Y},
  $$
  such that:
  \begin{enumerate}
  \item the morphism $\beta$ is unbranched away from the singularities
    of $X$ and $Y$, and
    
  \item the natural morphism
    $$
    \begin{array}{ccc}
      \factor \Aut^0 (Z) . {\rm Deck} (Z/Y) . & \to & \Hom_f (X,Y)\\
      g & \mapsto & \beta \circ g \circ \alpha
    \end{array}
    $$
    is an isomorphism of schemes, where ${\rm Deck} (Z/Y)$ is the
    group of Deck transformations, i.e., relative automorphisms.
  \end{enumerate}
  In particular, $f$ deforms unobstructedly, and the associated
  component $\Hom_f(X,Y)$ is a smooth abelian variety.
\end{theorem}

The following is an immediate corollary whose proof we omit for
brevity.

\begin{corollary}[\protect{\cite[cor.~1.3]{hkp03}}]
  In the setup of Theorem~\ref{thm:hkp-main}, if the target variety
  $Y$ is smooth, then $Y$ admits a finite, \'etale covering of the form
  $T \times W$, where $T$ is a torus of dimension $\dim T = h^0 (X,
  f^*(T_Y))$. Additionally, we have
  $$
  \dim \Hom_f (X,Y) \leq \dim Y - \kappa (Y),
  $$
  where $\kappa(Y)$ is the Kodaira dimension. \qed\index{Deformation space of a surjective morphism!dimension of}
\end{corollary}

We will later, in section~\ref{sec:HKP-idea} give an idea of the proof
of Theorem~\ref{thm:hkp-main}.
\index{Deformation space of a surjective morphism!when target is not uniruled|)}

\section{Description of $\Hom_f(X,Y)$ in the general case}\index{Deformation space of a surjective morphism!description of|(}
\label{sec:HKP2-result}

If $Y$ is rationally connected, partial descriptions of the
$\Hom$-scheme are known ---the results of \cite[thm.~1]{HM03} and
\cite[thm.~3]{HM04} assert that whenever $Y$ is a Fano manifold of
Picard number 1 whose variety of minimal rational tangents\index{Variety of minimal rational tangents} is finite,
or not linear, then all deformations of $f$ come from automorphisms of
$Y$. This covers all examples of Fano manifolds of Picard number one
that we have encountered in practice.

If $Y$ is covered by rational curves, but not rationally connected, we
consider the rationally connected quotient $q_Y : Y \dasharrow Q_Y$
which is explained in Definition~\ref{def:ratQuot}. It is shown in
\cite{KP05} that $f$ can be factored via an intermediate variety $Z$,
in a manner similar to that of Theorem~\ref{thm:hkp-main}, such that a
covering of $\Hom_f(X,Y)$ decomposes into
\begin{itemize}
\item an abelian variety, which comes from the automorphism group of
  $Z$, and
  
\item the space of deformations that are relative over the rationally
  connected quotient\index{Rationally connected quotient}, i.e., $H^f_{\vertical} := \{ f' \in
  \Hom_f(X,Y)_{\red} \,|\, q_Y \circ f' = q_Y \circ f \}$.\index{Deformation space of a surjective morphism!relative over the rationally connected quotient}
\end{itemize}

To formulate the result precisely, we recall a result that yields a
factorization of $f$ and may be of independent interest.

\begin{theorem}[\protect{\cite[thm.~1.4]{KP05}}]\label{thm:KP-main0}\index{Maximally \'etale factorization}
  Let $f : X \to Y$ be a surjective morphism between normal projective
  varieties. Then there exists a factorization
  \begin{equation}
    \label{eq:fact_of_f}
    \xymatrix{ {X} \ar[r]_{\alpha} \ar@/^0.3cm/[rr]^{f} & {Z}
      \ar[r]_{\beta} & {Y}}
  \end{equation}
  where $\beta$ is finite and \'etale in codimension
  one\footnote{I.e., where $\beta$ is finite and \'etale away from a
    set of codimension two.}, and where the following universal
  property\index{Maximally \'etale factorization!universal property} holds: for any factorization $f = \beta' \circ \alpha'$,
  where $\beta': Z' \to Y$ is finite and \'etale in codimension 1,
  there exists a morphism $\gamma: Z \to Z'$ such that
  $\beta = \beta' \circ \gamma$. \qed
\end{theorem}

It follows immediately from the universal property that the
factorization~\eqref{eq:fact_of_f} is unique up to isomorphism.  We
call~\eqref{eq:fact_of_f} the \emph{maximally \'etale factorization}\index{Maximally \'etale factorization}
of $f$. The maximally \'etale factorization can be seen as a natural
refinement of the Stein factorization\index{Refinement of Stein factorization}. More precisely, we can say that
a surjection $f: X \to Y$ of normal projective varieties decomposes as
follows.
$$
\xymatrix{ {X} \ar[rr]_{\text{conn.~fibers}}
  \ar@/^0.3cm/[rrrrr]^{f} & & {W} \ar[r]_{\text{finite}} & {Z}
  \ar[rr]_{\text{max.~\'etale}} & & {Y}}
$$

The paper \cite{KP05} discusses the maximally \'etale factorization\index{Maximally \'etale factorization} in
more detail. Its stability under deformations of $f$\index{Maximally \'etale factorization!stability under deformations} is shown
\cite[sect.~1.B]{KP05}, and a characterization\index{Maximally \'etale factorization!characterization} in terms of the
positivity of the push-forward sheaf $f_*(\O_X)$ is given,
\cite[sect.~4]{KP05}.

The main result is then formulated as follows.

\begin{theorem}[\protect{\cite[thm.~1.10]{KP05}}]\label{thm:KP-main}\index{Deformation space of a surjective morphism!description of}
  In the setup of Theorem~\ref{thm:KP-main0}, let $T \subset
  \Aut^0(Z)$ be a maximal compact abelian subgroup. Then there exists
  a normal variety $\tilde H$ and an \'etale morphism
    $$
    T \times \tilde H \to \Hom^n_f(X,Y)
    $$
    that maps $\{e\} \times \tilde H$ to the preimage of
    $H^f_{\vertical}$.  If $Y$ is smooth or if $f$ is itself maximally
    \'etale, then $\{e\} \times \tilde H$ surjects onto the preimage
    of $H^f_{\vertical}$. \qed
\end{theorem}

In Theorem~\ref{thm:KP-main}, we discuss the maximal compact abelian
subgroup of an algebraic group. Recall the following basic fact of
group theory.

\begin{remark}
  Let $G$ be an algebraic group. Then there exists a maximal compact
  abelian subgroup, i.e., an abelian variety $T \subset G$ which is a
  subgroup and such that no intermediate subgroup $T \subset S \subset
  G$, $T \not = S$, is an abelian variety. A maximal compact abelian
  subgroup is unique up to conjugation.
\end{remark}

\begin{remark}
  In the setup of Theorem~\ref{thm:KP-main}, it need not be true that
  $\Aut^0(Z)$ is itself an abelian variety.  Unlike
  Theorem~\ref{thm:hkp-main}, Theorem~\ref{thm:KP-main} does not make
  any statement about the scheme structure of $\Hom_f(X,Y)$.
\end{remark}

While the methods used to show Theorems~\ref{thm:hkp-main} and
\ref{thm:KP-main} are obviously related, the fact that the rational
quotient is generally not a morphism, and the lack of a good parameter
space for rational maps makes the proof of Theorem~\ref{thm:KP-main}
technically more involved. We have thus decided to restrict to a
sketch of a proof of Theorem~\ref{thm:hkp-main} only.
\index{Deformation space of a surjective morphism!description of|)}

\section{Sketch of proof of Theorem~4.1}\index{Deformation space of a surjective morphism!when target is not uniruled|)}
\label{sec:HKP-idea}

We give only a rough idea of a proof for Theorem~\ref{thm:hkp-main}
---our main intention is to show how the existence result for rational
curves comes into the picture. The interested reader is referred to
the rather short original article \cite{hkp03}, and perhaps to the
more detailed survey \cite{Kebekus04}.

\subsection*{Step 1: simplifying assumption}

As we are only interested in a presentation of the core of the
argumentation, we will show Theorem~\ref{thm:hkp-main} under the
simplifying assumption that the surjective morphism $f : X \to Y$ is a
finite morphism between complex-projective manifolds.

\subsection*{Step 2: The tangent map to $f^\circ$}

Consider again the natural morphism $f^\circ : \Aut^0 (Y) \to \Hom_f
(X,Y)$, as introduced in Equation~\eqref{eq:authom} on page
\pageref{eq:authom}. The universal properties of $\Hom(X,Y)$ and of
the automorphism group $\Aut^0(Y)$ yield the following description of
the tangent map $Tf^\circ$ at the point $e \in \Aut^0(Y)$,
\begin{equation}
  \label{eq:tanf}
  Tf^\circ|_e : T_{\Aut(Y)}|_e \to T_{\Hom}|_f.    
\end{equation}
Namely, the natural identifications
$$
T_{\Aut(Y)}|_e \cong H^0(Y, T_Y) \text{\quad and \quad} T_{\Hom}|_f
\cong H^0 (X,f^*(T_Y))
$$
associate the tangent morphism~\eqref{eq:tanf} with the pull-back
map
\begin{equation}
  \label{eq:tanfe}
  Tf^\circ|_e = f^* :H^0(Y, T_Y) \to  H^0 (X,f^* T_Y).  
\end{equation}

The following observation is now immediate:

\begin{lemma}
  The morphism~\eqref{eq:tanf} is injective. If the pull-back
  morphism~\eqref{eq:tanfe} is surjective, i.e., if any infinitesimal
  deformation of $f$ comes from a vector field on $Y$,
  then~\eqref{eq:tanf} yields an isomorphism of Zariski tangent
  spaces. \qed
\end{lemma}

Similar considerations also yield a description of the tangent
morphism at an arbitrary point $g \in \Aut^0(Y)$,
$$
Tf^\circ|_g = (g\circ f)^* :H^0(Y, T_Y) \to H^0 \bigl(X,(g\circ f)^*
(T_Y)\bigr).
$$
Since $g: Y \to Y$ has maximal rank everywhere, it can be shown
that the rank of $Tf^\circ$ is constant on all of $\Aut^0(Y)$. Using
that $\Aut^0(Y)$ is smooth, an elementary argumentation then yields
the following.

\begin{corollary}\label{cor:hkp-notfromdown}
  If the pull-back morphism~\eqref{eq:tanfe} is surjective, then the
  composition map $f^\circ : \Aut^0 (Y) \to \Hom_f (X,Y)$ is
  isomorphic. In this case, the proof of Theorem~\ref{thm:hkp-main} is
  finished by setting $Z := Y$. \qed
\end{corollary}

\subsection*{Step 3, central step in the proof: construction of a covering}

Corollary~\ref{cor:hkp-notfromdown} allows us to assume without loss
of generality that there exists a section $\sigma \in H^0 (X, f^*
(T_Y))$ that does {\it not} come from a vector field on $Y$.  Under
this assumption, we will then construct an unbranched covering of $Y$,
which factors the surjection $f$. The main tool is the following
negativity theorem of Lazarsfeld\index{Lazarsfeld, Robert}\index{Lazarsfeld, Robert!Negativity theorem} for the push-forward sheaf
$f_*({\mathcal O}_X)$ which can, in our setup, be seen as a competing
statement to Corollary~\ref{cor:Miyaoka}\index{Miyaoka, Yoichi!criterion of uniruledness} of
page~\pageref{cor:Miyaoka}, the characterization of uniruledness.

\begin{theorem}[{\cite{Laz80}, \cite[thm.~A]{PS00}}]\label{thm:laz}
  The trace map $tr: f_*(\mathcal O_X) \to \mathcal O_Y$ yields a
  natural splitting
  $$
  f_* ({\mathcal O}_X) \cong \mathcal O_Y \oplus \mathcal E^\vee,
  $$
  where $\mathcal E$ is a vector bundle with the following
  positivity property: if $C \subset Y$ is any curve not contained in
  the branch locus of $f$, then $\mathcal E|_C$ is nef. The curve $C$
  intersects the branch locus of $f$ if and only if $\deg (\mathcal
  E|_C) > 0$. \qed
\end{theorem}

\begin{corollary}\label{cor:hkp-constr}
  Choose an ample line bundle $H \in \Pic(Y)$ and let $C \subset Y$ be
  an associated general complete intersection curve\index{General complete intersection curve} in the sense of
  Mehta-Ramanathan\index{Mehta, Vikram B.}\index{Ramanathan, Annamalai}, as in Definition~\ref{def:GCIC} of
  page~\pageref{def:GCIC}.  Then the following holds:
  \begin{itemize}
  \item If $X \to Z \xrightarrow{\beta} Y$ is a factorization of $f$
    with $\beta$ \'etale, then $\beta_*(\sO_Z) \subset f_*(\sO_Y)$ is
    a subbundle which is closed under the multiplication map
    $$
    \mu : f_* (\mathcal O_X) \otimes f_* (\mathcal O_X) \to f_*
    (\mathcal O_X),
    $$
    and satisfies $\deg \beta_*(\sO_Z)|_C = 0$.
  \item Conversely, if $\sF \subset f_* (\mathcal O_X)$ is a subbundle
    that is closed under multiplication, and $\deg (\sF|_C) = 0$, then
    $f$ factors via $Z := \text{\bf Spec} (\sF)$, and $Z$ is \'etale
    over $Y$. \qed
  \end{itemize}
\end{corollary}

Observe that the projection formula gives an identification
$$
H^0 (X, f^* (T_Y)) = H^0 (Y, f_* (f^* (T_Y))) = H^0 (Y, T_Y) \oplus
\underbrace{H^0 (Y, \mathcal E^\vee \otimes T_Y)}_{= \Hom_Y(\mathcal
  E, T_Y)}
$$
Since $\sigma$ does not come from a vector field on $Y$, the
section $\sigma$ is not contained in the component $H^0 (Y, T_Y)$.
Thus, we obtain a non-trivial morphism $\sigma : \mathcal E \to T_Y$.

Now choose $H$ and $C$ as in Corollary~\ref{cor:hkp-constr}.
Lazarsfeld's Theorem~\ref{thm:laz}\index{Lazarsfeld, Robert!Negativity theorem} then implies that
Image$(\sigma)|_C$ is nef.  On the other hand,
Corollary~\ref{cor:Miyaoka}\index{Miyaoka, Yoichi!criterion of uniruledness} asserts that Image$(\sigma)|_C$ cannot be
ample. In summary we have the following.

\begin{lemma}\label{faz:NA}
  The restricted vector bundle $\mathcal E|_C$ is nef, but not ample.
  \qed
\end{lemma}

\begin{definition}
  Let $\mathcal V_C \subset \mathcal E|_C$ be the maximal ample
  subbundle\index{Maximally ample subbundle}, as discussed in Proposition~\ref{prop:HN1-new} and
  Definition~\ref{def:masb}. Further, let $\mathcal F_C \subset
  \mathcal E^\vee|_C$ be the kernel of the associated morphism
  $\mathcal E^\vee|_C \to \mathcal V_C^\vee$.
\end{definition}

\begin{remark}
  The bundle $\mathcal F_C$ is dual to the quotient $\factor \mathcal
  E|_C.\mathcal V_C.$. In particular, $\mathcal F_C$ is nef and has
  degree 0.
\end{remark}

With these preparations we will now construct the factorization of
$f$. We construct the factorization first over $C$, and then extend
it to all of $Y$.

\begin{lemma}
  The sub-bundle $\mathcal O_C \oplus \mathcal F_C \subset f_*
  (\mathcal O_x) |_C$ is closed under the multiplication map
  $$
  \mu : \bigl(f_* (\mathcal O_X) \otimes f_* (\mathcal O_X)\bigr)|_C \to f_*
  (\mathcal O_X)|_C,
  $$
  because the associated morphism
  $$
  \mu' : \underbrace{(\mathcal O_C \oplus \mathcal F_C) \oplus
    (\mathcal O_C \oplus \mathcal F_C)}_{\deg 0, {\rm nef}}
  \longrightarrow \underbrace{\factor \mathcal O_C \oplus \mathcal E^\vee
    |_C.\mathcal O_C \oplus \mathcal F_C.}_{\cong \mathcal V_C^\vee,
    \text{ anti-ample}}
  $$
  is necessarily trivial. \qed
\end{lemma}

By Corollary~\ref{cor:hkp-constr}, the subbundle $\mathcal O_C \oplus
\mathcal F_C \subset f_* (\mathcal O_x) |_C$ induces a factorization
of the restricted morphism $f|_C : f^{-1}(C) \to C$ via an \'etale
covering of $C$. In order to extend this covering from $C$ to all of $Y$, if
suffices to extend the maximal ample subbundle\index{Maximally ample subbundle} $\mathcal V_C \subset
\mathcal E|_C$ to a subbundle $\mathcal V \subset \mathcal E$ that has
the property that the restriction $\mathcal V|_{C'} \subset \mathcal
E|_{C'}$ to any complete intersection curve $C' \subset Y$ is exactly
the maximal ample subbundle of $\mathcal E|_{C'}$. But since the
maximal ample subbundle is by definition a term of the
Harder-Narasimhan filtration of $\mathcal E|_C$, Theorem~\ref{thm:fmr}
of Flenner\index{Flenner, Hubert} and Mehta-Ramanathan\index{Mehta, Vikram B.}\index{Ramanathan, Annamalai}\index{Mehta-Ramanathan theorem}, exactly asserts that this is possible.

\begin{corollary}\label
  If the pull-back morphism $H^0 (Y, T_Y) \to H^0 (X, f^*(T_Y))$ is
  not surjective, then there exists a factorization of $f$,
  $$
  \xymatrix{ X \ar[r]_{\alpha} \ar@/^0.5cm/[rr]^{f} & Y^{(1)}
    \ar[r]_{\beta} & Y},
  $$
  where $\beta$ is \'etale. \qed
\end{corollary}

\subsection*{Step 4: End of proof}

Assuming that the pull-back map $H^0 (Y, T_Y) \to H^0 (X, f^*(T_Y))$
was not surjective, we have in Step~3 constructed a factorization of
$f$ via an \'etale covering $\beta: Y^{(1)} \to Y$. The \'etalit\'e
obviously implies $f^*(T_Y) = \alpha^*(T_{Y^{(1)}})$. We ask again if
the pull-back morphism
$$
\alpha^* : H^0 \bigl( Y^{(1)}, T_{Y^{(1)}} \bigr) \to H^0 (X, f^*(T_Y))
$$
is surjective.

\begin{description}
\item[\bf Yes $\to$] Following the considerations of Step 2, the proof is
  finished if we set $Z := Y^{(1)}$.
  
\item[\bf No $\to$] repeat Step 3, using the morphism $\alpha: X \to
  Y^{(1)}$ rather than $f:X\to Y$.
\end{description}

Repeating this procedure, we construct a sequence of \'etale coverings
$$
\xymatrix{
  X  \ar[r] \ar@/^0.5cm/[rrrrr]^{f} &
  Y^{(d)} \ar[r] &
  Y^{(d-1)} \ar[r] &
  \ldots \ar[r] &
  Y^{(1)} \ar[r] &
  Y.
}
$$
The sequence, however, must terminate after finitely many steps,
simply because the number of leaves is finite. The same line of
argumentation that led to Corollary~\ref{cor:hkp-notfromdown} shows
that we can end the proof is we set $Z = Y^{(d)}$. This finishes the
proof of Theorem~\ref{thm:hkp-main} ---under the simplifying
assumption that $f$ is a finite morphism between complex manifolds.
\index{Deformation space of a surjective morphism!when target is not uniruled|)}

\section{Open Problems}

It is not quite clear to us if projectivity or if complex number field
is really used in an essential manner. We would therefore like to ask
the following.
\begin{question}
  Does Theorem~\ref{thm:hkp-main} hold in positive characteristic?
\end{question}
\begin{question}
  Does it hold for  K\"ahler manifolds or complex spaces?
\end{question}

Theorem~\ref{thm:hkp-main} can also be interpreted as follows: all
obstructions to deformations of surjective morphisms come from
rational curves in the target. Is it possible to make this statement
precise?
\index{Deformation space of a surjective morphism|)}
% Local Variables:
% TeX-master: "RC-arXive.tex"
% End:
 
\chapter{Families of canonically polarized varieties}\index{Family of canonically polarized varieties|(}
\label{chap:hyperbolicity}

Let $B^\circ$ be a smooth quasi-projective curve, defined over an
algebraically closed field of characteristic $0$ and $q>1$ a positive
integer. In his famous paper \cite{Shaf63}, Shafarevich\index{Shafarevich, Igor R.} considered the
set of families of curves of genus $q$ over $B^\circ$. More precisely,
he considered isomorphism classes of smooth proper morphisms $f: S \to
B^\circ$ whose fibers are connected curves of genus $q$. He
conjectured the following.

\begin{description}
\item[Finiteness conjecture:]\index{Shafarevich, Igor R.!Finiteness conjecture} There are only finitely many isomorphism
  classes of non-isotrivial families of smooth projective curves of
  genus $q$ over $B^\circ$ ---recall that a family is called
  \emph{isotrivial}\index{Family of canonically polarized varieties!isotrivial} if any two fibers are isomorphic.

\item[Hyperbolicity conjecture:]\index{Shafarevich, Igor R.!Hyperbolicity conjecture} If $\kappa(B^\circ) \leq 0$, then no such
  families exist.
\end{description}

These conjectures, which later played an important role in Faltings'
proof of the Mordell conjecture\index{Mordell conjecture}\index{Faltings, Gerd}, were confirmed by Parshin\index{Parshin, Alexei N.}
\cite{Parshin68} for projective bases $B^\circ$ and by Arakelov\index{Arakelov, Sergei J.}
\cite{Arakelov71} in general. We refer the reader to the survey
articles \cite{Viehweg01} and \cite{Kovacs03c} for a historical
overview and references to related results.

It is a natural and important question whether similar statements hold
for families of higher dimensional varieties over higher dimensional
bases.  Families over a curve have been studied by several authors in
recent years and they are now fairly well understood---the strongest
results known were obtained in \cite{Vie-Zuo01, VZ02}, and
\cite{Kovacs02}\index{Kov\'acs, S\'andor J.}.  For higher dimensional bases, however, a complete
picture is still missing and no good understanding of subvarieties of
the corresponding moduli stacks is available.  As a first step toward
a better understanding, Viehweg\index{Viehweg, Eckart} conjectured the following:

\begin{conjecture}[\protect{\cite[6.3]{Viehweg01}}]\label{conj:viehweg}\index{Viehweg, Eckart!conjecture relating variation and log Kodaira dimension}
  Let $f^\circ: X^\circ \to S^\circ$ be a smooth family of canonically
  pol\-arized varieties\index{Family of canonically polarized varieties}. If $f^\circ$ is of maximal variation\index{Family of canonically polarized varieties!of maximal variation}, then
  $S^\circ$ is of log general type ---see
  Definition~\ref{def:log-Kdim} on page \pageref{def:log-Kdim} for the
  logarithmic Kodaira dimension\index{Logarithmic Kodaira dimension} and log general type.
\end{conjecture}

For the reader's convenience, we briefly recall the definition of
\emph{variation} that was introduced by Koll\'ar\index{Koll\'ar, J\'anos} and Viehweg\index{Viehweg, Eckart}.

\begin{definition}
  Let $f^\circ: X^\circ \to S^\circ$ be a family of canonically
  polarized varieties\index{Family of canonically polarized varieties} and $\mu_{f^\circ} : S^\circ \dasharrow
  \mathfrak M$ the induced map to the corresponding moduli scheme. The
  \emph{variation}\index{Family of canonically polarized varieties!variation of} of $f^\circ$ is defined as $\Var(f^\circ) := \dim
  \mu_{f^\circ}(S^\circ)$.
  
  The family $f^\circ$ is called \emph{isotrivial}\index{Family of canonically polarized varieties!isotrivial} if
  $\Var(f^\circ)=0$. It is called \emph{of maximal variation}\index{Family of canonically polarized varieties!of maximal variation} if
  $\Var(f^\circ)=\dim S^\circ$.
\end{definition}

\section{Statement of result}

Using a result of Viehweg-Zuo\index{VZ Theorem, existence of pluri-log differentials on the base of a family}\index{Viehweg, Eckart}\index{Zuo, Kang} and Keel-McKernan's\index{Keel, Se\'an}\index{McKernan, James} \index{Miyanishi, Masayoshi!conjecture on affine lines on quasi-projective varieties}proof of the
Miyanishi\index{Miyanishi, Masayoshi} conjecture in dimension two, we describe families of
canonically polarized varieties over quasi-projective surfaces. We
relate the variation of the family to the logarithmic Kodaira
dimension of the base and give an affirmative answer to Viehweg's\index{Viehweg, Eckart}
Conjecture \ref{conj:viehweg} for families over surfaces.

\begin{theorem}[\protect{\cite[thm.~1.4]{KK05}}]\label{thm:KK-main}\index{Family of canonically polarized varieties!theorem relating variation and log.~Kodaira dimension of the base}
  Let $S^\circ$ be a smooth quasi-projective surface and $f^\circ:
  X^\circ \to S^\circ$ a smooth non-isotrivial family of canonically
  polarized varieties, all defined over $\mathbb C$.  Then the
  following holds.
  \begin{enumerate}
  \item If $\kappa(S^\circ)=-\infty$, then $\Var(f^\circ) \leq 1$.
    
  \item If $\kappa(S^\circ) \geq 0$, then $\Var(f^\circ) \leq
    \kappa({S^\circ})$.
  \end{enumerate}
  In particular, Viehweg's\index{Viehweg, Eckart} Conjecture holds true for families over
  surfaces.
\end{theorem}

\begin{remark}
  Notice that in the case of $\kappa(S^\circ)=-\infty$ one cannot
  expect a stronger statement. For an easy example take any
  non-isotrivial smooth family of canonically polarized varieties over
  a curve $g: Z\to C$, set $X : = Z\times \bP^1$, $S^\circ := C\times
  \bP^1$, and let $f^\circ := g \times \id_{\bP^1}$ be the obvious
  morphism. Then we clearly have $\kappa(S^\circ)=-\infty$ and
  $\Var(f)=1$.
\end{remark}

\section{Open problems}

Theorem~\ref{thm:KK-main} and its proof seem to suggest that the
logarithmic Kodaira dimension of a variety $S^\circ$ gives an upper
bound for the variation of any family of canonically polarized
varieties over $S^\circ$, unless $\kappa(S^\circ)=-\infty$. We would
thus like to propose the following generalization of Viehweg's\index{Viehweg, Eckart}
conjecture\index{Viehweg, Eckart!conjecture relating variation and log Kodaira dimension}.

\begin{conjecture}[\protect{\cite[conj.~1.6]{KK05}}]\index{Viehweg, Eckart!conjecture relating variation and log Kodaira dimension!generalization of}
  Let $f^\circ: X^\circ \to S^\circ$ be a smooth family of canonically
  polarized varieties. Then either $\kappa(S^\circ)=-\infty$ and
  $\Var(f^\circ)<\dim S^\circ$, or $\Var(f^\circ)\leq
  \kappa(S^\circ)$.
\end{conjecture}

\section{Sketch of proof of Theorem~5.3}

Again we give only an incomplete proof of Theorem~\ref{thm:KK-main}.
We restrict ourselves to the case where $\kappa({S^\circ})=0$ and show
only that $\Var(f^\circ) \leq 1$. We will then, at the end of the
present section, give a rough idea of how isotriviality can be concluded.

\subsection{Setup of Notation}

Let $S$ be a compactification of $S^\circ$ as in
Definition~\ref{def:log-Kdim}, and let $D := S \setminus S^\circ$ be
the boundary divisor, which has at worst simple normal crossings. We
assume $\kappa(S^\circ) = \kappa(K_S+D)=0$.
  
A part of the argumentation involves the log minimal model\index{Log minimal model of a surface} of $(S,D)$
---we refer to \cite{KM98} and \cite{Matsuki02} for details on the log
minimal model program for surfaces. If $\kappa(S^\circ) \not =
-\infty$, we denote the birational morphism from $S$ to its
logarithmic minimal model by
$$
\phi : (S,D) \to (S_{\lambda}, D_{\lambda}),
$$
where $D_\lambda$ is the cycle-theoretic image of $D$ in
$S_{\lambda}$. We briefly recall a few important facts of minimal
model theory in dimension two that can be found in the standard
literature, e.g.~\cite{KM98} or \cite{Matsuki02}.

\begin{proposition}\label{prop:MMP}\index{Log minimal model of a surface}
  The minimal model has the following properties:
  \begin{enumerate}
  \item The pair $(S_{\lambda}, D_{\lambda})$ has only log-canonical
    singularities and $S_{\lambda}$ itself has only log-terminal
    singularities. In particular, $S_{\lambda}$ has only quotient
    singularities and is therefore $\Q$-factorial.
    
  \item The log-canonical divisor $K_{S_{\lambda}} + D_{\lambda}$ is
    nef.
    
  \item The log Kodaira dimension remains unchanged, i.e.,
    $\kappa(K_{S_{\lambda}} + D_{\lambda} ) = \kappa(K_S + D )$.
    
  \item Logarithmic Abundance in Dimension 2\index{Logarithmic abundance}: The linear system
    $|n(K_{S_{\lambda}}+D_{\lambda})|$ is basepoint-free for
    sufficiently large and divisible $n \in \mathbb N$. \qed
  \end{enumerate}
\end{proposition}

\subsection{A result of Viehweg and Zuo}\index{VZ Theorem, existence of pluri-log differentials on the base of a family|(}

Before starting the proof we recall an important result of Viehweg\index{Viehweg, Eckart} and
Zuo\index{Zuo, Kang} that describes the sheaf of logarithmic differentials on the base
of a family of canonically polarized varieties in our setup.

\begin{theorem}[\protect{\cite[thm.~1.4(i)]{VZ02}}]\label{thm:VZ}
  In the setup introduced above, there exists an integer $n > 0$ and
  an invertible subsheaf $\sA \subset \Sym^n \Omega^1_S(\log D)$ of
  Kodaira dimension $\kappa(\sA) \geq \Var(f^\circ)$. \qed
\end{theorem}
\index{VZ Theorem, existence of pluri-log differentials on the base of a family|)}

\subsection{Reduction to the uniruled case} 

A surface $S$ with logarithmic Kodaira dimension zero need of course
not be uniruled. Using the result of Viehweg-Zuo\index{Viehweg, Eckart}\index{Zuo, Kang} we can show,
however, that any family of canonically polarized varieties over a
non-uniruled surface with Kodaira dimension zero is isotrivial.

\begin{proposition}\label{prop:reduction-to-uniruled}
  If $S$ is not uniruled, then $\Var(f^\circ)=0$.
\end{proposition}

We prove Proposition~\ref{prop:reduction-to-uniruled} using three
lemmas that describe the sheaves of log-differentials on the minimal
model $S_\lambda$.

\begin{lemma}\label{lem:ntrivksmin}
  If $n \in \mathbb N$ is sufficiently large and divisible, then
  \begin{equation}
    \label{eq:logmin}
    \O_{S_{\lambda}}(n(K_{S_{\lambda}}+D_{\lambda})) = \O_{S_{\lambda}}.  
  \end{equation}
  In particular, the log-canonical $\mathbb Q$-divisor
  $K_{S_{\lambda}}+D_{\lambda}$ is numerically trivial.
\end{lemma}
\begin{proof}
  Equation~\eqref{eq:logmin} is an immediate consequence of the
  assumption $\kappa(S^\circ)=0$ and the logarithmic abundance theorem
  in dimension 2\index{Logarithmic abundance}, which asserts that the linear system
  $|n(K_{S_{\lambda}}+D_{\lambda})|$ is basepoint-free. 
\end{proof}

\begin{lemma}\label{lem:kntriv}
  If $\kappa(S) \geq 0$, then $S_{\lambda}$ is $\mathbb
  Q$-Gorenstein, $K_{S_{\lambda}}$ is numerically trivial and
  $D_{\lambda} = \emptyset$.
\end{lemma}
\begin{proof}
  Lemma~\ref{lem:ntrivksmin} together with the assumption that
  $|nK_S|\neq\emptyset$ for large $n$ imply that $\phi$ contracts all
  irreducible components of $D$, and all divisors in any linear system
  $|nK_S|$, for all $n \in \bN$. The claim follows. 
\end{proof}

\begin{lemma}\label{lem:unstability}
  Assume that $\kappa(S) \geq 0$ and $\Var(f^\circ) \geq 1$. If $H \in
  \Pic(S_{\lambda})$ is any ample line bundle, then the (reflexive)
  sheaf of differentials $(\Omega^1_{S_{\lambda}})^{\vee \vee}$ has
  slope $\mu_H\bigl((\Omega^1_{S_{\lambda}})^{\vee \vee}\bigr) = 0$,
  but it is not semistable with respect to $H$.
\end{lemma}
\begin{proof}
  Fix a sufficiently large number $m > 0$ and a general curve
  $C_{\lambda} \in |m H|$.  Theorem \ref{thm:fmr} 
  ensures that if
  $(\Omega^1_{S_{\lambda}})^{\vee \vee}$ is semistable, then so is its
  restriction $\Omega^1_{S_{\lambda}}|_{C_{\lambda}}$.
  
  By general choice, $C_{\lambda}$ is contained in the smooth locus of
  $S_{\lambda}$ and stays off the fundamental points of $\phi^{-1}$.
  The birational morphism $\phi$ will thus be well-defined and
  isomorphic along $C := \phi^{-1}(C_{\lambda})$.
  Lemma~\ref{lem:kntriv} then asserts that
  $$
  \mu_H \bigl((\Omega^1_{S_{\lambda}})^{\vee \vee}\bigr) =
  \frac{K_{S_{\lambda}} \cdot C_{\lambda}}{2} = 0,
  $$
  which shows the first claim.
  
  Similarly, Lemma~\ref{lem:kntriv} implies that $\codim_{S_{\lambda}}
  \phi(D) \geq 2$, and so $C$ is disjoint from $D$. The unstability of
  $(\Omega^1_{S_{\lambda}})^{\vee \vee}$ can therefore be checked
  using the identifications
  \begin{equation}\label{eq:identlog}
    (\Omega^1_{S_{\lambda}})^{\vee \vee}|_{C_{\lambda}} \cong 
    \left.\Omega^1_{S_{\lambda}}\right|_{C_{\lambda}} \cong
    \left.\Omega^1_S\right|_C \cong \Omega^1_S(\log D)|_C.    
  \end{equation}
  Since symmetric powers of semistable vector bundles over curves are
  again semistable \cite[cor.~3.2.10]{HL97}, in order to prove
  Lemma~\ref{lem:unstability}, it suffices to show that there exists a
  number $n\in \mathbb N$ such that $\Sym^n \Omega^1_S(\log D)|_C$ is
  not semistable. For that, use the
  identifications~\eqref{eq:identlog} to compute
  \begin{align*}
    \deg_C \Sym^n \Omega^1_S(\log D)|_C & = const^+ \cdot \deg_C
    \left.\Omega^1_S\right|_C \\
    & = const^+ \cdot \deg_{C_{\lambda}}
    (\Omega^1_{S_{\lambda}})^{\vee
      \vee}|_{C_{\lambda}} && \text{Isomorphisms~\eqref{eq:identlog}}\\
    & = const^+ \cdot (K_{S_{\lambda}}\cdot C_{\lambda}) = 0. &&
    \text{Lemma~\ref{lem:kntriv}}
  \end{align*}
  Hence, to prove unstability it suffices to show that $\Sym^n
  \Omega^1_S(\log D)|_C$ contains a subsheaf of positive degree.
  
  Theorem~\ref{thm:VZ}\index{VZ Theorem, existence of pluri-log differentials on the base of a family} implies that there exists an integer $n > 0$
  such that $\Sym^n \Omega^1_S(\log D)$ contains an invertible
  subsheaf $\mathcal A$ of Kodaira dimension $\kappa(\mathcal A) \geq
  1$. But by general choice of $C_\lambda$, this in turn implies that
  $\deg_C(\mathcal A|_C) > 0$, which shows the required unstability.
  This ends the proof of Lemma~\ref{lem:unstability}. 
\end{proof} 

With these preparations, the proof of
Proposition~\ref{prop:reduction-to-uniruled} is now quite short.

\begin{proof}[Proof of Proposition~\ref{prop:reduction-to-uniruled}]
  We argue by contradiction and assume to the contrary that both $S$
  is not uniruled, i.e., $\kappa(S) \geq 0$, and $\Var(f^\circ) \geq
  1$. Again, let $H \in \Pic(S_{\lambda})$ be any ample line bundle.
  
  Lemma~\ref{lem:unstability} implies that
  $\Omega^1_{S_{\lambda}}|_{C_{\lambda}}$ has a subsheaf of positive
  degree or, equivalently, that it has a quotient of negative degree.
  On the other hand, Corollary~\ref{cor:Miyaoka}\index{Miyaoka, Yoichi!criterion of uniruledness} then asserts that $S$
  is uniruled, leading to a contradiction. This ends the proof of
  Proposition~\ref{prop:reduction-to-uniruled}. 
\end{proof}

\subsection{Images of $\C^*$ on $S^\circ$, end of proof}

In view of Proposition~\ref{prop:reduction-to-uniruled}, to prove that
$\Var(f^\circ) \leq 1$, we can assume without loss of generality that
$\Var(f^\circ)>0$ and therefore $S$ is uniruled. Since families of
canonically polarized varieties over $\C^*$ are always isotrivial,
\cite[thm.~0.2]{Kovacs00}\index{Family of canonically polarized varieties!over $\C^*$ are isotrivial}, the result then follows from the following
proposition.

\begin{proposition}\label{prop:uniLK0} 
  If $\Var(f^\circ)>0$, then $S^\circ$ is dominated by images of $\C^*$.
  In particular $\Var(f^\circ) \leq 1$.
\end{proposition}

As a first step in the proof of Proposition, recall the following
fact. If $X$ is a projective manifold, and $H \subset {\rm
  RatCurves}^n(X)$ is an irreducible component of the space of
rational curves such that the associated curves dominate $X$, it is
well understood that a general point of $H$ corresponds to a free
curve\index{Rational curve!free}, i.e., a rational curve whose deformations are not obstructed
(cf.~\cite[1.1]{KMM92}, \cite[II Thm.~3.11]{K96}, see also
Prop.~\ref{prop:gen=stand}). In particular, if $E \subset X$ is an
algebraic set of codimension $\codim_X E \geq 2$, then the subset of
curves that avoid $E$,
$$
H' := \{ \ell \in H \,|\, E\cap \ell = \emptyset \},
$$
is Zariski-open, not empty, and curves associated with $H'$ still
dominate $X$. A similar statement, which we quote as a fact without
giving a proof, also holds if $X$ is a surface with mild
singularities, and for quasi-projective varieties that are dominated
by $\C^*$ rather than complete rational curves.

\begin{proposition}[Small Set Avoidance, \protect{\cite[prop.~2.7]{KK05}}]\label{prop:smallsetavoidance}\index{Small set avoidance theorem}
  Let $X$ be a smooth projective surface, and $E \subset X$ a divisor
  with simple normal crossings.  Assume that $X \setminus E$ is
  dominated by images of $\C^*$, and let $F \subset X\setminus E$ be
  any finite set. Then $X \setminus (E \cup F)$ is also dominated by
  images of $\C^*$.  \qed

\end{proposition}

\begin{lemma}\label{lem:domination}
  In the setup of Proposition~\ref{prop:uniLK0}, the quasi-projective
  surface $S_{\lambda} \setminus (\Sing(S_{\lambda}) \cup D_\lambda)$
  is dominated by images of $\C^*$.
\end{lemma}
\begin{proof}
  We aim to apply Theorem~\ref{thm:KMcK2}.(2), and so we need to show
  that
  \begin{itemize}
  \item the log-canonical divisor $K_{S_{\lambda}}+D_{\lambda}$ is
    numerically trivial, and that
  \item the boundary divisor $D_{\lambda}$ is not empty.
  \end{itemize}
  The numerical triviality of $K_{S_{\lambda}}+D_{\lambda}$ has been
  shown in Lemma~\ref{lem:ntrivksmin} above. To show that $D_{\lambda}
  \not = \emptyset$, we argue by contradiction, and assume that
  $D_\lambda = \emptyset$. Set
  $$
  S_\lambda^\circ := S_\lambda \setminus
  \underbrace{\phi(\text{exceptional set of $\phi$})}_{\text{finite,
      contains $\phi(D)$}}.
  $$
  Then $S_\lambda^\circ$ is the complement of a finite set and
  $\phi^{-1}|_{S_\lambda^\circ}$ is a well-defined open immersion.  Let
  $f_\lambda:=\phi\circ f$. Then $X:=X^\circ|_{f_\lambda^{-1}(S^\circ_\lambda)} \to
  S^\circ_\lambda$ is a smooth family of canonically polarized varieties.
  Consider the following diagram:
  $$
  \xymatrix{ X \ar[d]_{f_\lambda} & & & **[r] \tilde X := X
    \times_{S_\lambda} \tilde S
    \ar[lll] \ar[d]^{\tilde f} \\
    S_\lambda & & \tilde S_\lambda \ar[ll]_{\alpha}^{\text{\ 
        index-one-cover}} & \tilde S \ar[l]_{\beta}^{\text{\ 
        log~resolution}} }
  $$
  where $\alpha$ is the index-one-cover\index{Index-one-cover} described in
  \cite[5.19]{KM98} or \cite[sect.~3.5]{Reid87}, and $\beta$ is the
  minimal desingularization of $\tilde S_\lambda$ composed with
  blow-ups of smooth points so that $\beta^{-1}(\tilde
  S_\lambda\setminus \alpha^{-1}(S^\circ_\lambda))$ is a divisor with at
  most simple normal crossings.
      
  By Lemma~\ref{lem:ntrivksmin}, $K_{S_\lambda}$ is torsion.  Since
  $\alpha$ is \'etale in codimension one this implies that $K_{\tilde
    S_\lambda}$ is trivial. Furthermore, $\tilde S_\lambda$ has only
  canonical singularities: we have already noted in
  Proposition~\ref{prop:MMP} that the singularities of $S_\lambda$ are
  log-terminal, i.e., they have minimal discrepancy $> -1$. Then by
  \cite[Prop.~5.20]{KM98} the minimal discrepancy of the singularities
  of $\tilde S_\lambda$ is also $> -1$, and as $K_{\tilde S_\lambda}$
  is Cartier, the discrepancies actually must be integral and hence
  $\geq 0$, cf.~\cite[proof of Cor.~5.21]{KM98}. Consequently,
  \begin{equation}\label{star}
    K_{\tilde S} =
    \underbrace{\beta^*(K_{\tilde S_\lambda})}_{\cong \mathcal O_{\tilde S}} +
    (\text{effective and $\beta$-exceptional}).        
  \end{equation}
  This in turn has two further consequences:
      \begin{enumerate}
      \item[i)] $\kappa(K_{\tilde S}) = 0$. In particular, $\tilde S$
        is not uniruled.
        
      \item[ii)] If we set $\tilde S^\circ := (\alpha \circ
        \beta)^{-1}(S^\circ_\lambda)$ and $\tilde X^\circ:=\tilde
        f^{-1}(\tilde S^\circ)$ then $\tilde X^\circ\to \tilde
        S^\circ$ is again a smooth family of canonically polarized
        varieties.  Letting $\tilde D := \tilde S \setminus \tilde
        S^\circ$ then $\tilde D$ is exactly the $\beta$-exceptional
        set, and \eqref{star} implies that
        $$
        \kappa(\tilde S^\circ) = \kappa (\underbrace{K_{\tilde S} +
          \tilde D}_{\hskip -1cm \text{effective, $\beta$-exceptional}
          \hskip -1cm}) = 0.
        $$
        In particular, Proposition~\ref{prop:reduction-to-uniruled}
        applies to $\tilde f: \tilde X \to \tilde S$ and shows that
        $Var(\tilde f|_{\tilde X^\circ})=0$.
      \end{enumerate}
      This is a contradiction and thus ends the proof of
      Lemma~\ref{lem:domination}. 
\end{proof}

Observe that Lemma~\ref{lem:domination} does \emph{not} immediately
imply Proposition~\ref{prop:uniLK0}. The problem is that the boundary
divisor $D \subset S$ can contain connected components that are
contracted by $\phi$ to points. These points do not appear in the
cycle-theoretic image divisor $D_\lambda$, and it is a priori possible
that all morphisms $\C^* \to S_{\lambda} \setminus \Sing(S_{\lambda})$
contain these points in the image. Taking the strict transforms would
then give morphisms $\C^* \to S \setminus \phi^{-1}(\Sing(S_{\lambda})
\cup D_\lambda)$, but $\phi^{-1}(\Sing(S_{\lambda}) \cup D_\lambda)
\ne D$. An application of Proposition~\ref{prop:smallsetavoidance}
will solve this problem.

\begin{proof}[Proof of Proposition~\ref{prop:uniLK0}]
  If $\phi(D) \subset D_{\lambda} \cup \Sing(S_{\lambda})$, i.e., if
  all connected components of $D$ are either mapped to singular
  points, or to divisors, Lemma~\ref{lem:domination} immediately
  implies Proposition~\ref{prop:uniLK0}. Likewise, if $S_{\lambda}$
  was smooth, Proposition~\ref{prop:smallsetavoidance} on small set
  avoidance\index{Small set avoidance theorem} would imply that almost all curves in the family stay off
  the isolated zero-dimensional components of $\phi(D)$, and
  Proposition~\ref{prop:uniLK0} would again hold. In the general case,
  where $S_{\lambda}$ is singular, and $d_1, \ldots, d_n$ are smooth
  points of $S_{\lambda}$ that appear as connected components of
  $\phi(D)$, a little more argumentation is required.
  
  If $D'$ is the union of connected components of $D$ which are
  contracted to the set of points $\{d_1, \ldots, d_n\} \subset
  S_{\lambda}$, it is clear that the birational morphism $\phi: S \to
  S_{\lambda}$ factors via the contraction of $D'$, i.e., there exists
  a diagram
  $$
  \xymatrix{ {S} \ar[r]_{\alpha} \ar@/^0.3cm/[rr]^{\phi} & {S'}
    \ar[r]_{\beta} & {S_{\lambda}}}
  $$
  where $S'$ is smooth, and $\alpha$ maps the connected components
  of $D'$ to points $d'_1, \ldots d'_2 \in S'$ and is isomorphic
  outside of $D'$.
  
  Now, if $D'' := D \setminus D'$, the above argument shows that $S'$
  is dominated by rational curves that intersect $\alpha(D'')$ in two
  points. Since $S'$ is smooth,
  Proposition~\ref{prop:smallsetavoidance}\index{Small set avoidance theorem} applies and shows that
  almost all of these curves do not intersect any of the $d'_i$. In
  summary, we have seen that most of the curves in question intersect
  $\alpha(D)$ in two points. This completes the proof of
  Proposition~\ref{prop:uniLK0}.
\end{proof}

\subsection{Sketch of further argument}

We have seen above that $S^\circ$ is dominated by images of $\C^*$,
which implies $\Var(f^\circ) \leq 1$. If $S^\circ$ is \emph{connected
  by $\C^*$}, i.e., if there exists an open subset $\Omega \subset
S^\circ$ such that any two points $x,y \in \Omega$ can be joined by a
chain of $\C^*$, then it is clear that $\Var(f^\circ) = 0$, and
Theorem~\ref{thm:KK-main} is shown in case $\kappa(S^\circ)=0$.

Since $\dim S^\circ = 2$, we can thus assume without loss of
generality that a general point of $S^\circ$ is contained in exactly
one curve which is an image of $\C^*$. Blowing up points on the
boundary $D$ if necessary, we can thus assume that $S$ is a
birationally ruled surface, with a map $\pi : S \to C$ to a curve $C$,
and that the boundary $D$ is a divisor that intersects the general
$\pi$-fiber $F$ in exactly two points, which implies that the
restriction of the vector bundle $\Omega^1_S(\log D)$ to $F$ is
trivial. If $\Var(f^\circ) > 0$, the result of Viehweg-Zuo\index{Viehweg, Eckart}\index{Zuo, Kang}\index{VZ Theorem, existence of pluri-log differentials on the base of a family},
Theorem~\ref{thm:VZ}, then implies that the restriction of
$\Omega^1_S(\log D)$ to non-fiber components of $D$ cannot be stable
of degree zero. However, a detailed analysis of the self-intersection
graph of $D$, and the standard description of the restriction of
$\Omega^1_S(\log D)$ to components of $D$ shows both stability and
zero-degree.
\index{Family of canonically polarized varieties|)}
% Local Variables:
% TeX-master: "RC-arXive.tex"
% End:

\appendix

%\listoffigures

\bibliographystyle{alpha}
\bibliography{../bibliography/general}

%\seealso{Rational connectivity}{Maximally rational connected fibration}
\index{VMRT|see{variety of minimal rational tangents}}
\index{Graber, Tom|see{GHS Theorem}}
\index{Harris, Joe|see{GHS Theorem}}
\index{Starr, Jason|see{GHS Theorem}}
\index{Keel, Se\'an|see{KMcK}}
\index{McKernan, James|see{KMcK}}
\index{$M^r$|see{Moduli space of vector bundles on a curve}}
\index{$\Hom_f(X,Y)$|see{Deformation space of a surjective morphism}}
\index{$\rat(X)$|see{Space of rational curves on a variety}}
\index{Rational curve|see{Family of rational curves, Space of rational curves on a variety}}
\index{Boucksom, S\'ebastien|see{BDPP}}
\index{Demailly, Jean-Pierre|see{BDPP}}
\index{Paun, Mihai@P\u aun, Mihai|see{BDPP}}
\index{Peternell, Thomas|see{BDPP}}
\index{Mori, Shigefumi!Bend-and-Break argument|see{Bend-and-Break}}
\index{Deformation of a surjective morphism|see{Deformation space of a surjective morphism}}
\index{Kodaira dimension of a quasi-projective variety|see{Logarithmic Kodaira dimension}}
\index{Cartan, \'Elie|seealso{Cartan-Fubini type theorem}}
\index{Fubini, Guido|seealso{Cartan-Fubini type theorem}}
\index{Family of canonically polarized varieties|seealso{Shafarevich conjecture, Viehweg conjecture}}
\index{Foliation|seealso{Leaf of a foliation}}
\index{Logarithmic Kodaira dimension!of the base of a family|see{Family of canonically polarized varieties, Shafarevich conjecture, Viehweg conjecture}}
\index{Miyaoka, Yoichi!foliations in positive characteristic|see{Foliation}}
%\index{MRC fibration|see{Maximally rationally connected fibration}}
%\index{Maximally rationally connected fibration|see{Rationally connected quotient}}
\index{MRCC fibration|see{Maximally rationally chain connected fibration}}
\index{Harder-Narasimhan filtration!of the tangent bundle|see{Partial rationally connected quotients}}
\index{Rational chain connectivity|seealso{Maximally rationally chain connected fibration}}
\index{Rational connectivity|seealso{Rationally connected quotient}}
\index{Rational connectivity!criterion for the leaf of a foliation|see{Leaf of a foliation}}
\index{Maximally rationally chain connected fibration|seealso{Rationally connected quotient}}
\index{Rational curve!existence conjecture of Mumford|see{Mumford}}
\index{Rational curve!existence result of Miyaoka|see{Miyaoka, criterion of uniruledness}}
\index{Rational curve!of minimal degree!is hardly ever singular|see{Family of rational curves}}
\index{Rationally connected quotient!non-uniruledness|see{GHS Theorem}}
\index{Rationally connected quotient!relative deformations over|see{Deformation space of a surjective morphism}}
\index{Rational quotient|see{Rationally connected quotient}}
\index{Space of rational curves on a variety|seealso{Family of rational curves}}
\index{Family of rational curves|seealso{Space of rational curves on a variety}}
\index{Spanning dimension|see{Family of rational curves}}
\index{Stability of the tangent bundle|see{Tangent bundle}}
\index{Stein factorization|see{Refinement of Stein factorization, Maximally \'etale factorization}}
\index{Tangent bundle!stability and partial rationally connected quotients|see{Partial rationally connected quotients}}
\index{Tangent bundle!of the variety of minimal rational tangents|see{Variety of minimal rational tangents}}
\index{Tangent bundle!of the moduli space of vector bundles on a curve|see{Moduli space of vector bundles on a curve}}
\index{Tangent map!for Hecke curves in the moduli space of vector bundles on a curve|see{Moduli space of vector bundles on a curve}}
\index{Tangent map|seealso{Tangent morphism}}
\index{Tangent morphism|seealso{Tangent map}}
\index{Uniqueness of contact structures|see{Contact manifold}}
\index{Uniruledness criterion of Miyaoka|see{Miyaoka, Yoichi}}
\index{Variation of a family|see{Family of canonically polarized varieties}}
\index{Variety of minimal rational tangents!of the moduli space of vector bundles on a curve|see{Moduli space of vector bundles on a curve}}
\index{Variety of minimal rational tangents!determines the variety|see{Cartan-Fubini type theorem}}
\index{Variety of minimal rational tangents!secant defect!for the Moduli space of vector bundles on a curve|see{Moduli space of vector bundles on a curve}}
\index{Viehweg, Eckart|seealso{VZ Theorem}}
\index{Zuo, Kang|seealso{VZ Theorem}}
\index{Viehweg, Eckart!conjecture relating variation and log Kodaira dimension!proof in dimension two|see{Family of canonically polarized varieties}}
\index{Vertical subsheaf with respect to the rationally connected quotient|see{Rationally connected quotient}}
\index{Mehta, Vikram B.|seealso{Metha-Ramathan theorem}}
\index{Ramanathan, Annamalai|seealso{Metha-Ramathan theorem}}
\printindex

\end{document}